\documentclass[11pt,a4paper]{article}
\usepackage[a4paper,margin=20mm]{geometry}
\usepackage[T1]{fontenc}
\usepackage{mathpazo, courier, microtype}
\usepackage{amsmath,amssymb,amsthm,mathtools,mathrsfs}
\usepackage[numbers,square,sort&compress]{natbib}
\usepackage{enumitem,xcolor,hyperref}
\usepackage{array,tabularx,booktabs,float,etoolbox,needspace}
\allowdisplaybreaks[4]
\hypersetup{
	hypertexnames=false,
	colorlinks,
	citecolor=blue,
	filecolor=blue,
	linkcolor=blue,
	urlcolor=black
}

\numberwithin{equation}{section}
\newtheorem{thm}{Theorem}[section]
\newtheorem{prop}{Proposition}[section]
\newtheorem{lem}{Lemma}[section]
\newtheorem{cor}{Corollary}[section]
\newtheorem{rem}{Remark}[section]

\newcommand{\R}{\mathbb R}
\newcommand{\N}{\mathbb N}
\renewcommand{\S}{\mathbb S}
\newcommand{\Om}{\Omega}

\newcommand{\al}{\alpha}

\newcommand{\dist}{\operatorname{dist}}
\newcommand{\supp}{\operatorname{supp}}

\newcommand{\Id}{\operatorname{Id}}
\newcommand{\Prob}{\mathcal P}
\newcommand{\Xr}{\mathcal X_r}
\newcommand{\norm}[1]{\lVert#1\rVert}
\newcommand{\weak}{\rightharpoonup}
\newcommand{\be}{\begin{equation}}
\newcommand{\ee}{\end{equation}}

\newcommand{\pa}{\partial}

\newcommand{\lam}{\lambda}

\newcommand{\Lam}{\Lambda}
\newcommand{\si}{\sigma}
\newcommand{\Si}{\Sigma}

\newcommand{\ga}{\gamma}

\newcommand{\Ga}{\Gamma}

\newcommand{\var}{\varepsilon}

\renewcommand{\deg}{\operatorname{deg}}

\renewcommand{\d}{\mathrm{d}}

\title{\bfseries Weighted barycenter sets and singular \(N\)-mean field equations with arbitrary exponents}
\author{
Mohameden Ahmedou\thanks{Mathematisches Institut, Justus-Liebig-Universit\"at Gie\ss en,
Arndtstra\ss e 2, D-35392 Gie\ss en, Germany.
E-mail: \href{mailto:Mohameden.Ahmedou@math.uni-giessen.de}{Mohameden.Ahmedou@math.uni-giessen.de}.},\quad ~
Zhengni Hu\thanks{School of Mathematical Sciences, Shanghai Jiao Tong University,
800 Dongchuan Road, Minhang District, Shanghai 200240, China.
E-mail: \href{mailto:zhengni_hu@sjtu.edu.cn}{zhengni\_hu@sjtu.edu.cn}.},\quad~
Heming Wang\thanks{State Key Laboratory of Mathematical Sciences, Academy of Mathematics and Systems Science,
Chinese Academy of Sciences, Beijing 100190, China.
E-mail: \href{mailto:hmwang@amss.ac.cn}{hmwang@amss.ac.cn}.},
\quad~
Lei Zhang\thanks{Department of Mathematics, University of Florida,
Gainesville, FL 32611, USA.
E-mail: \href{mailto:leizhang@ufl.edu}{leizhang@ufl.edu}.}
}

\date{}

\begin{document}
\maketitle
\begin{abstract}

In this paper, we study singular Dirichlet $N$-mean field equations with arbitrary exponents greater than $-1$ on smooth bounded domains. We compute the Leray--Schauder degree for all positive nonresonant parameters. 
In the explicit degree formula, positive singularities contribute differently in odd and even dimensions. The proofs combine weighted barycenter techniques, entropy duality method, and annular decomposition. Weighted barycenters also yield existence results in cases where the degree vanishes.

\end{abstract}



\tableofcontents

\section{Introduction}\label{sec:intro}
Let $\Om\subset\R^N$, $N\ge2$, be a smooth bounded domain, and let $P=\{p_1,\ldots,p_m\}\subset\Om$ consist of distinct points. We consider the following quasilinear mean field equation:
\be \label{maineq}
\left\{\begin{aligned}
   & -\Delta_Nu=\rho\frac{H(x)e^{u(x)}}{\int_\Om H(y)e^{u(y)}\, \d  y},&&x\in\Om,\\ &u=0,&&x\in\pa \Om,
\end{aligned}\right.
\ee
where $\Delta_Nu=\operatorname{div}(|\nabla u|^{N-2}\nabla u)$, $\rho\geq0$ is a parameter, and $H(x)=h(x)\prod_{i=1}^m|x-p_i|^{N\al _i}$ with
\be \label{eq:assumptions}
h\in C^1(\overline\Om),\quad \min_{\overline\Om}h>0,\quad \al _i>-1 \quad(1\leq i\leq m).
\ee
We take the empty product to be one when $m=0$ and set $|x-p_i|^{N\al_i}=1$ also at $p_i$ when $\al_i=0$. The domain $\Om$ may have arbitrary topology. We call \eqref{maineq} the singular $N$-mean field equation.

For $P=\varnothing$ or $\al_i=0$ for all $i$, \eqref{maineq} is the regular $N$-mean field equation. For $N=2$, the classical mean field equation arises in conformal geometry ( \cite{CGY93,CY87,KW74}), statistical mechanics (see, e.g., \cite{CLMP92,CLMP95,CK94,K93}),  turbulent Euler flows (see, e.g., \cite{DPEM10,S67}), and self-dual condensates in Chern--Simons--Higgs models (see, e.g., \cite{CLW04,DPEFM15,DJLW98,EF14,LW10,LY13,NT99}). 

Energy scaling also links the equation to the planar case.
For $u_\lam (x)=u(\lam  x)$, a change of variables yields
\[
 \int_{\lam ^{-1}\Om}|\nabla u_\lam |^N\, \d  x
 =\int_\Om|\nabla u|^N\, \d  x,\quad \lam >0.
\]
Scale invariance and the exponential integrability of $W^{1,N}_0(\Om)$ explain the choice of the $N$-Laplacian in \eqref{maineq}. Since $\al_i>-1$ in \eqref{eq:assumptions}, $H\in L^s(\Om)$ for some $s>1$. Hence the map
$u\mapsto\int_\Om He^u\, \d  x$ is continuously Fr\'echet
differentiable on $W^{1,N}_0(\Om)$, although $H$ may be
unbounded at the points of $P$.
In \eqref{maineq}, the normalized density has total mass $\rho$. As $\rho$ varies, solutions may develop
high peaks and concentrate at finitely many points.

The Leray--Schauder degree remains constant as $\rho$ and $H$ vary continuously, provided all solutions remain uniformly bounded. Moreover, a nonzero degree implies the existence of at least one solution. In the present paper,  we are interested  in  studying how the degree depends on $\rho$, the topology of $\Om$, and the singular exponents, and how concentration affects its value.

For the regular mean field equation on a closed surface, Li \cite{L99} first
proved the dependence of the nonresonant degree on the parameter and the Euler characteristic alone, and proposed computing the degree from its jumps at blow-up parameters. Subsequently, 
Chen--Lin \cite{CL02}  established sharp estimates for solutions with several
concentration points and obtained explicit degree
formulas on closed surfaces and planar Dirichlet domains in 
\cite{CL03}. Later,  Chen--Lin \cite{CL15}  extended the formulas to singular mean field equations with arbitrary exponents greater than $-1$.
Related degree formulas for two-dimensional mean field type equations
can be found, for instance, in \cite{ABL17,AB17,M26,JWY18,AHW26,M08}. 
In particular, Ahmedou--Ben Ayed--Lucia \cite{ABL17}  and Ahmedou and Ben Ayed \cite{AB17} obtained Leray--Schauder degree formulas
at resonant parameters for the regular Dirichlet mean field equation and its singular
counterpart with positive nonintegral singular exponents, respectively, via Bahri's theory
of critical points at infinity \cite{B89}.
For a recent survey, see 
Malchiodi \cite{M26} and the references therein.

The topology of energy sublevels provides another approach to proving existence results. 
Ding--Jost--Li--Wang \cite{DJLW97,DJLW99} used minimization and
min--max arguments to prove existence at and beyond the first
concentration threshold. Djadli--Malchiodi \cite{DM08} and Djadli \cite{D08} applied formal barycenters to describe concentration at several points.
Malchiodi \cite{M08} further related the degree to the Euler characteristic
of low energy sets, using the deformation method of Lucia
\cite{L07}.

For singular equations, Carlotto--Malchiodi \cite{CM12} introduced
different weights to ordinary and singular concentration points.
Weighted barycenters parametrize test functions and provide a target space for the projection from low energy sublevels. The test map followed by the projection is homotopic to the identity on the barycenter space, which suffices for proving existence results. However, a more detailed analysis is needed to determine the topology of low energy sublevels. For weighted barycenter spaces,  Kallel \cite{K19} computed the Euler
characteristic of weighted barycenter spaces, which also related the results to the degree formulas in \cite{CL03,CL15}. 
Similar barycenter and degree methods appear in the boundary problem of  Li--Sun--Yang \cite{LSY23} and the
generalized equation studied by  Sun--Wang--Yang \cite{SWY21}.
When  the parameter is at critical value, flow and variational methods can also yield
 existence and concentration results, see, e.g., 
\cite{LZ19,SZ21,SZ24,Z23,YY24}.

For the quasilinear mean field equation \eqref{maineq} with $N\ge2$, Esposito--Morlando \cite{EM15} proved mass quantization and nonresonant compactness results for smooth positive weights. In particular,  the noncontractibility of the corresponding barycenter space  yields the existence result  in  \cite[Theorem~1.3]{EM15}. 

Compared with the planar case $N=2$, the analysis for $N>2$ is more delicate. The reasons include that the $N$-Dirichlet energy is no longer  quadratic, $\Delta_N$ may degenerate if $\nabla u=0$, and the Dirichlet inverse is nonlinear. Moreover, potentials associated with different concentration points cannot be added as in the linear case, and gradients may be unbounded at negative singularities.  

Recently, Ahmedou--Tang--Wang \cite{ATW26} computed the Leray--Schauder degree of \eqref{maineq} for $-1<\al_i\le0$ and determined the homotopy type of sufficiently low sublevels. Here we study the degree for arbitrary exponents greater than $-1$.


The main difficulty in extending \cite{ATW26} arises from positive exponents. For $\al_i>0$, the concentration mass at $p_i$ exceeds the ordinary concentration mass. The local degree calculation must therefore distinguish concentration at $p_i$ from concentration near $p_i$. Moreover, a zero degree does not exclude the existence of solutions. We are thus led to the following two questions.
\begin{enumerate}[label=\textup{(Q\arabic*)},leftmargin=3.2em,
                  itemsep=2pt,topsep=4pt]
\item How can one compute the Leray--Schauder degree for all
nonresonant parameters when positive exponents are present?
\item When the degree is zero, which additional topological
condition can still guarantee a weak solution?
\end{enumerate}

To answer (Q1), we first determine the concentration masses and compute the local degree near each positive singularity. We then use entropy duality method, domain comparisons with shrinking holes and annular cuts, and a degree formula for disconnected domains to compute the degree on general domains.

To answer (Q2),  we first construct a trial map from the weighted barycenter space to low energy sublevels and a projection in the opposite direction. By Proposition~\ref{prop:homotopy}, the composition of the two maps is homotopic to the identity on the barycenter space. The noncontractibility of a nonempty weighted barycenter space therefore implies the existence of a weak solution, even when the degree vanishes. 

We omit the points with $\al_i=0$ before counting $m$ and indexing the singularities. Their factors are identically one, and their normalized mass is already included in the ordinary integer mass in \eqref{eq:sigma}. Thus $m=0$ corresponds to the regular case.

Denote
\[
 c_N=|\S ^{N-1}|\Big(\frac{N^2}{N-1}\Big)^{N-1},
 \quad \beta_i=(1+\al_i)^{N-1},\quad r=\frac{\rho}{c_N},
\]
where  $\S ^{N-1}$ denotes the unit sphere and
$|\S ^{N-1}|$ is its surface area. By Proposition~\ref{prop:blowup}, the mass is $c_N$ at an ordinary concentration point and $c_N\beta_i$ at $p_i$. In dimension two, $c_2=8\pi$.
Thus $r$ measures the total mass in units of $c_N$.
Write
\[
 I_-=\{i:\al_i<0\},\quad I_+=\{i:\al_i>0\},\quad
 \beta_I=\sum_{i\in I}\beta_i,\quad \beta_\varnothing=0.
\]
Here $I_-$ and $I_+$ index the negative and positive exponents, respectively, and $\beta_I$ sums the weights with indices in $I$.
We use $\N=\{1,2,\ldots\}$, $\N_0=\N\cup\{0\}$, and
$|I|$ or $\#I$ for the number of elements of a finite set $I$.

An unbounded sequence of solutions can have several concentration points. The possible sums of concentration masses form the locally finite set
\be\label{eq:sigma}
 \Si_H=c_N\{\ell+\beta_I:\ell\in\N_0,\,
 I\subset\{1,\ldots,m\},\, \ell+|I|>0\}.
\ee
In \eqref{eq:sigma}, $\ell$ counts ordinary concentration points
and $I$ indexes the prescribed concentration points. By Proposition~\ref{prop:blowup}, every blow-up limit of the parameter belongs to $\Si_H$. Membership in $\Si_H$ is only a necessary condition for blow-up. We call $\rho$
\emph{nonresonant} when $\rho\notin\Si_H$.

To define the degree, let $E=W^{1,N}_0(\Om)$ be the Sobolev
space of functions with zero boundary trace, with norm
$\norm u_E=\norm{\nabla u}_{L^N(\Om)}$. Write
\be\label{eq:solutionmap}
 \begin{aligned}
 K_{\rho,H}(u)&=(-\Delta_N)^{-1}
       \Big(\rho\frac{He^u}{\int_\Om He^u\, \d  x}\Big),\\
 d_\rho(\Om,H)&=\deg(\Id-K_{\rho,H},B_R^E,0),
 \quad B_R^E=\{u\in E:\norm u_E<R\}.
 \end{aligned}
\ee
In \eqref{eq:solutionmap}, $(-\Delta_N)^{-1}$ maps a right-hand side to the unique weak solution of the corresponding zero-boundary Dirichlet problem. The symbol $\Id$ denotes the identity on $E$. Fixed points of $K_{\rho,H}$ are
precisely weak solutions of \eqref{maineq}.
At a nonresonant parameter, Theorem~\ref{thm:compactness} and Proposition~\ref{prop:parametercompact} place all solutions in $B_R^E$ for sufficiently large $R$. By excision, $d_\rho$ in \eqref{eq:solutionmap} is independent of the choice of sufficiently large $R$.


For $n\in\mathbb Z$ and
$k\in\N_0$, we use the generalized binomial coefficients
\[
 \binom n0=1,\quad
 \binom nk=\frac{n(n-1)\cdots(n-k+1)}{k!}\quad(k\ge1).
\]
Here $\chi(\Om)$ is the Euler characteristic of $\Om$. Negative upper entries are included. For example,
$\binom{-1}{k}=(-1)^k$. Empty sums are zero.
For fixed $r$, only finitely many subsets and integer indices contribute to the degree in \eqref{eq:solutionmap}. The formula depends on the parity of $N$.

\begin{thm}\label{thm:degree}
Assume  the conditions in \eqref{eq:assumptions} hold,
and let $\rho>0$ satisfy $\rho\notin\Si_H$.
If $N$ is even, then
\be\label{eq:all-nonresonant-degree}
 d_\rho(\Om,H)=
 \sum_{\substack{I\subset\{1,\ldots,m\},\ \ell\in\N_0\\
                  \beta_I+\ell<r}}
 (-1)^{|I|+\ell}\binom{\chi(\Om)-m}{\ell}.
\ee
If $N$ is odd, then
\be\label{eq:odd-nonresonant-degree}
 d_\rho(\Om,H)=
 \sum_{\substack{I\subset I_-,\ J\subset I_+\\
                  k,\ell\in\N_0,\ k\le |I_+|-|J|\\
                  \beta_I+k+\ell+\sum_{j\in J}\min\{\beta_j,2\}<r}}
 (-1)^{|I|+k+\ell}2^k
 \binom{|I_+|-|J|}{k}\binom{\chi(\Om)-m}{\ell}.
\ee
\end{thm}

\begin{rem}
Equations \eqref{eq:all-nonresonant-degree} and \eqref{eq:odd-nonresonant-degree} depend only on $N$, $\chi(\Om)$, the weights $\beta_i$, and $r$. In particular, both degrees are independent of $h$ and the positions of the singularities. At $\rho=0$, $d_0(\Om,H)=1$.
\end{rem}

A nonzero value of the degree in Theorem~\ref{thm:degree} implies that  \eqref{maineq} has a solution.
\begin{cor}\label{cor:degree-existence}
Assume  the conditions in \eqref{eq:assumptions} hold.
If $\rho>0$, $\rho\notin\Si_H$, and
$d_\rho(\Om,H)\ne0$, then \eqref{maineq} admits a weak solution.
\end{cor}

To apply Corollary~\ref{cor:degree-existence}, one can simplify the degree formulas by summing over the integer index when $N$ is even or $I_+=\varnothing$.
Here $\lfloor s\rfloor$ denotes the largest integer not exceeding $s$.
\begin{cor}\label{cor:even-full-degree}
Assume  the conditions in \eqref{eq:assumptions} hold,
and let $\rho>0$ satisfy $\rho\notin\Si_H$.
If $N$ is even or $I_+=\varnothing$, then
\be\label{eq:full-degree-formula}
 d_\rho(\Om,H)=
 \sum_{\substack{I\subset\{1,\ldots,m\}\\\beta_I<r}}
 (-1)^{|I|}\binom{\lfloor r-\beta_I\rfloor-\chi(\Om)+m}
                       {\lfloor r-\beta_I\rfloor}.
\ee
\end{cor}

\begin{rem}
For $N=2$, \eqref{eq:full-degree-formula} is the singular
Dirichlet degree formula established by  Chen--Lin
\cite[Theorem~1.4]{CL15}. For nonpositive exponents in any dimension, \eqref{eq:full-degree-formula} recovers Ahmedou--Tang--Wang \cite[Theorem~1.1]{ATW26}. When $m=0$, the formula \eqref{eq:full-degree-formula} becomes
\[
 d_\rho(\Om,h)=
 \binom{\lfloor\rho/c_N\rfloor-\chi(\Om)}
       {\lfloor\rho/c_N\rfloor},\quad
 \rho>0,\quad \rho\notin c_N\N.
\]
Thus the regular case includes every positive nonresonant
parameter. For $N\ge3$, Theorem~\ref{thm:degree} extends \cite{ATW26} to positive exponents, with a different contribution in odd dimensions.
\end{rem}

For $r>1$ in the intervals below, the contributions of positive singularities in Theorem~\ref{thm:degree} take a simpler form.
For the next corollary and the variational result below, set
\[
 \beta_*=\min (\{1\}\cup\{\beta_i:i\in I_-\}).
\]
With $I_-=\varnothing$, the smallest normalized concentration mass is $\beta_*=1$. In general, $\beta_*$ is the minimum of the normalized ordinary and singular masses in \eqref{eq:sigma}.
Under the hypotheses of the next corollary, $\rho\notin\Si_H$ by \eqref{eq:sigma}.
\begin{cor}\label{prop:first-positive-band}
Assume  the conditions in \eqref{eq:assumptions} hold,
and suppose that $r\ne\beta_j$ for any $j\in I_+$.
\begin{enumerate}[label=\textup{(\roman*)},leftmargin=2em]
\item If $I_-=\varnothing$ and $1<r<2$, then
\be\label{eq:first-positive-band}
 d_\rho(\Om,H)=1-\chi(\Om)+(-1)^N\#\{j\in I_+:\beta_j>r\}.
\ee
\item If $I_-\ne\varnothing$ and
\[
 \max\{1,\beta_{I_-}\}<r<1+\beta_*,
\]
then
\be\label{eq:clean-supercritical-degree}
 d_\rho(\Om,H)=|I_-|-\chi(\Om)+(-1)^N\#\{j\in I_+:\beta_j>r\}.
\ee
\end{enumerate}
\end{cor}

A zero degree does not rule out solutions. To obtain a second
existence criterion, we study the energy
\[
 J_\rho(u)=\frac1N\int_\Om|\nabla u|^N\, \d  x
          -\rho\ln\int_\Om He^u\, \d  x,\quad u\in E.
\]
Critical points of $J_\rho$ solve \eqref{maineq}. On sufficiently low sublevels, the density concentrates near a finite set of points. Weighted barycenters describe the locations and mass fractions of the concentration points.

Set $P_\pm=\{p_i:i\in I_\pm\}$, and let
$Y\subset\overline\Om$ be a compact set containing $P_-$.
Denote by $\Prob(Y)$ the space of Borel probability measures on $Y$ with the weak topology. Every finitely supported measure $\si\in\Prob(Y)$ can be written as
\[
 \si=\sum_{j=1}^q t_j\delta_{x_j},\quad
 x_j\in Y,\quad t_j>0,\quad \sum_{j=1}^q t_j=1,
\]
where the points $x_j$ are distinct and $\delta_{x_j}$ denotes the Dirac measure at $x_j$. In particular, $\supp\si=\{x_1,\ldots,x_q\}$.

We define the weighted barycenter space by
\[
 \Xr(Y)=\Big\{\si\in\Prob(Y):\supp\si\text{ finite},\,
 \#(\supp\si\setminus P_-)
 +\sum_{\substack{i\in I_-\\p_i\in\supp\si}}\beta_i<r\Big\}.
\]
In $\Xr(Y)$, a negative point has weight $\beta_i<1$ and each other point has weight one. The weight condition depends only on the support of $\si$, independently of $t_j$. We use $\Xr(Y)$ in the existence argument. 
Recall that a space is noncontractible if its identity map is not homotopic to a constant map.

For $0\le\rho<c_N\beta_*$, one can minimize $J_\rho$ directly. At nonresonant parameters, a nonempty and noncontractible $\Xr(\overline\Om)$ yields a critical point through the topology of low sublevels.

\begin{thm}\label{thm:existence}
Assume  the conditions in \eqref{eq:assumptions} hold.
\begin{enumerate}[label=\textup{(\roman*)},leftmargin=2em]
\item If $0\le\rho<c_N\beta_*$, then equation \eqref{maineq}
admits a weak solution that minimizes $J_\rho$ on $E$.
\item If $\rho>0$, $\rho\notin\Si_H$, and
$\Xr(\overline\Om)$ is nonempty and noncontractible,
then \eqref{maineq} admits a weak solution.
\end{enumerate}
\end{thm}

\begin{rem}
By the regularity results  in Section~\ref{sec:2}, each solution $u$ in Theorem~\ref{thm:existence} satisfies  $u\in E\cap C^0(\overline\Om)$, and $u>0$ in $\Om$ for $\rho>0$. Moreover, $u\in C^{1,\eta}_{loc}(\overline\Om\setminus P_-)$ for some $\eta>0$. If $P_-=\varnothing$, then $u\in C^{1,\eta}(\overline\Om)$.
\end{rem}

The following family of examples satisfies the barycenter condition in Theorem~\ref{thm:existence}(ii) but  the degree is zero.
In the following,  $\simeq$ denotes homotopy equivalence, and
$\vee$ denotes the wedge sum. In particular, $\bigvee_{\nu=1}^M\S^d$ denotes $M$ spheres of dimension $d$ joined at a single point, and $\S^0$ consists of two points.  
A sufficiently deep sublevel means
$\{u\in E:J_\rho(u)\le-L\}$ with $L$ sufficiently large.

\begin{cor}\label{cor:eulerzero}
Let $n\ge k\ge2$ be integers with $k$ even, and let
$0<a\le b$ satisfy $kb+a<1$. 
Suppose the full list of negative weights consists of one copy of $a$, $n$ copies of $b$, and $\binom nk$ copies of $kb$. The positive exponents may be arbitrary.
Then, for any $\rho$ satisfying
\be\label{eq:homologyinterval}
 c_Nkb<\rho<c_N(kb+a),
\ee
equation \eqref{maineq} admits a solution but  $d_\rho(\Om,H)=0$.
Moreover,
\be\label{eq:eulerzerofamily}
 \Xr(\overline\Om)\simeq
 \Big(\bigvee_{\nu=1}^{\binom nk}\S ^0\Big)
 \vee\Big(\bigvee_{\nu=1}^{\binom nk}\S ^{k-1}\Big).
\ee
Over any field, the reduced homology of $\Xr(\overline\Om)$ has rank $\binom nk$ in degrees $0$ and $k-1$, and is zero in all other degrees.
The reduced homology of any sufficiently low energy sublevel has rank at least $\binom nk$ in degrees $0$ and $k-1$.
\end{cor}


We end this section  by explaining the main steps of the proof and how we extend the results of \cite{ATW26} to positive exponents.

The main difficulty comes from the different concentration masses at positive singularities and ordinary points. For $\al_i>0$, concentration at $p_i$ has mass $c_N\beta_i>c_N$, and concentration at nearby ordinary points has mass $c_N$. We therefore need a local comparison near each positive singularity; see \cite[Remark~1.8]{ATW26}. To achieve this, we first determine all possible concentration masses, then compute the degree on a ball with one positive singularity, and finally extend the degree calculation to general domains.

We begin with proving a compactness result. At a positive singularity, the main difficulty is to show that no mass remains outside the concentration points when blow-up occurs. If the limiting normalization is positive, comparison with a solution having the same Dirac mass yields the integrability condition \eqref{eq:strictgap}, which contradicts \eqref{eq:massidentity}. Proposition~\ref{prop:blowup} then determines the mass at each interior concentration point and rules out boundary concentration.  The compactness result in Theorem~\ref{thm:compactness} is used in the degree calculation.

Next, we consider a ball with one positive singularity $p_i$. For a nonresonant parameter $0<\rho<c_N\beta_i$, Proposition~\ref{prop:puncture-degree} shows that removing a sufficiently small ball around $p_i$ and imposing zero Dirichlet data on the new boundary do not change the degree. The resulting annulus is homotopy equivalent to $\S^{N-1}$, with
\[
 \chi(\S ^{N-1})=1+(-1)^{N-1}.
\]
The different values of $\chi(\S^{N-1})$ in even and odd dimensions explain the two local degree formulas. For $\rho>c_N\beta_i$, the Pohozaev identity rules out solutions on the original ball, so the degree is zero. Combining the two cases, we obtain Lemma~\ref{lem:positive-local-factor}, including the term $\min\{\beta_i,2\}$ in odd dimensions.

To extend the local degree calculation to general domains, we use the relation between the degree and the Euler characteristic established in \cite[Theorem~3.10]{ATW26}; see Proposition~\ref{prop:all-parameter-euler}. By \eqref{eq:all-parameter-euler}, computing the degree reduces to computing the Euler characteristic of sufficiently low energy sublevels. Lemma~\ref{lem:nonlinear-entropy-duality} provides an equivalent description of these sublevels in terms of probability densities. We can then remove a small amount of mass near suitable annuli. By comparing the energies before and after removing the annuli, we obtain the degree formula \eqref{eq:annular-cech-degree} in Section~\ref{sec:4}.

After removing the annuli, the equations on different connected components still share a single normalizing integral. To separate these equations, we consider $-\Delta_Nu=\lam He^u$ with the same small parameter $\lam>0$ on each component. The product property of the degree yields the formula in Proposition~\ref{prop:disconnected-convolution}. Combining the proposition with the annular comparison and induction on the number of positive singularities, we prove Theorem~\ref{thm:degree}.

Our degree calculation extends \cite{ATW26} to positive exponents. We use the degree--Euler relation and the degree formula for nonpositive exponents. Proposition~\ref{prop:blowup} determines the concentration masses at positive singularities, while the comparisons on domains with small holes and annular cuts determine their contributions to the degree. For positive exponents, \eqref{eq:annular-cech-degree} determines the Euler characteristic of low energy sublevels. To determine their full homotopy type, further analysis is needed because concentration at a positive singularity and concentration at nearby ordinary points have different energy costs.

For obtaining existence results, a simpler topological argument is sufficient. We choose a compact subset of $\Om$ that contains all negative singularities, avoids all positive singularities, and is a deformation retract of $\overline\Om$. On this subset, we construct trial functions as in \cite[Proposition~4.1]{ATW26}. By Proposition~\ref{prop:homotopy}, the composition of the trial map and the projection is homotopic to the identity on the weighted barycenter space. Thus, if the weighted barycenter space is noncontractible, the equation admits a solution. Corollary~\ref{cor:eulerzero} shows how this argument yields existence even when the degree is zero. See also \cite{JKM15} for similar homology arguments in Toda systems.

The remaining sections divide the proof as follows. In
Section~\ref{sec:2}, we prove a compactness result.
Section~\ref{sec:3} constructs the barycenter maps and
establishes the existence criteria. In
Section~\ref{sec:4},  we develop the local degree calculations
and the domain comparisons.  The proofs of the main results are completed in Section~\ref{sec:5}.

\section{Compactness and concentration of solutions}\label{sec:2}
\subsection{Exponential estimates and weak solutions}
For the variational arguments, denote 
\[
 Q(u)=\frac1N\int_\Om|\nabla u|^N\, \d  x,\quad
 M(u)=\int_\Om He^u\, \d  x.
\]
So  $J_\rho=Q-\rho\ln M$. We write $J_\rho^a=\{u\in E:J_\rho(u)\le a\}$. A weak solution of \eqref{maineq} is a function $u\in E$ satisfying
\[
 \int_\Om|\nabla u|^{N-2}\nabla u\cdot\nabla v\, \d  x
 =\frac{\rho}{M(u)}\int_\Om He^uv\, \d  x,
 \quad\forall\, v\in E.
\]

Set
\[
 q=\frac N{N-1},\quad \ga_N=\frac{N^2}{N-1},\quad
 S_N=|\S ^{N-1}|,
\]
so $c_N=S_N\ga_N^{N-1}$.
By the sharp Moser inequality \cite[Theorem~1]{M71} and Young's inequality, we have
\be \label{eq:MT}
 \ln \int_\Om e^u\, \d  x\le\frac{Q(u)}{c_N}+C.
\ee 
Since $\al_i>-1$, we have $H\in L^{s_0}(\Om)$ for any $1<s_0<\min_{i\in I_-}(-1/\al _i)$. If  $I_-=\varnothing$, we can take any  $s_0>1$.  By H\"older's inequality and \eqref{eq:MT}, for any $\var>0$,
\[
 \ln  M(u)\le\Big(\frac1{c_N\beta_*}+\var \Big) Q(u)+C_\var .
\]
Hence $J_\rho$ is coercive for $\rho<c_N\beta_*$.

On bounded subsets of $E$, the functions $e^u$ are uniformly bounded in $L^a(\Om)$ for any  finite $a\ge1$. If $u_j\weak u$ in $E$, compact Sobolev embeddings and uniform integrability yield $He^{u_j}\to He^u$ in $L^t(\Om)$ for some $t>1$. By Jensen's inequality applied to the probability measure $H\,\d x/M(0)$, the integral $M(u)$ has a uniform positive lower bound on bounded subsets of $E$. Moreover, $J_\rho\in C^1(E)$ and
\[
 J_\rho'(u)=Q'(u)-\rho(\ln M)'(u),\quad
 \langle Q'(u),v\rangle=\int_\Om|\nabla u|^{N-2}\nabla u\cdot\nabla v\, \d  x.
\]
For any  weakly convergent sequence $u_j\weak u$ in $E$, we have
$(\ln M)'(u_j)\to(\ln M)'(u)$ in $E^*$.
The Dirichlet operator $Q'=-\Delta_N:E\to E^*$ is a homeomorphism and satisfies
\be \label{eq:monotonicity}
 \langle Q'(u)-Q'(v),u-v\rangle\ge d_N\norm{u-v}_E^N,
\ee 
for some constant $d_N>0$.
Consequently, every bounded sequence with $\rho_j\to\rho$ and $J_{\rho_j}'(u_j)\to0$ has a strongly convergent subsequence in $E$. All applications below involve bounded sequences. We apply the direct method in the proof of Theorem~\ref{thm:existence}(i).

To study regularity near negative singularities, 
we compare two solutions with the same boundary values. 
Suppose $-\Delta_Nu=f$, $-\Delta_Nv=g$, $u-v\in E$, 
and $f-g\in L^t(\Om)$ for some $t>1$. Then
\be \label{eq:stability}
 \norm{u-v}_{L^\infty(\Om)}\le C_t\norm{f-g}_{L^t(\Om)}^{1/(N-1)}.
\ee 
To prove \eqref{eq:stability}, we test the difference 
of the two equations with $(u-v-k)^+$ for $k\ge0$. 
By \eqref{eq:monotonicity} and the Sobolev inequality, 
for any finite $a>Nt'$ with $t'=t/(t-1)$,
\[
 \norm{(u-v-k)^+}_{L^a(\Om)}^{N-1}
 \le C\norm{f-g}_{L^t(\Om)}|\{x\in\Om:u(x)-v(x)>k\}|^{1/t'-1/a}.
\]
Since $(a/t'-1)/(N-1)>1$, the standard level-set 
iteration yields an $L^\infty$ estimate for $(u-v)^+$. 
Interchanging $u$ and $v$ proves \eqref{eq:stability}.

For zero Dirichlet data and $f\in L^t(\Om)$, we approximate 
$f$ by the bounded functions 
$f_T=\max\{-T,\min\{f,T\}\}$. 
By \eqref{eq:stability} with $1<b<t$, the corresponding 
solutions converge uniformly as $T\to\infty$. 
The convergence is uniform over bounded families in 
$L^t(\Om)$. For each fixed $T$, the solutions satisfy 
uniform $C^{1,\eta}(\overline\Om)$ estimates. 
Hence the solutions with the original data are continuous 
on $\overline\Om$ and form a relatively compact family 
in $C^0(\overline\Om)$ whenever the right-hand sides 
are uniformly bounded in $L^t(\Om)$.

Let $u$ be a critical point of $J_\rho$. Testing 
\eqref{maineq} with the negative part of $u$ yields $u\ge0$.
Since the right-hand side of \eqref{maineq} is in 
$L^t(\Om)$ for some $t>1$, we have 
$u\in C^0(\overline\Om)$ by \eqref{eq:stability}. 
Away from $P_-$, the right-hand side is locally bounded. 
By the interior and boundary regularity results in 
\cite[Theorems~A.1 and A.4]{EM15}, we obtain
$u\in C^{1,\eta}_{loc}(\overline\Om\setminus P_-)$
for some $\eta>0$. If $\rho>0$, the strong maximum 
principle implies $u>0$ in $\Om$.

We also use an exponential estimate for solutions with 
general boundary values. Let $D$ be a smooth bounded domain 
and suppose $z\in W^{1,N}(D)$ solves $-\Delta_N z=f$ in $D$,
where $0\ne f\ge0$ and $f\in L^1(D)$. Let $\varphi$ be the
$N$-harmonic function with the same trace as $z$. Then
\be\label{eq:BM}
 \int_D\exp\Big(
 \frac{a|z-\varphi|}{\norm{f}_{L^1(D)}^{1/(N-1)}}
 \Big)\,\d x\le C(a,N,|D|),
 \quad 0<a<d_N^{1/(N-1)}N S_N^{1/(N-1)}.
\ee
If $z=0$ on $\pa D$, then $\varphi=0$, and 
\eqref{eq:BM} holds for the  range
$0<a<N S_N^{1/(N-1)}$ by \cite[Lemma~2.1]{EM15}.
For the estimate with general boundary data, set
$w=z-\varphi\ge0$, by comparison, and test the difference
of the weak equations with $\min\{w,k\}$ for $k>0$. Then
\[
 \int_D|\nabla\min\{w,k\}|^N\,\d x
 \le \frac{k\|f\|_{L^1(D)}}{d_N}.
\]
Apply the sharp Moser inequality to $\min\{w,k\}$ divided by
its gradient norm. On $\{w\ge k\}$ this gives
\[
 |\{w\ge k\}|\le C_N|D|
 \exp\Big(
 -N S_N^{1/(N-1)}
 \Big(\frac{d_N}{\|f\|_{L^1(D)}}\Big)^{1/(N-1)} k
 \Big).
\]
Integration in $k$ proves \eqref{eq:BM}. The same truncation
argument applies to two energy solutions with the same trace
and ordered right-hand sides, with $\|f\|_{L^1(D)}$ replaced
by the $L^1(D)$ norm of their difference.

Let $u_j$ be solutions with zero boundary values on $D$ 
and nonnegative right-hand sides $f_j$ satisfying
\[
 \sup_j\int_D f_j\,\d x<\infty.
\]
Up to a subsequence, there exists a function $u$ such that
\[
 u_j\to u\quad\text{ in }\,W^{1,b}(D),\quad 1\le b<N,
\]
and
\[
 |\nabla u_j|^{N-2}\nabla u_j
 \to|\nabla u|^{N-2}\nabla u
 \quad\text{ in }\,L^1(D;\R^N).
\]
Thus $u$ solves the equation with the limiting measure
in $D$, and $u\in W^{1,b}_0(D)$ for any  $1\le b<N$.
The zero Sobolev trace follows from strong convergence in
$W^{1,b}(D)$ and closedness of $W^{1,b}_0(D)$.
For fixed boundary values admitting an extension
$g\in W^{1,N+\var}(D)$ for some $\var>0$, the same
convergence holds and $u-g\in W^{1,b}_0(D)$ for any
$1\le b<N$.

For the gradient estimate, let $\varphi$ be the
$N$-harmonic function with the same boundary values
and set $w_j=u_j-\varphi\ge0$ by the comparison principle. By \eqref{eq:BM}, $w_j$ is uniformly
bounded in each finite $L^a(D)$. Testing the difference
equation with $1-(1+w_j)^{-1}$ yields
\[
 d_N\int_D\frac{|\nabla w_j|^N}{(1+w_j)^2}\,\d x
 \le\sup_\ell\|f_\ell\|_{L^1(D)}.
\]
H\"older's inequality gives a uniform bound for $\nabla w_j$
in $L^b(D)$ for $1\le b<N$, using the $L^{2b/(N-b)}$ bound
for $1+w_j$. Rellich compactness gives a subsequence Cauchy
in $L^1(D)$. For $\varepsilon>0$, test the equations for
$u_j$ and $u_\ell$ with the truncation of $u_j-u_\ell$ at
$\pm\varepsilon$. The common trace and monotonicity give
\[
 d_N\int_{\{|u_j-u_\ell|<\varepsilon\}}
       |\nabla u_j-\nabla u_\ell|^N\,\d x
 \le\varepsilon\|f_j-f_\ell\|_{L^1(D)}.
\]
First letting $j,\ell\to\infty$ and then $\varepsilon\downarrow0$
shows that the gradients are Cauchy in measure. For each $b<N$,
the uniform estimate at some $b'\in(b,N)$ and Vitali's theorem
give strong convergence in $L^b(D)$. Taking $b>N-1$ also
gives strong $L^1$ convergence of the vector fields.
See \cite[Lemma~2.2 and Appendix~A]{ATW26} for this truncation
argument.

For existence and regularity of solutions with measure 
data, we refer to \cite[Theorems~1.1--1.2 and 
Proposition~6.1]{DHM00}. To prove uniqueness, let $u$ 
and $v$ be two solutions with the same measure and 
boundary values. The common $W^{1,N+\var}(D)$ extension
can be extended to $W^{1,N+\var}(\R^N)$ because $D$ is smooth.
For $b<N$ sufficiently close to $N$,
\cite[Proposition~6.2]{DHM00} implies that both gradients
are in weak-$L^N(D)$. Moreover, $u-v\in W^{1,1}_0(D)$.
Since $D$ is smooth and the vector field 
$|\xi|^{N-2}\xi$ satisfies the required structure and 
monotonicity conditions, \cite[Theorem~1.2]{DHM00} 
implies $u=v$.

In the local comparison arguments, we approximate the
nonnegative measures by smooth nonnegative functions of
uniformly bounded total mass. The preceding compactness
then applies in $W^{1,s}(D)$ for $1\le s<N$.
The $C^1$ boundary traces in the local problems admit 
the required extensions. Uniqueness in the above 
solution class identifies the limit of the approximations.

\subsection{Local analysis of isolated singularities}
In the next lemma, we adapt the arguments in 
\cite[Theorem~1.1 and Proposition~2.3, pp.~9--10]{E21} 
to study isolated singularities with bounded positive coefficients.
\begin{lem}\label{lem:singularity}
Let $\al>-1$, $b>0$, and $\si>0$. Assume that
$V$ satisfies $0<v_-\le V(x)\le v_+<\infty$ in $B_R$.  Let $U\in\bigcap_{1\le s<N}W^{1,s}(B_R)\cap C^{1,\eta}_{loc}(B_R\setminus\{0\})$ be a solution of
\[ 
 -\Delta_N U=f+\si\delta_0,
\]
where $f=bV(x)|x|^{N\al}e^U\in L^1(B_R)$.
Set $a=(\si/S_N)^{1/(N-1)}$. Then $a<N(1+\al)$ and,
as $x\to0$,
\be \label{eq:residual-asymptotic}
  U(x)=a\ln \frac1{|x|}+O(1),\quad 
 |x|\nabla U(x)=-a\frac{x}{|x|}+o(1).
\ee 
\end{lem}
\begin{proof}
For $0<r<R$, let $Z_r$ solve
\[
 -\Delta_N Z_r=\si\delta_0\quad\text{ in }\,B_r,
 \quad Z_r=U\quad\text{ on }\,\pa B_r.
\]
For any constant $C$, the function $a\ln(r/|x|)+C$
satisfies the same equation in $B_r$ and has boundary
value $C$. By the comparison principle, we obtain
\[
 a\ln\frac r{|x|}+\min_{\pa B_r}U
 \le Z_r(x)\le
 a\ln\frac r{|x|}+\max_{\pa B_r}U.
\]
Hence $Z_r(x)=a\ln(1/|x|)+O_r(1)$. We next show that
$w_r:=U-Z_r\ge0$. To see this, approximate $\si\delta_0$
by the same smooth nonnegative functions $g_j$ in both
equations, and approximate $f$ by smooth nonnegative
functions $f_j$ in $L^1(B_r)$. The corresponding solutions
have the same boundary values and right-hand sides
$g_j+f_j$ and $g_j$, respectively. By the comparison
principle and uniqueness for measure data, passing to
the limit yields $U\ge Z_r$. Since $V\ge v_->0$, we obtain
\[
 \int_{B_r}|x|^{N\al-a}\,\d x
 \le C_r\int_{B_r}f\,\d x<\infty.
\]
Therefore, $a<N(1+\al)$. Choose $p>1$ sufficiently close
to one such that $p(N\al-a)>-N$. By the asymptotic
behavior of $Z_r$,
\be\label{eq:backgroundLp}
 V|x|^{N\al}e^{Z_r}\in L^p(B_r).
\ee
The corresponding $L^p$ norm may depend on $r$.

We next estimate $w_r$. Let $U_j$ and $Z_{r,j}$ be the 
approximate solutions constructed above, and set
$w_j=U_j-Z_{r,j}\ge0$. Testing the difference of the 
two equations with $\min\{t,(w_j-k)^+\}$ and using 
\eqref{eq:monotonicity}, we obtain
\[
 d_N\int_{\{x\in B_r:k<w_j(x)<k+t\}}
 |\nabla w_j|^N\,\d x
 \le t\int_{\{x\in B_r:w_j(x)>k\}}f_j\,\d x
\]
for all $k\ge0$ and $t>0$.
Arguing as in the proof of \eqref{eq:BM}, we have
\[
 \int_{B_r}\exp\Big(
 \frac{\kappa_N w_j}{\norm{f_j}_{L^1(B_r)}^{1/(N-1)}}
 \Big)\,\d x\le C_N|B_r|
\]
for some constant $\kappa_N>0$ independent of $j$ and $r$.
Since $w_j\to w_r$ almost everywhere and $f_j\to f$ 
in $L^1(B_r)$, Fatou's lemma yields
\[
 \int_{B_r}\exp\Big(
 \frac{\kappa_N w_r}{\norm{f}_{L^1(B_r)}^{1/(N-1)}}
 \Big)\,\d x\le C_N|B_r|.
\]
Since $\norm{f}_{L^1(B_r)}\to0$ as $r\to0$, for any 
finite $s\ge1$, we have $e^{w_r}\in L^s(B_r)$
provided $r$ is sufficiently small.
Moreover,
\[
 f=bV|x|^{N\al}e^{Z_r}e^{w_r}.
\]
Combining \eqref{eq:backgroundLp} with H\"older's 
inequality, we obtain $f\in L^t(B_r)$ for some $t>1$.

We now choose smooth approximations with $f_j\to f$
in $L^t(B_r)$. Since $w_j=0$ on $\pa B_r$,
\eqref{eq:stability} yields
\[
 \norm{w_j}_{L^\infty(B_r)}
 \le C_t\norm{f_j}_{L^t(B_r)}^{1/(N-1)}\le C.
\]
Passing to the limit, we obtain $w_r\in L^\infty(B_r)$.
Hence
\[
 U(x)=Z_r(x)+w_r(x)=a\ln\frac1{|x|}+O(1).
\]

To prove the gradient expansion, set
\[
 U_\var(y)=U(\var y)-a\ln\frac1\var.
\]
On every compact annulus in $\R^N\setminus\{0\}$,
$U_\var$ is uniformly bounded and satisfies
\[
 -\Delta_N U_\var
 =b\var^{N(1+\al)-a}V(\var y)|y|^{N\al}e^{U_\var}.
\]
Since $a<N(1+\al)$, the right-hand side converges
uniformly to zero on compact annuli as $\var\to0$.
By local $C^{1,\eta}$ estimates, up to a subsequence,
$U_\var$ converges in
$C^1_{loc}(\R^N\setminus\{0\})$
to an $N$-harmonic function $U_0$ satisfying
\[
 U_0(y)=a\ln\frac1{|y|}+O(1)
 \quad\text{ in }\, \R^N\setminus\{0\}.
\]
By \cite[Corollary~2.2]{KV86}, we have
$U_0(y)=a\ln(1/|y|)+C$ for some constant $C$;
see also \cite[Theorem~1.1]{E21}.
Since all possible limits have the same gradient,
we obtain
\[
 \var\nabla U(\var y)\to-ay
 \quad\text{uniformly on }|y|=1.
\]
The gradient expansion in \eqref{eq:residual-asymptotic}
then follows.
\end{proof}

\subsection{Concentration masses and parameter dependence}
We now determine the concentration masses using 
\eqref{eq:residual-asymptotic} and the Pohozaev identity.
\begin{prop}\label{prop:blowup}
Fix distinct points $p_1,\ldots,p_m\in\Om$ and exponents
$\al_i>-1$. Let $u_j$ be solutions of \eqref{maineq}
with bounded parameters $\rho_j\ge0$ and background weights
$h_j$ satisfying
\[
 c\le h_j\le C,\quad
 \norm{\nabla h_j}_{L^\infty(\Om)}\le C
\]
for some fixed constants $c,C>0$.
If $\norm{u_j}_{L^\infty(\Om)}\to\infty$, then, up to
a subsequence, there exists a nonempty finite set
$S\subset\Om$ such that
\[
 \rho_j\frac{H_je^{u_j}}{M_j(u_j)}\,\d x
 \weak\sum_{p\in S}\si_p\delta_p
 \quad\text{ on }\, \overline\Om,
\]
where
\[
 \si_p=
 \left\{\begin{aligned}
 &c_N &&\text{ if }\, p\notin P,\\
 &c_N\beta_i &&\text{ if }\, p=p_i.
 \end{aligned}\right.
\]
Moreover, $\rho_j/M_j(u_j)\to0$.
No nonresonance assumption is required.
\end{prop}
\begin{proof}
\textbf{Step 1:} The concentration set. 

Up to a subsequence, we have $\rho_j\to\rho_\infty\ge0$
and $h_j\to h_\infty$ uniformly on $\overline\Om$.
By the maximum principle, $u_j\ge0$. Since $u_j=0$
when $\rho_j=0$, we may assume $\rho_j>0$.
Write
\[
 H_j(x)=h_j(x)\prod_{i=1}^m|x-p_i|^{N\al _i},\quad
 M_j(w)=\int_\Om H_je^w\, \d  x,
\]
and set
\[
 b_j=\frac{\rho_j}{M_j(u_j)},\quad  v_j=u_j+\ln  b_j,\quad 
 f_j=H_j e^{v_j}.
\]
Then
\[
 0<b_j\le C,\quad
 \int_\Om f_j\,\d x=\rho_j,\quad
 -\Delta_Nu_j=f_j.
\]
By \eqref{eq:BM} with zero boundary values,
$u_j$ is uniformly bounded in $L^a(\Om)$ for any 
finite $a\ge1$, and $e^{\delta u_j}$ is uniformly
bounded in $L^1(\Om)$ for some $\delta>0$. Up to a further subsequence, we have
\[
 b_j\to b\ge0,\quad  f_j\, \d  x\weak\mu~\text{ on }\,\overline\Om,
 \quad  u_j\to U~\text{ in }\, W^{1,s}_0(\Om),\quad 1\le s<N.
\]
Moreover,
\[
 |\nabla u_j|^{N-2}\nabla u_j
 \to|\nabla U|^{N-2}\nabla U
 \quad\text{ in }\, L^1(\Om;\R^N).
\]

At a singular point $p_i\in P$, set
$\widehat\theta=\min\{1,1+\al_i\}$.
At an ordinary point or a boundary point, set
$\widehat\theta=1$. We first show that if
\[
 \limsup_{j\to\infty}\int_D f_j\,\d x
 <d_NS_N(N\widehat\theta)^{N-1}
\]
on a neighborhood $D$, then $u_j$ and $v_j$ are 
uniformly bounded above on a smaller neighborhood.

Indeed, choose $a>1/\widehat\theta$ such that
\[
 a\Big(\limsup_{j\to\infty}\int_D f_j\,\d x\Big)^{1/(N-1)}
 <d_N^{1/(N-1)}NS_N^{1/(N-1)}.
\]
We can also choose $b_0>1$ such that
\[
 \sup_j\norm{H_j}_{L^{b_0}(D)}<\infty,\quad
 \frac1{b_0}+\frac1a<1.
\]
Let $\varphi_j$ be the $N$-harmonic function in $D$
with the same boundary values as $v_j$. By comparison,
$\varphi_j\le v_j$. Since $v_j$ is uniformly bounded
above in $L^N(D)$, the positive part of $\varphi_j$
is uniformly bounded in $L^N(D)$.
The local estimates for $N$-harmonic functions and
\eqref{eq:BM} then imply that $e^{v_j}$ is locally
bounded in $L^a(D)$. By H\"older's inequality,
$f_j=H_je^{v_j}$ is locally bounded in $L^t(D)$
for some $t>1$.

Combining the local regularity estimates with the
uniform $L^N(\Om)$ estimate for $u_j$, we obtain
uniform upper estimates for $u_j$ and $v_j$ on a
smaller neighborhood. Near the boundary, the same
argument applies using $v_j=\ln b_j$ and $u_j=0$
on $\pa\Om$; see
\cite[Lemma~2.2 and Theorems~A.1--A.4]{EM15}.
In particular, when $\al_i>0$, the small-mass condition
is the same as that at an ordinary point.

Choose $\var_0>0$ smaller than all the thresholds above
and set
\[
 S=\{x\in\overline\Om:\mu(\{x\})\ge\var_0\}.
\]
Since $\mu$ is a finite measure, $S$ is finite.
For any  $x\notin S$, the local estimates above imply
that $u_j$ and $v_j$ are uniformly bounded above near $x$.
Moreover, $f_j\le CH_j$ near $x$. Since $H_j$ is uniformly
bounded in $L^{s_0}(\Om)$ for some $s_0>1$, no mass
concentrates outside $S$. Thus $S$ is exactly the
concentration set.

Set $F=bH_\infty e^U$. By Fatou's lemma, $F\in L^1(\Om)$.
Away from $S\cup P_-$, the right-hand sides are locally
uniformly bounded. The interior and boundary regularity
estimates yield, up to a subsequence,
\[
 u_j\to U\quad\text{ in }\, 
 C^1_{loc}(\overline\Om\setminus(S\cup P_-)).
\]
To prove uniform convergence near a negative singularity
$p_i\notin S$, choose a small ball $B$ centered at $p_i$
such that $\overline B\subset\Om\setminus S$ and $B$
contains no other singularities. On $\pa B$, we have
$g_j:=u_j|_{\pa B}\to g:=U|_{\pa B}$ uniformly.
Since $u_j$ is uniformly bounded in $B$ and $H_j$
is dominated by a fixed function in $L^{s_0}(B)$,
we obtain $f_j\to F$ in $L^t(B)$ for some $t>1$.

Let $z$ solve $-\Delta_N z=F$ in $B$ with boundary
value $g$. By the comparison principle and
\eqref{eq:stability},
\[
 \norm{u_j-z}_{L^\infty(B)}
 \le \norm{g_j-g}_{L^\infty(\pa B)}
 +C\norm{f_j-F}_{L^t(B)}^{1/(N-1)}\to0.
\]
Since $u_j\to U$ almost everywhere, we have $z=U$ in $B$.
Consequently,
\be\label{eq:away}
 \begin{aligned}
 u_j&\to U\quad\text{ in }\, 
 C^0_{loc}(\overline\Om\setminus S),\\
 u_j&\to U\quad\text{ in }\, 
 C^1_{loc}(\overline\Om\setminus(S\cup P_-)).
 \end{aligned}
\ee

Since $\norm{u_j}_{L^\infty(\Om)}\to\infty$, the local
estimates imply $S\ne\varnothing$. Outside $S$,
the limiting measure is $F\,\d x$.
Passing to the limit in $-\Delta_Nu_j=f_j$, we obtain
\be\label{eq:limit}
 -\Delta_N U=F+\sum_{p\in S\cap\Om}\si_p\delta_p
 \quad\text{ in }\, \Om,\quad
 \si_p=\mu(\{p\})>0.
\ee

\textbf{Step 2:} Boundary concentration.

Suppose $p\in S\cap\pa\Om$. Choose a small smooth domain
$D\subset\Om$ that coincides with $\Om$ near $p$,
contains no singular points, and satisfies
$\pa D\cap\Om\cap S=\varnothing$.
By \eqref{eq:limit}, we have
\[
 -\Delta_N U=F\quad\text{ in }\, D,\quad
 F=bH_\infty e^U\in L^1(D).
\]
The boundary values of $U$ on $\pa D$ admit a $C^1$
extension that vanishes near $p$. Let $\varphi$ be
the $N$-harmonic function with the same boundary
values as $U$, and let $U_T$ solve
\[
 -\Delta_N U_T=F_T:=\min\{F,T\}\quad\text{ in }\, D,
 \quad U_T=U\quad\text{ on }\, \pa D.
\]
Since $\norm{F}_{L^1(D)}$ can be made arbitrarily small
by choosing $D$ sufficiently small, \eqref{eq:BM} and
the boundedness of $\varphi$ imply
\[
 \int_D e^{aU_T}\,\d x\le C
\]
for some $a>N$, with $C$ independent of $T$.
By uniqueness for measure data, $U_T\to U$ in
$W^{1,s}(D)$ for $1\le s<N$.
Since $H_\infty$ is bounded in $D$, Fatou's lemma
implies $F\in L^a(D)$.

Let $\widehat U\in W^{1,N}(D)$ be the weak solution of
\[
 -\Delta_N\widehat U=F\quad\text{ in }\, D,\quad
 \widehat U=U\quad\text{ on }\, \pa D.
\]
Testing the difference of the equations for $U_T$
and $\widehat U$ with $U_T-\widehat U$, and using
\eqref{eq:monotonicity} and the Sobolev inequality,
we obtain
\[
 \norm{\nabla(U_T-\widehat U)}_{L^N(D)}^{N-1}
 \le C\norm{F_T-F}_{L^a(D)}\to0.
\]
Hence $U=\widehat U\in W^{1,N}(D)$.
By the boundary estimate in \cite[Theorem~A.1]{EM15},
$U$ is bounded near $p$. Consequently,
$F=bH_\infty e^U$ is also bounded near $p$.
By \cite[Theorem~A.4]{EM15}, we obtain
$C^{1,\eta}$ regularity of $U$ up to the boundary
near $p$.

We now apply the Pohozaev identity near $p$.
The main idea is to choose a suitable vector field so that
the boundary integral involving $\nabla u_j$ vanishes.
Set
\[
 D_R=\Om\cap B_R(p),\quad
 \Ga_R=\pa\Om\cap B_R(p),\quad
 S_R=\Om\cap\pa B_R(p).
\]
For sufficiently small $R>0$, choose
\[
 t_{j,R}=
 \frac{\int_{\Ga_R}(x-p)\cdot\nu(x)|\nabla u_j|^N\,\d S}
 {\int_{\Ga_R}\nu(p)\cdot\nu(x)|\nabla u_j|^N\,\d S},
 \quad y_{j,R}=p+t_{j,R}\nu(p),
\]
where $\nu$ denotes the outward unit normal.
If $\nabla u_j=0$ on $\Ga_R$, we set $t_{j,R}=0$.
By smoothness of $\pa\Om$, we have
\[
 |(x-p)\cdot\nu(x)|\le CR^2,\quad
 \nu(p)\cdot\nu(x)\ge\frac12,\quad x\in\Ga_R.
\]
Hence $|t_{j,R}|\le CR^2$ and, by the definition of
$y_{j,R}$,
\[
 \int_{\Ga_R}(x-y_{j,R})\cdot\nu(x)
 |\nabla u_j|^N\,\d S=0.
\]

Applying the Pohozaev identity on $D_R$ with the vector
field $x-y_{j,R}$, we obtain
\begin{align*}
 &\,\int_{D_R}
 [N+(x-y_{j,R})\cdot\nabla\ln H_j]f_j\,\d x\\
 =&\,\int_{\pa D_R}\Big[
 f_j(x-y_{j,R})\cdot\nu
 +|\nabla u_j|^{N-2}\pa_\nu u_j\,
 \nabla u_j\cdot(x-y_{j,R})-\frac1N|\nabla u_j|^N(x-y_{j,R})\cdot\nu
 \Big]\,\d S.
\end{align*}
Since $u_j=0$ on $\Ga_R$, we have
$\nabla u_j=(\pa_\nu u_j)\nu$ on $\Ga_R$.
By the choice of $y_{j,R}$, the integrals involving
$\nabla u_j$ over $\Ga_R$ vanish.
Moreover, $f_j=b_jH_j$ on $\Ga_R$, and hence
\[
 \Big|\int_{\Ga_R}
 f_j(x-y_{j,R})\cdot\nu\,\d S\Big|
 \le CR^{N+1}.
\]

For fixed sufficiently small $R$, we have
$S_R\cap S=\varnothing$. By \eqref{eq:away} and the
$C^{1,\eta}$ regularity of $U$ near $p$, the upper limit
of the integral over $S_R$ is at most $CR^N$ as
$j\to\infty$. Here we only use
$|y_{j,R}-p|\le CR^2$ and do not need the convergence
of $y_{j,R}$.

On the other hand, since $H_j$ is uniformly positive
with uniformly bounded first derivatives near $p$,
\[
 N+(x-y_{j,R})\cdot\nabla\ln H_j\ge N-CR
 \quad\text{ in }\, D_R.
\]
Passing to the limit in the Pohozaev identity yields
\[
 (N-CR)\mu(\overline\Om\cap B_R(p))\le CR^N.
\]
In particular,
\[
 (N-CR)\si_p\le CR^N.
\]
Letting $R\to0$, we obtain $\si_p=0$, which contradicts
$p\in S$.

The Pohozaev identity on $D_R$ can be justified by the
approximation in \cite[Lemma~2.3]{EM15} and by smoothing
the intersection of $\pa\Om$ and $\pa B_R(p)$.
For fixed $j$ and $R$, the boundary integrals converge
to those above since $D_R$ contains no singular points
and $u_j,H_j$ are $C^1$ near $\overline D_R$.

Therefore, $S\subset\Om$ and no concentration occurs
on $\pa\Om$.

\textbf{Step 3:} The limiting normalization.

We show that $b=0$. Suppose that $b>0$. Fix $p\in S$
and choose a small ball $B_R(p)$ containing no other
point of $S\cup P$. Write
\[
 H_j(x)=|x-p|^{N\al}V_j(x),
\]
where $\al=\al_i$ if $p=p_i$ and $\al=0$ otherwise.
The functions $V_j$ are uniformly positive and bounded
in $C^1(B_R(p))$. By Lemma~\ref{lem:singularity}
and \eqref{eq:limit}, we have
\be\label{eq:strictgap}
 \begin{gathered}
 a_p=(\si_p/S_N)^{1/(N-1)}<N(1+\al),\\
 U(x)=a_p\ln\frac1{|x-p|}+O(1),\\
 |x-p|\nabla U(x)=-a_p\frac{x-p}{|x-p|}+o(1).
 \end{gathered}
\ee

By the Pohozaev identity \cite[Lemma~2.3]{EM15},
for each fixed $j$, we have
 \be
\int_{B_R(p)}
 [N(1+\al)+(x-p)\cdot\nabla\ln V_j]f_j\,\d x
=R\int_{\pa B_R(p)}
 \Big[f_j+|\nabla u_j|^{N-2}(\pa_\nu u_j)^2
 -\frac1N|\nabla u_j|^N\Big]\,\d S.\label{eq:poho}
 \ee
If $H_j$ is not $C^1$ at $p$, we first apply the
Pohozaev identity on
$B_R(p)\setminus\overline B_{\var_\ell}(p)$.
Since $|\nabla u_j|^N,f_j\in L^1(\Om)$, we can choose
a sequence $\var_\ell\to0$ such that
\[
 \var_\ell\int_{\pa B_{\var_\ell}(p)}
 (|\nabla u_j|^N+f_j)\,\d S\to0.
\]
Thus all integrals over the inner boundary tend to zero,
and \eqref{eq:poho} follows.

We now let $j\to\infty$ and then $R\to0$ in
\eqref{eq:poho}. Since $V_j$ is uniformly positive
and bounded in $C^1$, we have
\[
\Big|\int_{B_R(p)}
 (x-p)\cdot\nabla\ln V_j\,f_j\,\d x\Big|
 \le CR\sup_j\rho_j.
\]
Moreover, $F\in L^1(\Om)$ implies
\[
 \int_{B_R(p)}F\,\d x\to0
 \quad\text{ as }\, R\to0.
\]
By \eqref{eq:away}, we can pass to the limit in the
boundary integrals of \eqref{eq:poho} for fixed $R$.
Using \eqref{eq:strictgap}, we obtain
\[
 R\int_{\pa B_R(p)}F\,\d S
 \le CR^{N(1+\al)-a_p}\to0.
\]
The gradient expansion in \eqref{eq:strictgap}
and the Pohozaev identity \eqref{eq:poho} then yield
\be\label{eq:massidentity}
 N(1+\al)\si_p=\frac{N-1}{N}S_Na_p^N,\quad
 \si_p=S_Na_p^{N-1}.
\ee
By \eqref{eq:massidentity}, we obtain
\[
 a_p=\frac{N^2}{N-1}(1+\al)
 >N(1+\al),
\]
which contradicts \eqref{eq:strictgap}. Hence $b=0$.

\textbf{Step 4:} The concentration masses.

Since $b=0$, the local estimates in Step 1 imply that
$f_j\to0$ in $L^t_{loc}(\overline\Om\setminus S)$
for some $t>1$.
For each $p_i\in P_-\setminus S$, choose a small ball
$B$ centered at $p_i$ such that
\[
 \overline B\subset\Om\setminus S,\quad
 \overline B\cap P=\{p_i\}.
\]
By \eqref{eq:away}, $u_j\to g:=U|_{\pa B}$ uniformly
on $\pa B$. Let $V$ be the $N$-harmonic function in $B$
with boundary value $g$. By the comparison principle
and \eqref{eq:stability}, we have
\[
 \norm{u_j-V}_{L^\infty(B)}
 \le \norm{u_j-g}_{L^\infty(\pa B)}
 +C\norm{f_j}_{L^t(B)}^{1/(N-1)}\to0.
\]
Hence $U=V$ in $B$.

Since $b=0$, we have $F=0$ in \eqref{eq:limit}, and thus
\[
 -\Delta_N U=\sum_{p\in S}\si_p\delta_p
 \quad\text{ in }\, \Om.
\]
In particular, $U$ is $N$-harmonic in $\Om\setminus S$.
By \cite[Theorem~A.5]{EM15}, for each $p\in S$, we have
\[
 U(x)=a_p\ln\frac1{|x-p|}+O(1),\quad
 |x-p|\nabla U(x)=-a_p\frac{x-p}{|x-p|}+o(1),
\]
where $a_p=(\si_p/S_N)^{1/(N-1)}$.
As in Step 3, passing to the limit in \eqref{eq:poho}
yields \eqref{eq:massidentity}. Therefore,
\[
 \si_p=\left\{\begin{aligned}
 &c_N &&\text{ if }\, p\notin P,\\
 &c_N(1+\al_i)^{N-1} &&\text{ if }\, p=p_i.
 \end{aligned}\right.
\]
The mass calculation is carried out directly near each
concentration point, without any assumption on the number
of bubbles near that point.
Since $\mu=\sum_{p\in S}\si_p\delta_p$ on $\overline\Om$,
we obtain
\[
 \rho_\infty
 =\lim_{j\to\infty}\int_\Om f_j\,\d x
 =\sum_{p\in S}\si_p\in\Si_H.
\]
\end{proof}


By Proposition~\ref{prop:blowup}, we can obtain the following
compactness result at nonresonant parameters.
\begin{thm}\label{thm:compactness}
Fix the singular points in $P$ and the exponents $\al_i>-1$.
Let $\Lam\subset[0,\infty)\setminus\Si_H$ be compact
and let $C_0>1$. Then all weak solutions of \eqref{maineq}
with
\[
 \rho\in\Lam,\quad C_0^{-1}\le h\le C_0,\quad
 \norm{\nabla h}_{L^\infty(\Om)}\le C_0
\]
are uniformly bounded in $E\cap L^\infty(\Om)$ and
relatively compact in $E\cap C^0(\overline\Om)$.
Moreover, for some $\eta\in(0,1)$, the solutions satisfy
uniform $C^{1,\eta}$ estimates on every compact subset
of $\overline\Om\setminus P_-$.

For fixed $h$, the set of solution pairs $(\rho,u)$
is compact in $\Lam\times(E\cap C^0(\overline\Om))$.
\end{thm}

\begin{proof}
Suppose that there exists a sequence of solutions $u_j$
with $\norm{u_j}_{L^\infty(\Om)}\to\infty$.
Since $\Lam$ is compact, up to a subsequence,
$\rho_j\to\rho_\infty\in\Lam$.
By Proposition~\ref{prop:blowup}, we have
$\rho_\infty\in\Si_H$, which contradicts
$\Lam\cap\Si_H=\varnothing$.
Hence all solutions are uniformly bounded in $L^\infty(\Om)$.

Testing \eqref{maineq} with $u$, we obtain
\[
 \norm{u}_E^N\le\rho\norm{u}_{L^\infty(\Om)}\le C.
\]
Set $f=\frac{\rho He^u}{M(u)}$.
Since $f\le CH$, we have $\norm{f}_{L^t(\Om)}\le C$
for some $t>1$.
By the compactness argument following \eqref{eq:stability},
every sequence of solutions $u_j$ has a subsequence
converging in $C^0(\overline\Om)$.
Set $f_j=\rho_jH_je^{u_j}/M_j(u_j)$.
By \eqref{eq:monotonicity},
\[
 d_N\norm{u_j-u_k}_E^N
 \le\int_\Om(f_j-f_k)(u_j-u_k)\,\d x
 \le C\norm{u_j-u_k}_{L^\infty(\Om)}\to0
\]
as $j,k\to\infty$.
Hence $u_j$ also converges strongly in $E$.

On every compact subset of $\overline\Om\setminus P_-$,
$f$ is uniformly bounded.
By the interior and boundary regularity estimates in
\cite[Theorems~A.1 and A.4]{EM15}, we obtain uniform
$C^{1,\eta}$ estimates for $u$.
For fixed $h$, passing to the limit in \eqref{maineq}
shows that the set of solution pairs $(\rho,u)$ is closed
and hence compact in $\Lam\times(E\cap C^0(\overline\Om))$.
\end{proof}

For the degree homotopies, we also need compactness for varying singularities and exponents.
\begin{prop}\label{prop:parametercompact}
Let $\mathcal K$ be a compact set of parameters
\[
 (\lam,(p_1,\ldots,p_m),(\al_1,\ldots,\al_m))
\]
satisfying
\[
 \lam\in[0,\infty)\setminus\Si_H,\quad
 p_i\in\Om,\quad p_i\ne p_j\ (i\ne j),\quad \al_i>-1,
\]
where $\Si_H$ is defined by \eqref{eq:sigma}. Assume that 
\[
 h\ge c,\quad \norm{h}_{C^1(\overline\Om)}\le C
\]
for some fixed constants $c,C>0$.
Then all solutions of \eqref{maineq} with parameters in
$\mathcal K$ are uniformly bounded and relatively compact
in $E\cap C^0(\overline\Om)$.

The solution maps $K_{\lam,H}$ are jointly continuous
in the parameters and $u$, and collectively compact on
bounded subsets of $E$. Hence the degree remains constant
along any continuous path in $\mathcal K$ for fixed $h$.
\end{prop}
\begin{proof}
Since $\mathcal K$ is compact, there exist constants
$d,\var,A>0$ such that
\[
 \dist(p_i,\pa\Om)\ge d,\quad
 |p_i-p_j|\ge d\quad(i\ne j),\quad
 -1+\var\le\al_i\le A.
\]
Choose $s_0>1$ such that $s_0(-1+\var)>-1$.
On each ball $B_{d/4}(p_i)$, only the factor
$|x-p_i|^{N\al_i}$ may be unbounded, while all other
factors are uniformly positive and bounded.
By integration in polar coordinates, we obtain a uniform
$L^{s_0}(\Om)$ estimate for $H$.

Consider a sequence of parameters in $\mathcal K$.
Up to a subsequence, we have
\[
 \lam_j\to\lam_\infty,\quad
 p_{i,j}\to p_i,\quad
 \al_{i,j}\to\al_{i,\infty},\quad
 h_j\to h_\infty\quad\text{uniformly on }\overline\Om.
\]
The function $h_\infty$ is positive and Lipschitz continuous.
Set
\[
 H_\infty(x)=h_\infty(x)
 \prod_{i=1}^m|x-p_i|^{N\al_{i,\infty}}.
\]
Since $H_j\to H_\infty$ almost everywhere in $\Om$,
Vitali's theorem yields
\be\label{eq:parameterweight}
 H_j\to H_\infty\quad\text{ in }\, L^s(\Om),
 \quad 1<s<s_0.
\ee
In particular, $\int_\Om H_j\,\d x\ge c_0>0$
for some constant $c_0$ independent of $j$.

Suppose now that there exists a sequence of solutions
$u_j$ with $\norm{u_j}_{L^\infty(\Om)}\to\infty$.
Since
\[
 \min\{1,1+\al_{i,j}\}\ge\min\{1,\var\}>0,
\]
the local estimates in Step 1 of the proof of
Proposition~\ref{prop:blowup} remain uniform.
Moreover, the singular points stay a fixed positive
distance from $\pa\Om$. Hence the argument in Step 2
of the same proof rules out boundary concentration.

Near a limiting singular concentration point $p_i$, write
\[
 H_j(x)=|x-p_{i,j}|^{N\al_{i,j}}V_{i,j}(x).
\]
For a sufficiently small fixed $R>0$, we have
\[
 V_{i,j}\ge c,\quad
 \norm{V_{i,j}}_{C^1(B_R(p_{i,j}))}\le C.
\]
We apply the Pohozaev identity \eqref{eq:poho} on
$B_R(p_{i,j})$ with $x-p_{i,j}$ in place of $x-p$.
The term involving $\nabla\ln V_{i,j}$ satisfies
\[
\Big|\int_{B_R(p_{i,j})}
 (x-p_{i,j})\cdot\nabla\ln V_{i,j}\,f_j\,\d x\Big|
 \le CR\sup_j\lam_j.
\]

Choose $R$ such that $\mu(\pa B_R(p_i))=0$.
For any $\delta\in(0,R)$ and sufficiently large $j$,
\[
 B_{R-\delta}(p_i)\subset B_R(p_{i,j})
 \subset B_{R+\delta}(p_i).
\]
Letting $j\to\infty$ and then $\delta\to0$, we obtain
\[
 \int_{B_R(p_{i,j})}f_j\,\d x
 \to\mu(B_R(p_i)).
\]
Since the spheres $\pa B_R(p_{i,j})$ stay away from
the concentration points and the limiting singularities,
the boundary integrals in \eqref{eq:poho} converge
by the local $C^1$-convergence as in \eqref{eq:away}.

If the limiting normalization is positive, then
Lemma~\ref{lem:singularity} implies
\[
 a_{p_i}<N(1+\al_{i,\infty}).
\]
On the other hand, letting $j\to\infty$ and then
$R\to0$ in the Pohozaev identity yields
\[
 a_{p_i}=\ga_N(1+\al_{i,\infty})
 >N(1+\al_{i,\infty}),
\]
which is a contradiction. At an ordinary concentration
point, the same argument applies with exponent zero.
Hence the limiting normalization is zero.

As in Step 4 of the proof of Proposition~\ref{prop:blowup},
we obtain
\[
 \si_p=\left\{\begin{aligned}
 &c_N &&\text{ if }\, p\notin P,\\
 &c_N(1+\al_{i,\infty})^{N-1}
 &&\text{ if }\, p=p_i.
 \end{aligned}\right.
\]
The singular points outside $S$ are removable by the
same argument as in the proof of Proposition~\ref{prop:blowup}.
Consequently,
\[
 \lam_\infty=\sum_{p\in S}\si_p\in\Si_{H_\infty},
\]
which contradicts the assumption on $\mathcal K$.
Hence all solutions are uniformly bounded in
$L^\infty(\Om)$.

Testing \eqref{maineq} with $u$ yields a uniform estimate
for $\norm u_E$. By \eqref{eq:parameterweight},
\eqref{eq:stability}, and \eqref{eq:monotonicity},
we obtain relative compactness in
$E\cap C^0(\overline\Om)$ and can pass to the limit
in \eqref{maineq}.

Finally, consider converging parameters and a sequence
$u_j\weak u$ in $E$. By the exponential integrability
of $E$ and \eqref{eq:parameterweight}, we have
\[
 \frac{\lam_jH_je^{u_j}}{M_j(u_j)}
 \to
 \frac{\lam_\infty H_\infty e^u}{M_\infty(u)}
 \quad\text{ in }\, L^b(\Om)
\]
for some $b>1$. By continuity of $(-\Delta_N)^{-1}$,
the solution maps in \eqref{eq:solutionmap} are jointly
continuous and collectively compact on bounded subsets
of $E$. The uniform estimates above allow us to choose
one ball in $E$ containing all fixed points along a
continuous path in $\mathcal K$. Hence the degree
remains constant along the path.

If $\al_i=0$, then $\beta_i=1$, which can be included
in the integer $\ell$ in \eqref{eq:sigma}.
Removing a point with zero exponent therefore does not
change either \eqref{maineq} or the resonance set
\eqref{eq:sigma}.
\end{proof}

By Proposition~\ref{prop:parametercompact}, we obtain
the following invariance property of the degree.
\begin{cor}\label{cor:degreepositions}
For fixed exponents $\al_i$ and a fixed nonresonant
parameter $\rho$, the degree $d_\rho(\Om,H)$ is
independent of the positions of the singular points
and the choice of the positive $C^1$-function $h$.
\end{cor}
\begin{proof}
Since $\Om$ is connected and $N\ge2$, we can move
the singular points continuously to any other distinct
interior points without collisions or touching $\pa\Om$.
Along such a path, the distances between the singular
points and from the singular points to $\pa\Om$ have
a uniform positive lower bound.

For any two positive functions
$h_0,h_1\in C^1(\overline\Om)$, consider
\[
 h_t=(1-t)h_0+th_1,\quad t\in[0,1].
\]
Since the exponents and $\rho$ are fixed, the resonance
set $\Si_H$ remains unchanged. By
Proposition~\ref{prop:parametercompact} and homotopy
invariance, the degree remains constant.
\end{proof}

\begin{rem}\label{rem:collision}
The uniform estimate in Proposition~\ref{prop:parametercompact}
can fail when the points are allowed to collide.
Translate so that $0\in\Om$ and, for sufficiently small $t>0$, set
$H_t(x)=|x|^{-3N/4}|x-te_1|^{-3N/4}$, where $e_1$ is the first
coordinate vector. Theorem~\ref{thm:existence}(i) gives a
solution for each $t$ at any fixed $0<\rho<c_N(1/4)^{N-1}$.
This parameter lies below the first resonance for any  $t$.

Suppose solutions were bounded in $L^\infty$ along $t\downarrow0$.
Since they are nonnegative, their normalizing integrals diverge:
$\int_{B_{t/2}(0)}H_t\,\d x\ge c t^{-N/2}$. On the complement
of any fixed ball about zero, the weights stay uniformly bounded.
The normalized right-hand sides therefore converge as measures
to $\rho\delta_0$. On the other hand, testing with the solutions
gives a uniform energy estimate. By H\"older's inequality, evaluation
against any $\phi\in C_c^\infty(\Om)$ then gives
\[
 \rho|\phi(0)|\le C\|\nabla\phi\|_{L^N(\Om)}.
\]
This contradicts the zero $N$-capacity of a point, since one
can choose $\phi(0)=1$ with arbitrarily small gradient norm.
Thus the solutions cannot remain uniformly bounded when
the points collide.
\end{rem}

\section{Weighted barycenters and variational existence}\label{sec:3}
\subsection{Concentration on low sublevels}
Recall that in the weighted barycenter space, each negative
singularity $p_i$ has weight $\beta_i$, while every other
point has weight one. The corresponding set of critical
values is
\[
 \widehat\Si=c_N\{\ell+\beta_I:\ell\in\N_0,\ 
 I\subset I_-,\ \ell+|I|>0\}\subset\Si_H.
\]
In particular, $\rho\notin\Si_H$ implies
$\rho\notin\widehat\Si$.
We use $\widehat\Si$ to study concentration on low energy
sublevels and $\Si_H$ in \eqref{eq:sigma} to study
compactness of solutions.

Let $Y$ be a compact subset of $\R^N$. For
$\mu,\nu\in\Prob(Y)$, define the bounded--Lipschitz distance
\[
 d_{\mathrm{BL}}(\mu,\nu)
 =\sup_{\substack{\varphi:Y\to\R\text{ Lipschitz}\\
 \norm{\varphi}_{L^\infty(Y)}+\operatorname{Lip}(\varphi)\le1}}
 \Big|\int_Y\varphi\,\d(\mu-\nu)\Big|.
\]
The distance $d_{\mathrm{BL}}$ induces the weak topology
on $\Prob(Y)$. For $u\in E$, define the probability measure
\[
 \mu_u=\frac{He^u}{M(u)}\,\d x
 \quad\text{ on }\, \overline\Om.
\]

When $\mu_u$ has a fixed amount of mass in several
disjoint regions, we obtain a better estimate for $\ln M(u)$.

\begin{lem}\label{lem:improved}
Let $A_j\subset\overline\Om$, $1\le j\le a$, be Borel sets with
pairwise disjoint relatively open neighborhoods $U_j$. Choose
fixed Lipschitz functions $\zeta_j$ satisfying
\[
 0\le\zeta_j\le1,\quad
 \zeta_j=1\text{ on }\,A_j,\quad
 \supp\zeta_j\subset U_j.
\]
Assume that each $\overline U_j$ contains at most one
negative singularity. Define
\[
 w_j=\left\{\begin{aligned}
 &1 &&\text{ if }\, \overline U_j\cap P_-=\varnothing,\\
 &\beta_i &&\text{ if }\, \overline U_j\cap P_-=\{p_i\},
 \end{aligned}\right.
 \quad W=\sum_{j=1}^a w_j.
\]
If $\mu_u(A_j)\ge\eta>0$ for any  $j=1,\ldots,a$,
then for any $\var>0$, there exists a constant $C>0$,
independent of $u$, such that
\be\label{eq:improved}
 \ln M(u)\le\Big(\frac1{c_NW}+\var\Big)Q(u)+C.
\ee
The estimate also holds when $U_j$ meets $\pa\Om$
or contains a positive singularity.
\end{lem}
\begin{proof}
We follow the argument in \cite[Lemma~4.3]{B15}.
Let $U$ be a neighborhood of weight $w$. We first show
that, for any $z\in E$ and $\delta>0$,
\be\label{eq:localMT}
 \ln\int_{U\cap\Om}He^z\,\d x
 \le\Big(\frac1{c_Nw}+\delta\Big)Q(z)+C_\delta.
\ee
If $\overline U\cap P_-=\varnothing$, then $H$ is
bounded on $U\cap\Om$, and \eqref{eq:localMT} follows
directly from \eqref{eq:MT} with $w=1$.
If $\overline U\cap P_-=\{p_i\}$, write
$H=V|x-p_i|^{N\al_i}$ with $V$ bounded above.
Choose $1<b<-1/\al_i$ and set $b'=b/(b-1)$.
By H\"older's inequality and \eqref{eq:MT} applied
to $b'z$, we have
\[
 \ln\int_{U\cap\Om}He^z\,\d x
 \le\frac{(b')^{N-1}}{c_N}Q(z)+C_b.
\]
Since
\[
 \frac{(b')^{N-1}}{c_N}\to\frac1{c_N\beta_i}
 \quad\text{ as }\, b\uparrow-1/\al_i,
\]
we obtain \eqref{eq:localMT} with $w=\beta_i$.
The same argument applies when $U$ meets $\pa\Om$
or contains a positive singularity.

We now apply \eqref{eq:localMT} on each $U_j$.
Set $R_u=\norm{u}_E$. If $R_u=0$, the conclusion is
immediate. For $R_u>0$, choose $A>0$ and set
$T=AR_u$ and $v=(u-T)^+$. By the Sobolev embedding,
\[
 \norm{v}_{L^N(\Om)}^N
 \le T^{-N}\norm{u}_{L^{2N}(\Om)}^{2N}
 \le CA^{-N}R_u^N.
\]
Since the supports of $\zeta_j$ are pairwise disjoint,
for any $\delta_1>0$, we have
\[
 \sum_{j=1}^a Q(\zeta_jv)
 \le(1+\delta_1)Q(u)
    +C_{\delta_1}A^{-N}Q(u).
\]

By $\mu_u(A_j)\ge\eta$, $\zeta_j=1$ on $A_j$,
and $u\le T+v$, we have
\[
 \eta M(u)
 \le\int_{A_j\cap\Om}He^u\,\d x
 \le e^T\int_{U_j\cap\Om}He^{\zeta_jv}\,\d x.
\]
Applying \eqref{eq:localMT} with $z=\zeta_jv$ yields
\[
 w_j\ln M(u)
 \le w_jAR_u+
 \Big(\frac1{c_N}+\delta w_j\Big)Q(\zeta_jv)+C.
\]
Summing over $j$ and using $w_j\le1$, we obtain
\[
 W\ln M(u)
 \le WAR_u+
 \Big(\frac1{c_N}+\delta\Big)
 ( 1+\delta_1+C_{\delta_1}A^{-N})Q(u)+C.
\]
First choose $\delta$ and $\delta_1$ sufficiently small
and then choose $A$ sufficiently large.
By Young's inequality, for any $\delta_2>0$,
\[
 AR_u\le\delta_2Q(u)+C_{\delta_2,A}.
\]
Dividing by $W$ and choosing $\delta_2$ sufficiently
small, we obtain \eqref{eq:improved}.
\end{proof}

By \eqref{eq:improved}, $J_\rho(u)$ is bounded below
if $\mu_u$ has a fixed positive amount of mass in each
of several separated regions with total weight greater
than $r$. Following the argument in
\cite[Theorem~4.1]{B15}, we obtain the following
concentration result.
\begin{lem}\label{lem:concentration}
Let $\rho>0$ and $\rho\notin\widehat\Si$.
If $\Xr(\overline\Om)\ne\varnothing$, then
\be\label{eq:concentration}
 \sup_{J_\rho(u)\le-L}
 \dist_{\mathrm{BL}}(\mu_u,\Xr(\overline\Om))
 \to0\quad\text{ as }\, L\to\infty.
\ee
If $\Xr(\overline\Om)=\varnothing$, then $r<\beta_*$
and $J_\rho$ is coercive.
\end{lem}
\begin{proof}
We first show that $\Xr(\overline\Om)$ is compact.
Since $\beta_i<1$ for $i\in I_-$, merging support points
or letting some masses tend to zero cannot increase
the total weight. Moreover, every measure in
$\Xr(\overline\Om)$ has at most $r/\beta_*$ support points.
Hence $\Xr(\overline\Om)$ is a finite union of continuous
images of compact parameter spaces, and is therefore compact.

Fix $\var_0>0$ and set
\[
 \mathcal A_{\var_0}
 =\{\mu\in\Prob(\overline\Om):
 \dist_{\mathrm{BL}}(\mu,\Xr(\overline\Om))\ge\var_0\}.
\]
By compactness of $\Prob(\overline\Om)$,
the set $\mathcal A_{\var_0}$ is compact.

For any $\mu\in\mathcal A_{\var_0}$, we can choose
finitely many distinct points $x_1,\ldots,x_a\in\supp\mu$
with total weight $W>r$. If $\supp\mu$ is infinite,
we simply choose sufficiently many points.
If $\supp\mu$ is finite, the conclusion follows from
$\mu\notin\Xr(\overline\Om)$ and $\rho\notin\widehat\Si$.

Around each $x_j$, choose relatively open neighborhoods
$A_j$ and $U_j$ with $\overline A_j\subset U_j$ such that
the closures $\overline U_j$ are pairwise disjoint and
contain no negative singularities other than $x_j$.
Let $w_j$ be the corresponding weight as in
Lemma~\ref{lem:improved}, so that
$W=\sum_{j=1}^a w_j>r$.
Since $x_j\in\supp\mu$, we can choose $\eta>0$ such that
$\mu(A_j)>2\eta$ for any  $j$.

By the Portmanteau theorem,
\[
 \mathcal U_\mu
 =\{\nu\in\Prob(\overline\Om):\nu(A_j)>\eta
 \text{ for any  }\,j\}
\]
is an open neighborhood of $\mu$ in the weak topology.
By compactness, finitely many such neighborhoods cover
$\mathcal A_{\var_0}$.

For each of these neighborhoods, choose $\var>0$
sufficiently small such that
\[
 \rho\Big(\frac1{c_NW}+\var\Big)<1.
\]
By Lemma~\ref{lem:improved}, whenever
$\mu_u\in\mathcal A_{\var_0}$, we obtain
\[
 J_\rho(u)
 \ge\Big[1-\rho\Big(\frac1{c_NW}+\var\Big)\Big]Q(u)-C
 \ge-C.
\]
Since only finitely many neighborhoods are involved,
the constant $C$ can be chosen independently of $u$.
Therefore, for sufficiently large $L$,
\[
 J_\rho(u)\le-L
 \quad\text{ implies that }\quad
 \dist_{\mathrm{BL}}(\mu_u,\Xr(\overline\Om))<\var_0.
\]
Since $\var_0>0$ is arbitrary, we obtain
\eqref{eq:concentration}.

Finally, if $\Xr(\overline\Om)=\varnothing$, then
$r<\beta_*$ since $\rho\notin\widehat\Si$.
By the weighted exponential estimate following
\eqref{eq:MT}, for any $\var>0$,
\[
 J_\rho(u)\ge
 \Big(1-\frac{\rho}{c_N\beta_*}-\rho\var\Big)Q(u)-C_\var.
\]
Since $\rho<c_N\beta_*$, choosing $\var>0$
sufficiently small proves the coercivity of $J_\rho$.
\end{proof}

\subsection{Barycenter projections and trial functions}
To construct trial functions for weighted barycenters,
we first move their support into the interior of $\Om$.
Following \cite[Lemma~2.5]{B15}, we keep the negative
singularities fixed and choose a deformation retract
that contains no positive singularities.

\begin{lem}\label{lem:spine}
There exists a compact set $\Ga$, homeomorphic to a
finite polyhedron, such that
\[
 P_-\subset\Ga\Subset\Om\setminus P_+.
\]
Moreover, $\overline\Om$ admits a strong deformation
retraction onto $\Ga$. The deformation fixes $P_-$
and does not increase the total weight of any barycenter.
Consequently,
\be\label{eq:spinehomotopy}
 \Xr(\Ga)\simeq\Xr(\overline\Om).
\ee
\end{lem}
\begin{proof}
We first construct a compact subset of $\Om$ that is
a strong deformation retract of $\overline\Om$.
By moving the boundary of $\Om$ slightly inward,
we obtain a compact smooth manifold $M\Subset\Om$
diffeomorphic to $\overline\Om$.

Choose a finite triangulation of $M$. Consider the dual
graph of the $N$-simplices, add a vertex representing
the exterior of $M$, and choose a spanning tree rooted
at this vertex. Starting from the boundary, we remove
the $N$-simplices one by one along the spanning tree.
At each step, the face shared with the parent simplex
is free, since the parent simplex has already been
removed. Each elementary collapse defines a strong
deformation retraction. After removing all $N$-simplices,
we obtain a finite polyhedron $\Ga_0$ of dimension
at most $N-1$.

We next show that $\Ga_0$ contains sufficiently many
distinct points. Let $T$ be the number of $N$-simplices
and $b$ the number of boundary faces of dimension $N-1$
in the original triangulation. Since each collapse
removes one $N$-simplex and one $(N-1)$-dimensional face,
the number of remaining $(N-1)$-dimensional faces is
\[
 \frac{(N+1)T+b}{2}-T
 =\frac{(N-1)T+b}{2}>0.
\]
Since $N\ge2$, $\Ga_0$ contains infinitely many points
and has empty interior in $\Om$.

Choose distinct points
\[
 q_i\in\Ga_0\quad\text{ for }\, i\in I_-,
 \quad
 q_i\in\Om\setminus\Ga_0\quad\text{ for }\, i\in I_+,
\]
with $q_i\notin P$ for any  $i$.
Since $\Om$ is connected and $N\ge2$, we can move
the points $q_i$ to $p_i$ along paths without
intersecting the other prescribed points.
Thus there exists a smooth diffeomorphism
$F:\Om\to\Om$, equal to the identity near $\pa\Om$,
such that $F(q_i)=p_i$ for any  $i$.
Set $\Ga=F(\Ga_0)$. Then
\[
 P_-\subset\Ga\Subset\Om\setminus P_+.
\]
The deformation retraction onto $\Ga_0$, together
with $F$, defines a strong deformation retraction
$\Pi_t$ of $\overline\Om$ onto $\Ga$, with
$\Pi_0=\Id$ and $\Pi_1=\Pi$.

Finally, for any $\si=\sum_jt_j\delta_{x_j}
\in\Xr(\overline\Om)$, define
\[
 (\Pi_t)_*\si=\sum_jt_j\delta_{\Pi_t(x_j)}.
\]
Since $\Pi_t$ fixes every point of $P_-$ and the weight
of every other point is at most one, the total weight
cannot increase during the deformation.
Merging support points also cannot increase the total
weight. Hence $(\Pi_t)_*\si\in\Xr(\overline\Om)$
for any  $t\in[0,1]$.
Therefore, $\Xr(\Ga)$ is a strong deformation retract
of $\Xr(\overline\Om)$, which proves
\eqref{eq:spinehomotopy}.
\end{proof}

To construct the projection, we first show that the
weighted barycenter space admits a continuous retraction
from a neighborhood in $\Prob(\Ga)$.
\begin{lem}\label{lem:polyhedron}
For $r>0$, $\Xr(\Ga)$ is homeomorphic to a compact
finite polyhedron. Moreover, there exist an open
neighborhood $\mathcal U_r$ of $\Xr(\Ga)$ in $\Prob(\Ga)$
and a continuous map
\[
 R_r:\mathcal U_r\to\Xr(\Ga)
\]
such that $R_r(\si)=\si$ for any  $\si\in\Xr(\Ga)$.
\end{lem}

\begin{proof}
If $\Xr(\Ga)=\varnothing$, the conclusion is immediate.
Otherwise, we identify $\Ga$ with a finite polyhedron
in $\R^d$, with all negative singularities as vertices.
Since every measure in $\Xr(\Ga)$ has at most $r/\beta_*$
support points, we can choose an integer $n_0\ge1$
such that each support contains at most $n_0$ points.

For $\mu\in\Prob(\Ga)$, define the moment map
\[
 \mathcal M(\mu)=
 \Big(\int_\Ga x^\gamma\,\d\mu\Big)_{|\gamma|\le4n_0-2}.
\]
We first show that $\mathcal M$ is injective on $\Xr(\Ga)$.
Indeed, for any two distinct measures
$\mu,\nu\in\Xr(\Ga)$, the union of their supports
contains at most $2n_0$ points. Choose a point $x_0$
at which the two measures have different masses.
The polynomial
\[
 P(x)=\prod_{\substack{y\in\supp\mu\cup\supp\nu\\y\ne x_0}}
 |x-y|^2
\]
has degree at most $4n_0-2$, vanishes at all other
support points, and satisfies $P(x_0)>0$. Hence
\[
 \int_\Ga P\,\d\mu\ne\int_\Ga P\,\d\nu,
\]
which implies $\mathcal M(\mu)\ne\mathcal M(\nu)$.

By the definition of $\Xr(\Ga)$, every barycenter
can be represented by a point in
\[
 \Ga^k\times\Delta^{k+|I|-1},\quad
 k\in\N_0,\quad I\subset I_-,\quad
 k+\beta_I<r,\quad k+|I|>0,
\]
where $k$ points vary in $\Ga$ and the remaining
points are fixed at $p_i$ with $i\in I$.
If some support points coincide or some masses vanish,
the total weight cannot increase.

Since there are only finitely many choices of $k$ and $I$,
$\mathcal M(\Xr(\Ga))$ is a compact semialgebraic set
and therefore admits a finite triangulation.
Since $\Xr(\Ga)$ is compact and $\mathcal M$ is
continuous and injective on $\Xr(\Ga)$, the map
$\mathcal M$ is a homeomorphism onto its image.
Thus $\Xr(\Ga)$ is homeomorphic to a compact
finite polyhedron.

We next construct the retraction. Set
$K=\mathcal M(\Xr(\Ga))$ and choose a closed Euclidean
ball containing $\mathcal M(\Prob(\Ga))$.
Triangulate the ball so that $K$ is a subcomplex.
After a barycentric subdivision, $K$ is a full
subcomplex. On an open neighborhood $\mathcal V$ of $K$,
we keep only the barycentric coordinates corresponding
to vertices in $K$ and divide them by their sum.
This defines a continuous retraction
$R:\mathcal V\to K$.

Finally, set
\[
 \mathcal U_r=\mathcal M^{-1}(\mathcal V),\quad
 R_r=\mathcal M^{-1}\circ R\circ\mathcal M:
 \mathcal U_r\to\Xr(\Ga).
\]
By construction, $R_r(\si)=\si$ for any 
$\si\in\Xr(\Ga)$. Hence $R_r$ is the required
neighborhood retraction.
\end{proof}

Let $\mathcal U_r$ and $R_r:\mathcal U_r\to\Xr(\Ga)$
be the neighborhood and retraction given in
Lemma~\ref{lem:polyhedron}, respectively.
We now use $R_r$ to construct a projection from
low energy sublevels to $\Xr(\Ga)$.

\begin{prop}\label{prop:projection}
Let $\rho>0$, $\rho\notin\widehat\Si$, and
$\Xr(\Ga)\ne\varnothing$. Then, for all sufficiently
large $L$, the map
\be\label{eq:Psi}
 \Psi:J_\rho^{-L}\to\Xr(\Ga),\quad
 \Psi(u)=R_r(\Pi_*\mu_u)
\ee
is well defined and continuous.
\end{prop}
\begin{proof}
By Lemma~\ref{lem:spine}, the map
$\Pi_*:\Prob(\overline\Om)\to\Prob(\Ga)$ is uniformly
continuous and maps $\Xr(\overline\Om)$ into $\Xr(\Ga)$.
By \eqref{eq:concentration}, for sufficiently large $L$,
$\Pi_*\mu_u\in\mathcal U_r$ whenever $J_\rho(u)\le-L$.
Since both $u\mapsto\Pi_*\mu_u$ and $R_r$ are continuous,
the map $\Psi$ in \eqref{eq:Psi} is well defined
and continuous.
\end{proof}

The projection in \eqref{eq:Psi} maps low energy
sublevels to the weighted barycenter space.
We now construct a map in the opposite direction
using trial functions supported near $\Ga$.
Since $\Ga\cap P_+=\varnothing$, the factors with
positive exponents can be included in the positive
background weight near $\Ga$.

\begin{lem}\label{lem:bubbles}
Let $\rho>0$, $\rho\notin\widehat\Si$, and
$\Xr(\Ga)\ne\varnothing$. For sufficiently large $\lam$,
there exists a continuous map
$\Phi_\lam:\Xr(\Ga)\to E$ such that
\begin{align}
 M(\Phi_\lam(\si))
 &=M(0)+\lam^q,\label{eq:exactmass}\\
 Q(\Phi_\lam(\si))
 &\le qc_N\Big(\#(\supp\si\setminus P_-)
 +\sum_{\substack{i\in I_-\\p_i\in\supp\si}}\beta_i\Big)
 \ln\lam+C,\label{eq:bubbleenergy}\\
 d_{\mathrm{BL}}(\mu_{\Phi_\lam(\si)},\si)
 &\to0\quad\text{ as }\,\lam\to\infty.\label{eq:bubblemeasure}
\end{align}
The estimates and convergence are uniform in
$\si\in\Xr(\Ga)$.
\end{lem}
\begin{proof}
Choose $\kappa\in C_c^\infty(\Om\setminus P_+)$
such that $\kappa=1$ on a neighborhood of $\Ga$.
For each $i\in I_+$, replace the factor
$|x-p_i|^{N\al_i}$ in $H$ by
\[
 \exp\big(N\al_i\kappa(x)\ln|x-p_i|) ,
\]
which is equal to one near $p_i$.
Denote the resulting weight by $\widetilde H$.
Then $\widetilde H=H$ near $\Ga$, and $\widetilde H$
has only negative singularities.
The corresponding resonance set and weighted barycenter
space are $\widehat\Si$ and $\Xr(\Ga)$, respectively.

We apply \cite[Proposition~4.1 and formula~(4.4)]{ATW26}
with weight $\widetilde H$ and compact set $K=\Ga$.
Choose the radius of the trial functions sufficiently
small so that $\Phi_\lam(\si)=0$ outside the neighborhood
where $H=\widetilde H$. Hence
\[
 H(e^{\Phi_\lam(\si)}-1)
 =\widetilde H(e^{\Phi_\lam(\si)}-1).
\]
By \cite[formula~(4.4)]{ATW26}, we obtain
\[
 \int_\Om H(e^{\Phi_\lam(\si)}-1)\,\d x
 =\int_\Om\widetilde H(e^{\Phi_\lam(\si)}-1)\,\d x
 =\lam^q.
\]
Thus \eqref{eq:exactmass} holds.
Combining the energy estimate in
\cite[Proposition~4.1]{ATW26} with the mass identity
above, we obtain \eqref{eq:bubbleenergy}.

For $u\in E$, set
\[
 \widetilde\mu_u
 =\frac{\widetilde H e^u}
 {\int_\Om\widetilde H e^u\,\d x}\,\d x.
\]
The mass identity above implies
\[
 \norm{\mu_{\Phi_\lam(\si)}
 -\widetilde\mu_{\Phi_\lam(\si)}}_{L^1(\Om)}
 \le\frac{2\norm{H-\widetilde H}_{L^1(\Om)}}
 {M(0)+\lam^q}\to0
\]
uniformly in $\si\in\Xr(\Ga)$.
Since
$d_{\mathrm{BL}}(\widetilde\mu_{\Phi_\lam(\si)},\si)\to0$
uniformly by \cite[Proposition~4.1]{ATW26},
we obtain \eqref{eq:bubblemeasure}.

Finally, every barycenter in $\Xr(\Ga)$ has at most
$r/\beta_*$ support points. Hence there are only finitely
many possible values of the total weight, all strictly
less than $r$. By \eqref{eq:exactmass} and
\eqref{eq:bubbleenergy}, there exists $\var>0$,
independent of $\si$, such that
\be\label{eq:negativeestimate}
 J_\rho(\Phi_\lam(\si))\le-\var\ln\lam+C.
\ee
\end{proof}

\begin{rem}
(1) The coercive range of $J_\rho$ is sharp on each fixed
domain. Fix an auxiliary parameter
$\bar\rho>c_N\beta_*$ with $\bar\rho\notin\widehat\Si$,
and apply Lemma~\ref{lem:bubbles} at $\bar\rho$ to a
single Dirac measure $\si$ of weight $\beta_*$. Keep
the resulting family $\Phi_\lam(\si)$ fixed while
varying the parameter $\rho$ in the energy.
By \eqref{eq:exactmass},
\[
 M(\Phi_\lam(\si))=M(0)+\lam^q\to\infty
 \quad\text{ as }\,\lam\to\infty.
\]
The weighted exponential estimate following \eqref{eq:MT}
implies $\norm{\Phi_\lam(\si)}_E\to\infty$.
Moreover, by \eqref{eq:exactmass} and
\eqref{eq:bubbleenergy}, we have
\[
 J_\rho(\Phi_\lam(\si))
 \le q(c_N\beta_*-\rho)\ln\lam+C.
\]
Hence $J_\rho$ is not coercive for $\rho=c_N\beta_*$ and
is unbounded below for $\rho>c_N\beta_*$.

(2) We next explain why positive singularities have weight
one in the weighted barycenter space.
Fix $p_i\in P_+$ and $1<r<\beta_i$.
Choose a sequence $y_n\in\Om\setminus P$ with $y_n\to p_i$.
For each $n$, we apply
\cite[Proposition~4.1 and formula~(4.4)]{ATW26}
with $K=\{y_n\}$, a noninteger normalized parameter
in $(1,r)$, and a positive $C^1$ weight equal to $H$
near $y_n$. Such a weight is obtained by interpolating
with the constant~$1$ using a cutoff in a ball
about $y_n$ disjoint from $P$. As in the proof of
Lemma~\ref{lem:bubbles},
we choose trial functions $u_{n,\lam}$ supported
in a small ball around $y_n$ such that
\[
 M(u_{n,\lam})=M(0)+\lam^q,\quad 
 Q(u_{n,\lam})\le qc_N\ln\lam+C_n,\quad 
 d_{\mathrm{BL}}(\mu_{u_{n,\lam}},\delta_{y_n})
 \to0\quad\text{ as }\, \lam\to\infty.
\]
Since $\rho=c_Nr>c_N$, for each $n$ we can choose
$\lam_n$ sufficiently large such that
\[
 J_\rho(u_{n,\lam_n})\le-n,\quad
 d_{\mathrm{BL}}(\mu_{u_{n,\lam_n}},\delta_{y_n})
 \le\frac1n.
\]
Since $y_n\to p_i$, we obtain
\[
 J_\rho(u_{n,\lam_n})\to-\infty,\quad
 \mu_{u_{n,\lam_n}}\weak\delta_{p_i}.
\]
However, $\delta_{p_i}$ would not satisfy the weight
condition if $p_i$ had weight $\beta_i>r$.
This explains why we use weight one at positive
singularities in the weighted barycenter space.
\end{rem}

By \eqref{eq:negativeestimate}, the trial functions
$\Phi_\lam(\si)$ lie in $J_\rho^{-L}$ for sufficiently
large $\lam$. Using \eqref{eq:bubblemeasure}, we now
show that the composition $\Psi\circ\Phi_\lam$ is
homotopic to the identity on $\Xr(\Ga)$.

\begin{prop}\label{prop:homotopy}
Let $\rho>0$, $\rho\notin\widehat\Si$, and
$\Xr(\Ga)\ne\varnothing$. For sufficiently large $L$
and then sufficiently large $\lam$, we have
\[
 \Phi_\lam:\Xr(\Ga)\to J_\rho^{-L},\quad
 \Psi\circ\Phi_\lam\simeq\Id_{\Xr(\Ga)}.
\]
For any field $\mathbb F$ and integer $j\ge0$,
the induced homology maps satisfy
\be\label{eq:splithomology}
 \Psi_*\circ(\Phi_\lam)_*
 =\Id_{\widetilde H_j(\Xr(\Ga);\mathbb F)}.
\ee
In particular, $\widetilde H_j(\Xr(\overline\Om);\mathbb F)$
is isomorphic to a direct summand of
$\widetilde H_j(J_\rho^{-L};\mathbb F)$.

If $\Xr(\overline\Om)$ is noncontractible, then
$\Phi_\lam$ is not homotopic to a constant map
in $J_\rho^{-L}$.
\end{prop}

\begin{proof}
By Proposition~\ref{prop:projection}, $\Psi$ is well
defined for sufficiently large $L$. By
\eqref{eq:negativeestimate}, we have
$\Phi_\lam(\Xr(\Ga))\subset J_\rho^{-L}$ for
sufficiently large $\lam$.

Since $\Pi$ is the identity on $\Ga$,
\eqref{eq:bubblemeasure} implies
\[
 \Pi_*\mu_{\Phi_\lam(\si)}\to\si
 \quad\text{uniformly in }\si\in\Xr(\Ga).
\]
By compactness of $\Xr(\Ga)$, for sufficiently large
$\lam$, we have
\[
 (1-t)\si+t\Pi_*\mu_{\Phi_\lam(\si)}\in\mathcal U_r,
 \quad t\in[0,1].
\]
Hence the map
\[
 \mathcal H(t,\si)
 =R_r\big((1-t)\si+t\Pi_*\mu_{\Phi_\lam(\si)}) 
\]
is continuous and satisfies
\[
 \mathcal H(0,\si)=\si,\quad
 \mathcal H(1,\si)=\Psi\circ\Phi_\lam(\si).
\]
Therefore, $\Psi\circ\Phi_\lam\simeq\Id_{\Xr(\Ga)}$.
Passing to the induced homology maps proves
\eqref{eq:splithomology}; see also
\cite[Corollary~5.4]{B15}.

If $\Phi_\lam$ is homotopic to a constant map in
$J_\rho^{-L}$, then $\Psi\circ\Phi_\lam$ is also
homotopic to a constant map. Hence $\Xr(\Ga)$ is
contractible. By \eqref{eq:spinehomotopy},
$\Xr(\overline\Om)$ is also contractible, which
contradicts the assumption.
\end{proof}

\begin{rem}
The homology estimates in \eqref{eq:splithomology} do not
determine the number of solutions. Battaglia used weak Morse inequalities
and generic nondegeneracy to obtain multiplicity
\cite[Theorem~1.2, pp.~52--53]{B15}.
See also \cite[Lemma~2.14 and Theorem~2.15, p.~62]{B15}.
For $N>2$, although the $N$-energy is of class $C^2$,
its second derivative degenerates at points where
$\nabla u=0$. To obtain similar multiplicity results,
we need further analysis of the Fredholm property
and generic nondegeneracy. We do not study these
questions or assume differentiability of
$(-\Delta_N)^{-1}$ in the present paper.
\end{rem}

\subsection{Barycenter topology and existence criteria}

A nonempty barycenter space is noncontractible if its
Euler characteristic is different from one.
We compute the Euler characteristic by considering
support points in $P_-$ and outside $P_-$ separately.

\begin{prop}\label{prop:Euler}
Let $r>0$.
\be\label{eq:Euleridentity}
 1-\chi(\Xr(\Ga))
 =\sum_{\substack{I\subset I_-,\ k\in\N_0\\
                  \beta_I+k<r}}
 (-1)^{|I|+k}\binom{\chi(\Om)-|I_-|}{k}.
\ee
We use the convention $\chi(\varnothing)=0$.
\end{prop}

\begin{proof}
Since $\Ga$ is compact and homotopy equivalent to
$\overline\Om$, we have
\[
 \chi_c(\Ga\setminus P_-)=\chi(\Om)-|I_-|,
\]
where $\chi_c$ denotes the Euler characteristic with
compact support.

We divide $\Xr(\Ga)$ into strata according to the number
$k$ of support points outside $P_-$ and the set
$I\subset I_-$ of negative singularities in the support.
By the definition of $\Xr(\Ga)$, each stratum satisfies
\[
 k+\beta_I<r,\quad k+|I|>0.
\]
The space of $k$ distinct unordered points in
$\Ga\setminus P_-$ has Euler characteristic with
compact support
\[
 \binom{\chi(\Om)-|I_-|}{k}.
\]
The formula follows from the inclusion--exclusion
principle and division by $k!$, since the symmetric
group acts freely on ordered configurations of
distinct points.

Over ordered configurations of the $k$ free points, the masses
vary in an open simplex of dimension $k+|I|-1$. Permuting the
free points and their masses acts freely on this product, and
the quotient is the stratum. The covering has $k!$ sheets.
The product formula for compactly supported Euler characteristics
and division by $k!$ give the stratum contribution
\[
 (-1)^{k+|I|-1}\binom{\chi(\Om)-|I_-|}{k}
\]
to the Euler characteristic with compact support of
$\Xr(\Ga)$.

Since the strata are semialgebraic in a polyhedral
realization, we can sum their Euler characteristics
with compact support. The entire space is compact for any  $r>0$.
Indeed, the strict weight inequality selects finitely many pairs $(k,I)$.
For each such pair, take the measures parametrized by $\Ga^k$ and the
closed simplex of masses at the corresponding locations, allowing
coincident points and zero masses. The weight cannot increase when
points merge, when a free point reaches $P_-$, or when a mass vanishes.
These compact images form $\Xr(\Ga)$. Hence
$\chi_c(\Xr(\Ga))=\chi(\Xr(\Ga))$. Therefore,
\[
 \chi(\Xr(\Ga))
 =\sum_{\substack{I\subset I_-,\ k\in\N_0\\
                  \beta_I+k<r,\ k+|I|>0}}
 (-1)^{k+|I|-1}\binom{\chi(\Om)-|I_-|}{k}.
\]
Including the term $(k,I)=(0,\varnothing)$ in the sum
for $1-\chi(\Xr(\Ga))$ proves \eqref{eq:Euleridentity}.
See also \cite[Theorem~1.1 and Corollary~1.2]{K19}
for the Euler characteristic of weighted barycenter
spaces.
\end{proof}

If $r<1$, only negative singularities can occur in the
support of a weighted barycenter. In this case, the
weighted barycenter space is a finite simplicial complex,
whose homotopy type can be determined explicitly.

\begin{prop}\label{thm:threshold-wedge}
Suppose there are $s\ge1$ negative singularities with
weights
\[
 0<\beta_1\le\cdots\le\beta_s<1.
\]
Let $0<r<1$.
If $r\le\beta_1$, then $\Xr(\overline\Om)=\varnothing$.
If $r>\beta_1$, define
\[
 \mathcal C_r=\{\varnothing\ne I\subset\{2,\ldots,s\}:
 \beta_I<r\le\beta_I+\beta_1\}.
\]
Then
\[
 \Xr(\overline\Om)\simeq
\bigvee_{I\in\mathcal C_r}\S^{|I|-1}.
\]
Here the wedge sum is a single point if
$\mathcal C_r=\varnothing$, and each zero-dimensional
sphere adds one connected component.

Consequently, for any field $\mathbb F$ and integer
$j\ge0$, we have
\[
 \dim_{\mathbb F}\widetilde H_j
 (\Xr(\overline\Om);\mathbb F)
 =\#\{I\in\mathcal C_r:|I|=j+1\}.
\]
In particular, for $r>\beta_1$, the weighted barycenter
space $\Xr(\overline\Om)$ is noncontractible if and only
if $\mathcal C_r\ne\varnothing$.
\end{prop}

\begin{proof}
We follow the argument in \cite{PW13}.
Since $r<1$, only negative singularities can occur
in the support of a weighted barycenter.
Hence $\Xr(\overline\Om)$ is the geometric realization
of the finite simplicial complex
\[
 K_r=\{I\subset\{1,\ldots,s\}:\beta_I<r\}.
\]
If $r\le\beta_1$, then $\Xr(\overline\Om)=\varnothing$.

Suppose now that $r>\beta_1$. Consider the closed star
of vertex $1$, defined by
\[
 \operatorname{st}(1)
 =\{I\in K_r:I\cup\{1\}\in K_r\}.
\]
Since $\operatorname{st}(1)$ is a cone with vertex $1$,
it is contractible.

By the definition of $\mathcal C_r$, a simplex
$I\in K_r$ lies outside $\operatorname{st}(1)$
if and only if $I\in\mathcal C_r$.
For each $I\in\mathcal C_r$ and $j\in I$, we have
\[
 \beta_{I\setminus\{j\}}+\beta_1
 \le\beta_I<r,
\]
since $\beta_j\ge\beta_1$. Hence every proper face
of $I$ lies in $\operatorname{st}(1)$.
Moreover, adding any new vertex to $I$ increases
the total weight to at least $\beta_I+\beta_1\ge r$.
Thus $I$ is a maximal simplex of $K_r$.

Therefore, $K_r$ is obtained from $\operatorname{st}(1)$
by attaching one $(|I|-1)$-dimensional simplex along
its entire boundary for each $I\in\mathcal C_r$.
Since $\operatorname{st}(1)$ is a contractible subcomplex,
we can collapse it to a single point without changing
the homotopy type of $K_r$. Each attached simplex
then becomes a sphere of dimension $|I|-1$.
For $|I|=1$, the additional vertex together with
the base point forms $\S^0$. Consequently,
\[
 \Xr(\overline\Om)\simeq
\bigvee_{I\in\mathcal C_r}\S^{|I|-1}.
\]
The homology formula and the noncontractibility
criterion follow directly from the wedge sum above.
\end{proof}

To prove existence results, we use a deformation argument with
nearby parameters, following 
\cite[Propositions~3.3--3.4]{ATW26}.
\begin{lem}\label{lem:descent}
Let $\rho>0$ and $\rho\notin\Si_H$. Suppose that
$J_\rho$ has no critical point. Then, for any $a<b$
with $J_\rho^a\ne\varnothing$, the sublevel $J_\rho^a$
is a strong deformation retract of $J_\rho^b$.
\end{lem}

\begin{proof}
Choose $\delta_0\in(0,\rho)$ such that
$[\rho-\delta_0,\rho]\cap\Si_H=\varnothing$.
Since $J_\rho$ has no critical point, Theorem~\ref{thm:compactness}
implies that there exists $\delta\in(0,\delta_0]$ such that
\[
 J_\lam'(u)\ne0,\quad
 \forall\,u\in E,\quad \lam\in[\rho-\delta,\rho].
\]
Otherwise, a sequence of critical points of $J_{\lam_j}$
with $\lam_j\to\rho$ would converge to a critical point
of $J_\rho$.

For each $u\in E$, set
\[
 m(u)=\min_{\lam\in[\rho-\delta,\rho]}
 \norm{J_\lam'(u)}_{E^*}>0.
\]
Since $J_\lam'(u)$ is affine in $\lam$, the set
\[
 \{J_\lam'(u):\lam\in[\rho-\delta,\rho]\}
\]
is compact and convex and does not contain zero.
By the Hahn--Banach separation theorem and the
reflexivity of $E$, we can choose a vector $w\in E$
with $\norm w_E\le1$ such that
\[
 \langle J_\lam'(u),w\rangle\le-\frac12m(u),
 \quad \forall\,\lam\in[\rho-\delta,\rho].
\]
Using a locally Lipschitz partition of unity, we construct
a locally Lipschitz vector field $W:E\to E$ satisfying
\[
 \norm{W(u)}_E\le1,\quad
 \langle J_\lam'(u),W(u)\rangle\le-\frac14m(u)
\]
for all $u\in E$ and $\lam\in[\rho-\delta,\rho]$.

We claim that $m(u)$ has a positive lower bound on
every bounded subset of $E$. Otherwise, there exist
a bounded sequence $u_j\in E$ and parameters
$\lam_j\in[\rho-\delta,\rho]$ such that
$J_{\lam_j}'(u_j)\to0$ in $E^*$.
By the compactness of $(\ln M)'$ and \eqref{eq:monotonicity},
up to a subsequence, $u_j\to u$ strongly in $E$ and
$\lam_j\to\lam\in[\rho-\delta,\rho]$.
Hence $J_\lam'(u)=0$, a contradiction.

For $u_0\in J_\rho^b$ with $J_\rho(u_0)>a$, consider
the flow
\[
 \frac{\d u}{\d t}=W(u),\quad u(0)=u_0.
\]
The bound $\norm{W(u)}_E\le1$ gives global existence of the flow.
Indeed, at any finite maximal time the trajectory is Cauchy in $E$,
and the local existence theorem continues it from its limit.

Along the flow, we have
\[
 \frac{\d}{\d t}J_\lam(u(t))
 \le-\frac14m(u(t))<0,
 \quad \lam\in[\rho-\delta,\rho].
\]
In particular, both $J_\rho$ and $J_{\rho-\delta}$
decrease along the flow.
Using the identity
\[
 Q(u)=\frac{\rho}{\delta}J_{\rho-\delta}(u)
 -\Big(\frac{\rho}{\delta}-1\Big)J_\rho(u),
\]
we obtain
\[
 Q(u(t))\le
 \frac{\rho}{\delta}J_{\rho-\delta}(u_0)
 -\Big(\frac{\rho}{\delta}-1\Big)a
\]
as long as $J_\rho(u(t))\ge a$.
Hence $u(t)$ remains bounded in $E$ before reaching
the level $a$.

Since $m(u)$ has a positive lower bound on bounded
subsets of $E$, there exists $c>0$ such that
\[
 \frac{\d}{\d t}J_\rho(u(t))\le-c
\]
as long as $J_\rho(u(t))\ge a$.
Consequently, the flow reaches the level $a$ in finite
time. Denote the first time by $T(u_0)$ and set
$T(u_0)=0$ if $J_\rho(u_0)\le a$.

By continuous dependence of the flow on $u_0$ and
the strict decrease of $J_\rho$, the function $T$
is continuous for $J_\rho(u_0)>a$.
Moreover, the estimates above are locally uniform
in $u_0$. Thus
\[
 T(u_0)\le C(J_\rho(u_0)-a)\to0
 \quad\text{ as }\,J_\rho(u_0)\downarrow a
\]
locally near the level $a$.
Therefore, $T$ is continuous on $J_\rho^b$.

Let $u(t;u_0)$ denote the solution of the flow with
initial value $u_0$. Define
\[
 \mathcal H:[0,1]\times J_\rho^b\to J_\rho^b,\quad
 \mathcal H(s,u_0)=u(sT(u_0);u_0).
\]
By construction,
 \begin{gather*}
 \mathcal H(0,u_0)=u_0,\quad
 \mathcal H(1,u_0)\in J_\rho^a,\\
 \mathcal H(s,u_0)=u_0
 \quad\text{ for }\,u_0\in J_\rho^a,\ s\in[0,1].
 \end{gather*}
Hence $\mathcal H$ defines a strong deformation
retraction of $J_\rho^b$ onto $J_\rho^a$.
\end{proof}

By Theorem~\ref{thm:existence} and Proposition~\ref{prop:Euler},
we obtain the following explicit existence criterion.
\begin{cor}
Assume the conditions in \eqref{eq:assumptions} hold.
If $\rho>0$, $\rho\notin\Si_H$, and
\[
 \sum_{\substack{I\subset I_-,\ k\in\N_0\\
                  \beta_I+k<r}}
 (-1)^{|I|+k}\binom{\chi(\Om)-|I_-|}{k}\ne0,
\]
then \eqref{maineq} admits a weak solution.
\end{cor}

\begin{proof}
If $\Xr(\Ga)=\varnothing$, then $r<\beta_*$, and
existence follows from Theorem~\ref{thm:existence}(i).
If $\Xr(\Ga)\ne\varnothing$, then
$\chi(\Xr(\Ga))\ne1$ by \eqref{eq:Euleridentity}.
Hence $\Xr(\Ga)$ is noncontractible.
By \eqref{eq:spinehomotopy}, $\Xr(\overline\Om)$
is also noncontractible. Theorem~\ref{thm:existence}(ii)
then yields a solution of \eqref{maineq}.
\end{proof}

\section{Degree and domain comparison}\label{sec:4}
\subsection{Degree, entropy and potential estimates}

In Section~\ref{sec:3}, we studied the topology of low
energy sublevels to prove existence. We now use the
Euler characteristic of low energy sublevels to compute
the Leray--Schauder degree.

\begin{prop}\label{prop:all-parameter-euler}
Assume the conditions in \eqref{eq:assumptions} hold.
Let $\rho>0$ and $\rho\notin\Si_H$. Then there exists
$L_0>0$ such that, for any  $L\ge L_0$,
\[
 \sum_{j\ge0}\dim_{\mathbb Q}
 H_j(J_\rho^{-L};\mathbb Q)<\infty
\]
and
\be\label{eq:all-parameter-euler}
 d_\rho(\Om,H)=1-\chi(J_\rho^{-L}),
\ee
where
\[
 \chi(J_\rho^{-L})
 =\sum_{j\ge0}(-1)^j
 \dim_{\mathbb Q}H_j(J_\rho^{-L};\mathbb Q).
\]
We use the convention $\chi(\varnothing)=0$.
\end{prop}

\begin{proof}
Set $B=\ln M$. By the estimates in Section~\ref{sec:2},
$B\in C^1(E)$ and $B'$ maps weakly convergent sequences
in $E$ to strongly convergent sequences in $E^*$.
Moreover, $Q'$ satisfies \eqref{eq:monotonicity},
and bounded Palais--Smale sequences are compact.

Choose $0<\delta<\rho$ such that
$[\rho-\delta,\rho]\cap\Si_H=\varnothing$.
Set $J=J_\rho$ and $F=J_{\rho-\delta}$.
By Theorem~\ref{thm:compactness}, the set of critical
points of $J_\lam$ for $\lam\in[\rho-\delta,\rho]$
is compact. Hence, for sufficiently large $L$, $b$, and $C$,
every such critical point $u$ satisfies
\[
 a+3<J(u)<b-3,\quad F(u)<C-3,
 \quad a=-L.
\]
Using the identity
\[
 Q=\frac{\rho}{\delta}F
 -\Big(\frac{\rho}{\delta}-1\Big)J,
\]
we obtain
\[
 Q(u)\le\frac{\rho}{\delta}C
 -\Big(\frac{\rho}{\delta}-1\Big)a
\]
on $\{u\in E:F(u)\le C,\ a\le J(u)\le b\}$.
Thus this set is bounded in $E$.

Since $E=W^{1,N}_0(\Om)$ is isometric to a closed subspace
of $L^N(\Om;\R^N)$ through the gradient map, $E$ is separable,
reflexive, and strictly convex.
The estimates in Section~\ref{sec:2} also imply that
$Q'$ and $B'$ are uniformly bounded on bounded subsets
of $E$. Hence $Q$ and $B$ are Lipschitz on each bounded ball.

For completeness, $J_\lam'$ also has the $(S)_+$ property.
Suppose $u_j\weak u$ in $E$, $\lam_j\to\lam$, and
\[
 \limsup_{j\to\infty}
 \langle J_{\lam_j}'(u_j),u_j-u\rangle\le0.
\]
Since $B'(u_j)\to B'(u)$ in $E^*$, the pairing with
$\lam_jB'(u_j)$ tends to zero. Also,
$\langle Q'(u),u_j-u\rangle\to0$. Therefore
\eqref{eq:monotonicity} implies
\[
 0\le d_N\limsup_{j\to\infty}\norm{u_j-u}_E^N
 \le\limsup_{j\to\infty}
 \langle Q'(u_j)-Q'(u),u_j-u\rangle\le0,
\]
so $u_j\to u$ in $E$.

We apply the variational comparison in
\cite[Lemmas~3.6--3.9, Theorem~3.10 and Appendix~A]{ATW26}.
To check the compactness used in the distance-penalty step,
let $Y_j$ be increasing finite-dimensional subspaces with dense
union in $E$. For a bounded sequence $u_j\weak u$, take
$v_j\in Y_j$ with $v_j\to u$. If
$\xi_j\in\partial\dist(u_j,Y_j)$ is a convex subgradient, then
$\|\xi_j\|_{E^*}\le1$, $\xi_j|_{Y_j}=0$, and
\[
 \langle\xi_j,u_j-u\rangle
 =\dist(u_j,Y_j)+\langle\xi_j,v_j-u\rangle\ge-o(1).
\]
Thus if $J_{\lam_j}'(u_j)+s_j\xi_j\to0$ in $E^*$ for bounded
$s_j\ge0$ and $\lam_j\to\lam$, the $(S)_+$ estimate above
gives $u_j\to u$. Tests in the dense union of the $Y_j$ then
give $J_\lam'(u)=0$. If $u_j\in Y_j$ and
$\|J_{\lam_j}'(u_j)|_{Y_j}\|_{Y_j^*}\to0$, test first with
$u_j-v_j\in Y_j$ and use boundedness of $J_{\lam_j}'(u_j)$
to reach the same conclusion. In the penalized case,
$\dist(u_j,Y_j)\le\|u_j-v_j\|_E\to0$, so the modified energy
values have the same limits. Therefore neither type of
approximate critical sequence can approach the constraint levels
$a,b,C$, since all family critical points satisfy the strict
inequalities above. These are the compactness checks for the
continuation and the finite-dimensional section in Appendix~A.
Set $\mathcal M=\{F\le C\}$ and choose
$a_0<a<a_1<b$ with $a_1$ below the $J$-values of all the
family critical points. By \cite[Lemma~3.7]{ATW26}, the relative
homology with lower member $\{J<a_1\}$ agrees with that of the
corresponding pair in $\mathcal M$. In both pairs the lower
member can be replaced by the closed sublevel at $a$.
Indeed, every finite chain in $\{J<a_1\}$ lies in a sublevel
at some $a<s<a_1$. Proposition~3.3 and Lemma~3.8 of \cite{ATW26}
give the level deformation from $s$ to $a$ in the ambient
and constrained sublevels, respectively. Passing to the direct
limit and using the exact sequence of triples gives
\[
 H_*(J^b,J^a;\mathbb Q)
 \cong H_*(J^b\cap\mathcal M,J^a\cap\mathcal M;\mathbb Q).
\]

The set $U=\{a<J<b,F<C\}$ is bounded by the cap estimate and
contains every fixed point of $K_{\rho,H}$. For each $j$, let
$K_j(u)\in Y_j$ be defined by
$Q'(K_j(u))|_{Y_j}=\rho B'(u)|_{Y_j}$.
The maps $K_j$ are uniformly bounded on $\overline U$ by
testing their equations with $K_j(u)$. For $u\in\overline U$,
$Q'(K_j(u))-Q'(K_{\rho,H}(u))$ annihilates $Y_j$.
Using \eqref{eq:monotonicity} and boundedness of $Q'$ on balls,
we obtain
\[
 d_N\|K_j(u)-K_{\rho,H}(u)\|_E^N
 \le C_1\,\dist(K_{\rho,H}(u),Y_j).
\]
The right-hand side tends to zero uniformly on $\overline U$,
since $K_{\rho,H}(\overline U)$ is relatively compact and the
spaces $Y_j$ are increasing with dense union. Also
$\inf_{\partial U}\|u-K_{\rho,H}(u)\|_E>0$, since otherwise
compactness of the images would produce a boundary fixed point.
Thus the straight homotopy is admissible for large $j$, and
finite-dimensional reduction gives the Brouwer degree on
$U\cap Y_j$. On $Y_j$ the restricted gradient of $J$ and
$u-K_j(u)$ have positive pairing with $u-K_j(u)$ away from
their common zeros, by \eqref{eq:monotonicity}. This fixes
the degree sign in \cite[Lemma~3.6]{ATW26}.
Lemma~3.9 of \cite{ATW26} computes the resulting Brouwer
degree as the Euler characteristic of the constrained pair.
Excision and the pair isomorphism above yield
\[
 d_\rho(\Om,H)=\chi(J_\rho^b,J_\rho^{-L}),
\]
where
\[
 \sum_{j\ge0}\dim_{\mathbb Q}
 H_j(J_\rho^b,J_\rho^{-L};\mathbb Q)<\infty.
\]
The compactness and boundedness estimates used in this
comparison have been established under \eqref{eq:assumptions}.

Moreover, the deformation argument in
\cite{ATW26} shows that $J_\rho^b$ is contractible for
sufficiently large $b$. By the long exact sequence
in homology for the pair $(J_\rho^b,J_\rho^{-L})$,
$J_\rho^{-L}$ has finite total homology rank over $\mathbb Q$.
Consequently,
\[
 d_\rho(\Om,H)
 =\chi(J_\rho^b)-\chi(J_\rho^{-L})
 =1-\chi(J_\rho^{-L}),
\]
which proves \eqref{eq:all-parameter-euler}.

The same argument applies to a finite disjoint union
of smooth domains. In this case, we use the product
Sobolev space, the sum of the Dirichlet energies,
and a single normalizing integral over all components.
The concentration argument applies to each component;
see the discussion before Lemma~\ref{lem:small-gelfand-cutoff}.
Thus \eqref{eq:all-parameter-euler} also holds on
disconnected domains.
\end{proof}

To compare the energy sublevels on different domains,
we first introduce a functional on probability densities
whose sublevels have the same homotopy type as those
of $J_\rho$.

\begin{lem}\label{lem:nonlinear-entropy-duality}
Assume the conditions in \eqref{eq:assumptions} hold.
Choose $s_0>1$ such that $H\in L^{s_0}(\Om)$.
For $1<p<s_0$, set $p'=p/(p-1)$ and define
the space of probability densities
\[
 \mathscr P_p=\Big\{\mu\in L^p(\Om):
 \mu\ge0\ \text{a.e.},
 \int_\Om\mu\,\d x=1\Big\}
\]
with the strong $L^p$ topology.
The entropy functional
\[
 \operatorname{Ent}_H(\mu)
 =\int_\Om\mu\ln\frac{\mu}{H}\,\d x
\]
is finite and continuous on $\mathscr P_p$,
where we use the convention $0\ln0=0$.

For $\rho\ge0$, define
\[
 \mathscr S_\rho(\mu)=(-\Delta_N)^{-1}(\rho\mu),
 \quad \Theta(u)=\frac{He^u}{M(u)},
\]
and
\[
 \mathscr F_\rho(\mu)
 =\rho\operatorname{Ent}_H(\mu)
 -(N-1)Q(\mathscr S_\rho(\mu)).
\]
Equivalently,
\be\label{eq:dual-entropy-functional}
 \mathscr F_\rho(\mu)
 =\rho\operatorname{Ent}_H(\mu)
 -\frac{\rho^q}{q}\norm{\mu}_{E^*}^q,
 \quad q=\frac{N}{N-1},
\ee
where $\norm{\cdot}_{E^*}$ is the dual norm of
$\norm{\cdot}_E$.

The maps $\mathscr S_\rho:\mathscr P_p\to E$,
$\Theta:E\to\mathscr P_p$, and the functional
$\mathscr F_\rho$ are continuous.
For any $a\in\R$, define
\[
 \mathscr F_\rho^a
 =\{\mu\in\mathscr P_p:\mathscr F_\rho(\mu)\le a\}.
\]
Then $J_\rho^a$ and $\mathscr F_\rho^a$ are homotopy
equivalent through the maps
\[
 J_\rho^a
 \mathop{\rightleftarrows}^{\Theta}_{\mathscr S_\rho}
 \mathscr F_\rho^a.
\]
More precisely, the homotopies are given by
\[
 \mathcal H(t,u)
 =(1-t)u+t\mathscr S_\rho(\Theta(u)),\quad 
 \mathcal G(t,\mu)
 =(1-t)\mu+t\Theta(\mathscr S_\rho(\mu)).
\]
For any  $t\in[0,1]$, we have
\[
 J_\rho(\mathcal H(t,u))\le J_\rho(u),\quad
 \mathscr F_\rho(\mathcal G(t,\mu))
 \le\mathscr F_\rho(\mu).
\]
The homotopy equivalence holds for all $\rho\ge0$
and all exponents $\al_i>-1$.
\end{lem}

\begin{proof}
We first prove the continuity of the functionals and maps.
By the definition of $H$, we have
\[
 \ln H=\ln h+N\sum_{i=1}^m\al_i\ln|x-p_i|
 \in L^b(\Om)
\]
for any  finite $b\ge1$. Hence the map
$\mu\mapsto\int_\Om\mu\ln H\,\d x$ is continuous on
$L^p(\Om)$.
If $\mu_j\to\mu$ in $L^p(\Om)$, then
$\mu_j\ln\mu_j\to\mu\ln\mu$ in measure.
Since
\[
 |t\ln t|\le C_\var(1+t^{1+\var}),
 \quad t\ge0,\quad 0<\var<p-1,
\]
the functions $\mu_j\ln\mu_j$ are uniformly integrable.
By Vitali's theorem,
\[
 \mu_j\ln\mu_j\to\mu\ln\mu
 \quad\text{ in }\, L^1(\Om).
\]
Thus $\operatorname{Ent}_H$ is finite and continuous
on $\mathscr P_p$.

The Sobolev embedding $E\hookrightarrow L^{p'}(\Om)$
implies $L^p(\Om)\hookrightarrow E^*$.
Since $Q':E\to E^*$ is a homeomorphism,
$\mathscr S_\rho=(Q')^{-1}(\rho\,\cdot)$ is continuous.
Consequently, $\mathscr F_\rho$ is also continuous.

To prove the continuity of $\Theta$, let $u_j\to u$ in $E$.
By the exponential integrability estimates in Section~\ref{sec:2},
$e^{u_j}\to e^u$ in $L^c(\Om)$ for any  finite $c\ge1$.
Choose $p<b<s_0$ and $c$ such that $1/p=1/b+1/c$.
By H\"older's inequality,
\[
 He^{u_j}\to He^u\quad\text{ in }\, L^p(\Om).
\]
Since $M(u_j)\to M(u)>0$, we obtain
$\Theta(u_j)\to\Theta(u)$ in $L^p(\Om)$.

We next compare $J_\rho$ and $\mathscr F_\rho$.
For $u\in E$ and $\mu\in\mathscr P_p$, define
\[
 \mathscr L_\rho(u,\mu)
 =Q(u)-\rho\int_\Om u\mu\,\d x
 +\rho\operatorname{Ent}_H(\mu).
\]
Since
$\ln\Theta(u)=\ln H+u-\ln M(u)\in L^{p'}(\Om)$,
Gibbs' inequality yields
\[
 \mathscr L_\rho(u,\mu)
 =J_\rho(u)+\rho\int_\Om
 \mu\ln\frac{\mu}{\Theta(u)}\,\d x
 \ge J_\rho(u).
\]
Equality holds when $\mu=\Theta(u)$.

On the other hand, let $v=\mathscr S_\rho(\mu)$.
Since $Q'(v)=\rho\mu$ and $Q$ is convex, $v$ minimizes
the functional
\[
 u\mapsto Q(u)-\rho\int_\Om u\mu\,\d x.
\]
Testing $Q'(v)=\rho\mu$ with $v$, we obtain
\[
 \rho\int_\Om v\mu\,\d x=NQ(v).
\]
Consequently,
\[
 \mathscr L_\rho(u,\mu)
 \ge\mathscr L_\rho(v,\mu)
 =\rho\operatorname{Ent}_H(\mu)-(N-1)Q(v)
 =\mathscr F_\rho(\mu).
\]

By H\"older's inequality,
$\norm{Q'(v)}_{E^*}\le\norm v_E^{N-1}$.
For $v\ne0$, testing $Q'(v)$ with $v/\norm v_E$
gives equality. Thus
\[
 \rho\norm{\mu}_{E^*}
 =\norm{Q'(v)}_{E^*}=\norm v_E^{N-1}.
\]
Since $q=N/(N-1)$, we have
\[
 (N-1)Q(v)=\frac{\rho^q}{q}\norm{\mu}_{E^*}^q.
\]
Hence \eqref{eq:dual-entropy-functional} holds,
including the case $\rho=0$.

The above inequalities imply
\[
 \mathscr F_\rho(\Theta(u))
 \le\mathscr L_\rho(u,\Theta(u))=J_\rho(u),\quad 
 J_\rho(\mathscr S_\rho(\mu))
 \le\mathscr L_\rho(\mathscr S_\rho(\mu),\mu)
 =\mathscr F_\rho(\mu).
\]
Therefore, $\Theta$ maps $J_\rho^a$ into
$\mathscr F_\rho^a$, while $\mathscr S_\rho$ maps
$\mathscr F_\rho^a$ into $J_\rho^a$.

We now prove the homotopy equivalence.
For $u\in J_\rho^a$, set $\nu=\Theta(u)$ and
$v=\mathscr S_\rho(\nu)$. Since
$\mathscr L_\rho(\cdot,\nu)$ is convex and has a minimum
at $v$, we obtain, for any  $t\in[0,1]$,
\[
 J_\rho((1-t)u+tv)
 \le\mathscr L_\rho((1-t)u+tv,\nu)
 \le(1-t)\mathscr L_\rho(u,\nu)
       +t\mathscr L_\rho(v,\nu)
 \le J_\rho(u).
\]
Thus $\mathcal H(t,u)\in J_\rho^a$.

Similarly, for $\mu\in\mathscr F_\rho^a$, set
$v=\mathscr S_\rho(\mu)$ and $\nu=\Theta(v)$.
By convexity of $\operatorname{Ent}_H$,
$\mathscr L_\rho(v,\cdot)$ is convex and has a minimum
at $\nu$. Hence
\[
 \mathscr F_\rho((1-t)\mu+t\nu)
 \le\mathscr L_\rho(v,(1-t)\mu+t\nu)
 \le(1-t)\mathscr L_\rho(v,\mu)
       +t\mathscr L_\rho(v,\nu)
 \le\mathscr F_\rho(\mu).
\]
Thus $\mathcal G(t,\mu)\in\mathscr F_\rho^a$.

By continuity of $\Theta$ and $\mathscr S_\rho$,
both $\mathcal H$ and $\mathcal G$ are continuous.
Moreover,
\[
 \mathcal H(0,u)=u,\quad 
 \mathcal H(1,u)=\mathscr S_\rho(\Theta(u)),\quad
 \mathcal G(0,\mu)=\mu,\quad
 \mathcal G(1,\mu)=\Theta(\mathscr S_\rho(\mu)).
\]
Therefore,
\[
 \mathscr S_\rho\circ\Theta\simeq\Id_{J_\rho^a},
 \quad
 \Theta\circ\mathscr S_\rho\simeq\Id_{\mathscr F_\rho^a},
\]
which proves $J_\rho^a\simeq\mathscr F_\rho^a$.

Finally, if $\rho>0$ and $\rho\notin\Si_H$, then
\eqref{eq:all-parameter-euler} and the above homotopy
equivalence yield
\[
 d_\rho(\Om,H)=1-\chi(\mathscr F_\rho^{-L})
\]
for all sufficiently large $L$.
\end{proof}

For a smooth domain $D$, let $Q_D$ be the Dirichlet energy
on $W^{1,N}_0(D)$ and let $Q_D^*$ be its convex conjugate.
We now study how $Q_D^*$ in \eqref{eq:dual-entropy-functional}
changes when a domain is divided into several subdomains.
For this purpose, we consider densities that vanish
in fixed neighborhoods of the new boundaries.

\begin{lem}\label{lem:buffered-potential-comparison}
Let $\Om$ be a smooth bounded domain and
$D=\bigsqcup_{i=1}^kD_i\subset\Om$, where $D_i$ are
the connected components of $D$ with smooth boundaries.
Set $\Ga=\pa D\cap\Om$ and assume $\Ga\Subset\Om$.
Choose $d>0$ such that
\[
 U_d(\Ga)=\{x\in\Om:\dist(x,\Ga)<d\}\Subset\Om.
\]
Let $\nu\in L^p(\Om)$ for some $p>1$ and assume
\[
 \nu\ge0,\quad
 \nu=0\quad\text{a.e. on }(\Om\setminus D)\cup U_d(\Ga).
\]
Set
\[
 M=\int_\Om\nu\,\d x,\quad \nu_i=\nu|_{D_i}.
\]
Then
\[
 0\le Q_\Om^*(\nu)-\sum_{i=1}^kQ_{D_i}^*(\nu_i)
 \le CM^{N/(N-1)},
\]
where $C>0$ depends only on $\Om$, $D_1,\ldots,D_k$,
$d$, and $N$. In particular, $C$ is independent of $p$
and $\norm{\nu}_{L^p(\Om)}$.
If $\Ga=\varnothing$, no condition near $\Ga$ is needed.
\end{lem}

\begin{proof}
If $M=0$, then $\nu=0$ and the conclusion is immediate.
Assume $M>0$ and set
\[
 u=(Q_\Om')^{-1}(\nu),\quad
 v_i=(Q_{D_i}')^{-1}(\nu_i).
\]
By the maximum principle, $u\ge0$ and $v_i\ge0$.
By \eqref{eq:BM} with zero boundary values, for any
finite $a\ge1$, we have
\[
 \norm{u}_{L^a(\Om)}\le C_aM^{1/(N-1)}.
\]
Since $\nu=0$ in $U_d(\Ga)$, the function $u$ is
$N$-harmonic in $U_d(\Ga)$. By the local estimates
for $N$-harmonic functions,
\[
 \sup_\Ga u\le CM^{1/(N-1)}.
\]
Moreover, $u=0$ on $\pa\Om$. Hence
\[
 0\le u\le CM^{1/(N-1)}\quad\text{ on }\, \pa D_i.
\]
Since $u$ and $v_i$ satisfy
\[
 -\Delta_Nu=-\Delta_Nv_i=\nu_i\quad\text{ in }\, D_i,
 \quad v_i=0\quad\text{ on }\, \pa D_i,
\]
we can compare $u$ with $v_i$ and
$v_i+CM^{1/(N-1)}$. By the comparison principle,
\[
 0\le u-v_i\le CM^{1/(N-1)}
 \quad\text{ in }\, D_i.
\]

Recall that
\[
 Q_G^*(f)=\frac{N-1}{N}
 \int_G f(Q_G')^{-1}(f)\,\d x.
\]
Since $\nu=0$ outside $D$, we obtain
\[
 Q_\Om^*(\nu)-\sum_{i=1}^kQ_{D_i}^*(\nu_i)
 =\frac{N-1}{N}\sum_{i=1}^k
 \int_{D_i}\nu_i(u-v_i)\,\d x
 \le CM^{1/(N-1)}
 \sum_{i=1}^k\int_{D_i}\nu_i\,\d x=CM^{N/(N-1)}.
\]
Since $\nu_i\ge0$ and $u-v_i\ge0$, the difference
on the left is nonnegative. The proof is complete.
\end{proof}

To apply Lemma~\ref{lem:buffered-potential-comparison},
we first remove a small amount of mass near the new
boundaries and then normalize the remaining density.
The next lemma estimates the change in energy in terms
of the mass removed.

\begin{lem}\label{lem:small-mass-trimming}
Assume \eqref{eq:assumptions} and fix $0<\rho<\rho'$.
Let $\mu\in\mathscr P_p$ and let $w$ be a measurable
function with $0\le w\le1$. Set
\[
 m=\int_\Om(1-w)\mu\,\d x,\quad
 T_w\mu=\frac{w\mu}{1-m},
\]
where $m<1$. Then there exist $\eta_0\in(0,1)$
and constants $C_\eta\to0$ as $\eta\to0$ such that,
for $m\le\eta\le\eta_0$,
\[
 \frac{\mathscr F_{\rho'}(T_w\mu)}{\rho'}
 \le\frac{\mathscr F_\rho(\mu)}{\rho(1-m)}+C_\eta.
\]
In particular, if $\mathscr F_\rho(\mu)\le0$, then
\be\label{eq:small-trimming-negative-bound}
 \frac{\mathscr F_{\rho'}(T_w\mu)}{\rho'}
 \le\frac{\mathscr F_\rho(\mu)}{\rho}+C_\eta,
\ee
where the constant $C_\eta$ is independent of $w$ and $\mu$.
Moreover, for $m\le\eta_0$, the map
$(w,\mu)\mapsto T_w\mu$ is continuous with respect to
the strong topology of $L^\infty(\Om)\times L^p(\Om)$.
\end{lem}
\begin{proof}
Write $f=w\mu$ and $g=(1-w)\mu$. Use the same entropy formula for the
nonnegative densities $f$ and $g$, with $0\ln  0=0$. Integrating the pointwise
inequality $(f+g)\ln (f+g)\ge f\ln  f+g\ln  g$ proves
$\operatorname{Ent}_H(f)\le\operatorname{Ent}_H(\mu)-\operatorname{Ent}_H(g)$.
Set $\kappa=\rho'/\rho$ and choose $\var >0$ so that
$1+\var \le\kappa^{q-1}$. The formula
$Q_\Om^*(\nu)=\|\nu\|_{E^*}^q/q$, the triangle inequality and
Young's inequality yield
\[
 Q_\Om^*(\rho\mu)
 \le (1+\var )Q_\Om^*(\rho f)+C_\var  Q_\Om^*(\rho g).
\]
Since $(\kappa/(1-m))^{q-1}\ge1+\var $, we use homogeneity and the
entropy inequality to obtain
\[
 \frac{\mathscr F_{\rho'}(T_w\mu)}{\rho'}
 \le \frac{\mathscr F_\rho(\mu)
        +C_\var  Q_\Om^*(\rho g)-\rho\operatorname{Ent}_H(g)}
          {\rho(1-m)}
     -\ln (1-m).
\]

For the terms $C_\var Q_\Om^*(\rho g)$ and $\rho\operatorname{Ent}_H(g)$, choose $a_0>0$ and $C_0\ge0$ satisfying $\ln  M(u)\le a_0Q(u)+C_0$ for all $u\in E$.
Gibbs' inequality and the definition of the convex conjugate imply
\[
 Q_\Om^*(\tau v)\le\tau(\operatorname{Ent}_H(v)+C_0)\quad \text{ for }\,
v\in\mathscr P_p,\ 0<\tau a_0\le1.
\]
When $m>0$, apply this estimate with $v=g/m$ and
$\tau=C_\var ^{1/(q-1)}\rho m$, after reducing $\eta_0$ if needed.
Homogeneity and
$m\operatorname{Ent}_H(g/m)=\operatorname{Ent}_H(g)-m\ln  m$ show that
\[
 C_\var  Q_\Om^*(\rho g)-\rho\operatorname{Ent}_H(g)
 \le \rho m(C_0-\ln  m).
\]
At $m=0$ each term on the left vanishes. For
$\eta_0\le\min\{e^{-1},1/2\}$, the total error is therefore at most
$\eta(C_0-\ln \eta)/(1-\eta)-\ln (1-\eta)$, which tends to zero.
This proves the estimate and \eqref{eq:small-trimming-negative-bound}.
Continuity follows from bounded multiplication and the positive lower
estimate $1-m\ge1-\eta_0>0$.
For a product $w=\prod_j(1-t_j\vartheta_j)$ with
$0\le t_j,\vartheta_j\le1$,
\[
 \int_\Om(1-w)\mu\, \d  x
 \le\int_\Om\Big[1-\prod_j(1-\vartheta_j)\Big]\mu\, \d  x.
\]
The estimate is therefore uniform on the cutoff cube.
\end{proof}

\subsection{Domains with holes and disconnected domains}

We first compare the degree on $\Om$ with the degree
on the domain obtained by removing small balls around
positive singularities. The comparison will be used
to compute the contribution of a positive singularity
to the degree.

Fix the data in \eqref{eq:assumptions} and let $J\subset I_+$.
Choose $\var_0>0$ such that the closed balls
$\overline B_{2\var_0}(p_j)$, $j\in J$, are pairwise
disjoint, contained in $\Om$, and contain no other
singular points. For $0<\var\le\var_0$, set
\[
 \Om_\var=\Om\setminus
\bigcup_{j\in J}\overline B_\var(p_j),
 \quad \Om_0=\Om.
\]
For sufficiently small $\var_0$, all domains $\Om_\var$
are connected. By extending functions by zero outside
$\Om_\var$, we regard $E_\var=W^{1,N}_0(\Om_\var)$
as a closed subspace of $E$, with $E_0=E$.

For $u\in E$, define
\[
 M_\var(u)=\int_{\Om_\var}He^u\,\d x
\]
and let $K_\var(u)$ be the zero extension of
\[
 (Q_{\Om_\var}')^{-1}
 \Big(\rho\frac{He^u}{M_\var(u)}\Big).
\]
Then $K_0=K_{\rho,H}$, and the fixed points of $K_\var$
are exactly the solutions of \eqref{maineq} on
$\Om_\var$, extended by zero to $\Om$.

Since a point has zero $N$-capacity, we expect the
Dirichlet solution operators on $\Om_\var$ to converge
to the operator on $\Om$ as $\var\to0$.
The following lemma makes this convergence precise.

\begin{lem}\label{lem:puncture-resolvent}
As $\var\to0$, the spaces $E_\var$ converge to $E$
in the sense of Mosco. More generally, let
$\var_n\to\var_*\in[0,\var_0]$ and $f_n\to f$ in $E^*$.
Then the minimizers
\[
 z_n=\operatorname*{argmin}_{v\in E_{\var_n}}
 \{Q(v)-\langle f_n,v\rangle\}
\]
converge strongly in $E$ to
\[
 z_*=\operatorname*{argmin}_{v\in E_{\var_*}}
 \{Q(v)-\langle f,v\rangle\}.
\]
Moreover, if $u_n\weak u$ in $E$ and
$\var_n\to\var_*$, then
\[
 K_{\var_n}(u_n)\to K_{\var_*}(u)
 \quad\text{ strongly in }\,E.
\]
In particular, the map $(\var,u)\mapsto K_\var(u)$
is continuous, and the family $\{K_\var\}_{0\le\var\le\var_0}$
is collectively compact on bounded subsets of $E$.
Furthermore, for any  $R>0$,
\be\label{eq:puncture-map-convergence}
 \sup_{\norm u_E\le R}\norm{K_\var(u)-K_0(u)}_E
 \to0\quad\text{ as }\, \var\to0.
\ee
\end{lem}
\begin{proof}
We first prove the Mosco convergence of $E_\var$.
For each $j\in J$, choose a logarithmic cutoff
$\zeta_{\var,j}$ such that
\[
 \zeta_{\var,j}=1\quad\text{ in }\,B_\var(p_j),\quad
 \zeta_{\var,j}=0\quad\text{ in }\,\R^N\setminus B_{\sqrt\var}(p_j),
\]
and
\[
 \int_\Om|\nabla\zeta_{\var,j}|^N\,\d x
 \le C|\ln\var|^{1-N}\to0
 \quad\text{ as }\, \var\to0.
\]
For any $v\in C_c^\infty(\Om)$, set
\[
 v_\var=v\prod_{j\in J}(1-\zeta_{\var,j}).
\]
Then $v_\var\in E_\var$ and $v_\var\to v$ strongly
in $E$. By density, every $v\in E$ can be approximated
in $E$ by functions in $E_\var$ as $\var\to0$.

More generally, suppose $\var_n\to\var_*$ and
$v_n\in E_{\var_n}$ with $v_n\weak v$ in $E$.
If $\var_*=0$, then $v\in E_0=E$.
If $\var_*>0$, then $v=0$ in $B_{\var_*}(p_j)$
for any  $j\in J$, and hence $v\in E_{\var_*}$.
Conversely, every function in $C_c^\infty(\Om_{\var_*})$
lies in $E_{\var_n}$ for sufficiently large $n$.
By density, every $v\in E_{\var_*}$ can be approximated
strongly in $E$ by a sequence $v_n\in E_{\var_n}$.
Together with the case $\var_*=0$ proved above,
we obtain Mosco convergence.

We next prove the convergence of the minimizers.
Since $z_n$ minimizes $Q(v)-\langle f_n,v\rangle$
on $E_{\var_n}$, we have
\[
 Q(z_n)-\langle f_n,z_n\rangle\le0.
\]
Since $f_n\to f$ in $E^*$, the sequence $z_n$
is bounded in $E$. Up to a subsequence, let
$z_n\weak z$ in $E$. By the Mosco convergence,
$z\in E_{\var_*}$.

For any $v\in E_{\var_*}$, choose $v_n\in E_{\var_n}$
such that $v_n\to v$ strongly in $E$.
By weak lower semicontinuity and the definition of $z_n$,
\[
 Q(z)-\langle f,z\rangle
 \le\liminf_{n\to\infty}
(Q(z_n)-\langle f_n,z_n\rangle) \le\limsup_{n\to\infty}
(Q(z_n)-\langle f_n,z_n\rangle) \le Q(v)-\langle f,v\rangle.
\]
Hence $z=z_*$, where $z_*$ is the unique minimizer
on $E_{\var_*}$. Taking $v=z_*$ also yields
$Q(z_n)\to Q(z_*)$. By uniform convexity of
$L^N(\Om;\R^N)$, we obtain $z_n\to z_*$ strongly in $E$.

We now prove the convergence of $K_{\var_n}(u_n)$.
Suppose $u_n\weak u$ in $E$ and $\var_n\to\var_*$.
By the exponential estimates in Section~\ref{sec:2},
\[
 He^{u_n}\to He^u\quad\text{ in }\, L^b(\Om)
\]
for some $b>1$. Moreover,
$\mathbf1_{\Om_{\var_n}}\to\mathbf1_{\Om_{\var_*}}$
almost everywhere in $\Om$. Hence
\[
 M_{\var_n}(u_n)\to M_{\var_*}(u)>0
\]
and
\[
 \frac{\rho He^{u_n}\mathbf1_{\Om_{\var_n}}}
 {M_{\var_n}(u_n)}
 \to
 \frac{\rho He^u\mathbf1_{\Om_{\var_*}}}
 {M_{\var_*}(u)}
 \quad\text{ in }\, L^c(\Om)
\]
for some $c>1$. By the convergence of the minimizers
proved above,
\[
 K_{\var_n}(u_n)\to K_{\var_*}(u)
 \quad\text{strongly in }E.
\]
Since $E$ is reflexive, the same argument proves
continuity of $(\var,u)\mapsto K_\var(u)$ and collective
compactness on bounded subsets of $E$.

Finally, suppose that \eqref{eq:puncture-map-convergence}
does not hold. Then there exist $\eta>0$,
$\var_n\to0$, and $u_n\in E$ with $\norm{u_n}_E\le R$
such that
\[
 \norm{K_{\var_n}(u_n)-K_0(u_n)}_E\ge\eta.
\]
Up to a subsequence, $u_n\weak u$ in $E$.
By the convergence proved above,
\[
 K_{\var_n}(u_n)\to K_0(u),\quad
 K_0(u_n)\to K_0(u)
 \quad\text{ in }\, E,
\]
which is a contradiction. Hence
\eqref{eq:puncture-map-convergence} holds.
\end{proof}

We use a local Pohozaev inequality to rule out
concentration near a shrinking hole when the mass
$c_N\beta_j$ at its center is greater than $\rho$.
The condition $\rho\notin\Si_H$ then rules out
concentration at all other points.

\begin{prop}\label{prop:puncture-compactness}
Let $\rho>0$, $\rho\notin\Si_H$, and assume that
\be\label{eq:puncture-mass-condition}
 \beta_j>\frac{\rho}{c_N},\quad j\in J.
\ee
Then, for all sufficiently small $\var>0$, the solutions
of \eqref{maineq} on $\Om_\var$ are uniformly bounded
in $L^\infty(\Om_\var)$, and their zero extensions
are uniformly bounded in $E$.

Moreover, for any sequence $\var_n\to0$ and corresponding
solutions $u_n$ on $\Om_{\var_n}$, the zero extensions
of $u_n$ have a subsequence converging strongly in $E$
to a solution of \eqref{maineq} on $\Om$.
\end{prop}

\begin{proof}
\textbf{Step 1:} Uniform estimates and convergence away from
the concentration set.

Suppose that there exist $\var_n\downarrow0$ and solutions
$u_n$ on $\Om_{\var_n}$ such that
$\norm{u_n}_{L^\infty(\Om_{\var_n})}\to\infty$.
We extend $u_n$ by zero to $\Om$ and set
\[
 b_n=\frac{\rho}{M_{\var_n}(u_n)},\quad
 f_n=b_nHe^{u_n}\mathbf1_{\Om_{\var_n}}.
\]
Since $u_n\ge0$ and $\int_{\Om_{\var_n}}H\,\d x$
has a uniform positive lower bound, we have
\[
 0<b_n\le C,\quad \int_\Om f_n\,\d x=\rho.
\]
Let $w_n$ be the solution of
\[
 -\Delta_N w_n=f_n\quad\text{ in }\, \Om,\quad
 w_n=0\quad\text{ on }\, \pa\Om.
\]
By the comparison principle, $0\le u_n\le w_n$ in $\Om$.
By \eqref{eq:BM}, for some $\delta>0$,
\[
 \int_\Om e^{\delta u_n}\,\d x
 \le\int_\Om e^{\delta w_n}\,\d x\le C.
\]
Testing the equation for $u_n$ with
$\min\{1,(u_n-k)^+\}$ for each integer $k\ge0$, we obtain
\[
 \int_{\{x\in\Om:k<u_n(x)<k+1\}}
 |\nabla u_n|^N\,\d x\le\rho.
\]
Moreover,
\[
 |\{x\in\Om:u_n(x)>k\}|\le Ce^{-\delta k}.
\]
By H\"older's inequality, for $1\le s<N$,
\[
 \int_\Om|\nabla u_n|^s\,\d x
 \le \rho^{s/N}\sum_{k=0}^\infty
 |\{x\in\Om:k<u_n(x)<k+1\}|^{1-s/N}
 \le C\sum_{k=0}^\infty e^{-\delta k(1-s/N)}
 \le C.
\]
Hence $u_n$ is uniformly bounded in $W^{1,s}_0(\Om)$
for any  $1\le s<N$. Up to a subsequence, we have
 \begin{align*}
 u_n&\to U\ge0\quad\text{ in }\, L^1(\Om)
 \text{ and a.e. in }\Om,\\
 b_n&\to b\ge0,\\
 f_n\,\d x&\weak\mu\quad\text{ on }\, \overline\Om.
 \end{align*}

As in Step 1 of the proof of Proposition~\ref{prop:blowup},
choose a sufficiently small $\var_0'>0$ and define
\[
 S=\{x\in\overline\Om:\mu(\{x\})\ge\var_0'\}.
\]
Since $\mu$ is a finite measure, $S$ is finite.
Away from the points $p_j$ with $j\in J$, the local
estimates in Proposition~\ref{prop:blowup} remain valid
for sufficiently large $n$.

For $j\in J$ with $p_j\notin S$, choose a small fixed
ball centered at $p_j$. Since $\al_j>0$, the weight $H$
is bounded on the ball. Moreover,
\[
 f_n\le CH e^{w_n}.
\]
Applying the local estimates in Step 1 of the proof of
Proposition~\ref{prop:blowup} to $w_n$, we obtain a
uniform $L^t$ estimate for $f_n$ on a smaller ball,
for some $t>1$. The local regularity estimates for $w_n$
and the inequality $u_n\le w_n$ then imply that $u_n$
is uniformly bounded near $p_j$.
If $S=\varnothing$, finitely many local neighborhoods
cover $\overline\Om$, and $u_n$ is uniformly bounded
in $L^\infty(\Om)$. Thus $S\ne\varnothing$.

Set $P_J=\{p_j:j\in J\}$. On compact subsets of
$\overline\Om\setminus(S\cup P_J\cup P_-)$,
$f_n$ is uniformly bounded. By the interior and boundary
regularity estimates, up to a subsequence,
\[
 u_n\to U\quad\text{ in }\, 
 C^1_{loc}
 (\overline\Om\setminus(S\cup P_J\cup P_-)).
\]
Near each negative singularity $p_i\notin S$, the
comparison argument in Step 1 of the proof of
Proposition~\ref{prop:blowup}, together with
\eqref{eq:stability}, yields local uniform convergence
of $u_n$ to $U$.

In particular, $\nabla u_n\to\nabla U$ almost everywhere
in $\Om$. For any $1\le s<N$, choose $s<t<N$.
Since $\nabla u_n$ is uniformly bounded in $L^t(\Om)$,
Vitali's theorem implies
\[
 u_n\to U\quad\text{ in }\, W^{1,s}_0(\Om),
 \quad 1\le s<N.
\]
Taking $s>N-1$, we also obtain
\[
 |\nabla u_n|^{N-2}\nabla u_n
 \to|\nabla U|^{N-2}\nabla U
 \quad\text{ in }\, L^1(\Om;\R^N).
\]

\textbf{Step 2:} The limiting equation near the holes.

Let $\nu$ be the outward unit normal to $\Om_{\var_n}$.
Since $u_n\ge0$ in $\Om_{\var_n}$ and $u_n=0$ on
$\pa\Om_{\var_n}$, we have $\pa_\nu u_n\le0$ on
$\pa\Om_{\var_n}$. Define the nonnegative measure
\[
 \nu_n=\sum_{j\in J}|\pa_\nu u_n|^{N-1}
 \mathcal H^{N-1}\big|_{\pa B_{\var_n}(p_j)}.
\]
For any $\phi\in C_c^\infty(\Om)$, integration by parts
on $\Om_{\var_n}$ yields
\[
 \int_\Om|\nabla u_n|^{N-2}\nabla u_n\cdot\nabla\phi\,\d x
 =\int_\Om f_n\phi\,\d x-\int_\Om\phi\,\d\nu_n.
\]
Hence the zero extension of $u_n$ satisfies
\be\label{eq:puncture-signed-equation}
 -\Delta_Nu_n=f_n\,\d x-\nu_n
 \quad\text{ in }\, \Om.
\ee
Integrating the equation on $\Om_{\var_n}$, we also obtain
\[
 \rho=\int_{\pa\Om}|\pa_\nu u_n|^{N-1}\,\d S
 +\nu_n(\Om),\quad 0\le\nu_n(\Om)\le\rho.
\]
The two identities can be justified by using cutoff
functions near $\pa\Om_{\var_n}$ and then passing to
the limit. Near negative singularities, the integrals
over small spheres vanish along a suitable sequence
of radii, as in the proof of \eqref{eq:poho}.

Since the holes shrink to the points $p_j$, $j\in J$,
up to a subsequence, we have
\[
 \nu_n\weak\sum_{j\in J}\eta_j\delta_{p_j},
 \quad \eta_j\ge0.
\]
By the strong convergence of
$|\nabla u_n|^{N-2}\nabla u_n$ in $L^1(\Om;\R^N)$
proved in Step 1, we can pass to the limit in
\eqref{eq:puncture-signed-equation}.

Since the holes stay away from $\pa\Om$, the argument
in Step 2 of the proof of Proposition~\ref{prop:blowup}
rules out concentration on $\pa\Om$.
Moreover, the local estimates outside $S$ imply
\[
 \mu=F\,\d x\quad\text{ on }\, \overline\Om\setminus S,
 \quad F=bHe^U.
\]
By Fatou's lemma, $F\in L^1(\Om)$. Consequently,
\[
 -\Delta_NU
 =F+\sum_{p\in S}\mu(\{p\})\delta_p
 -\sum_{j\in J}\eta_j\delta_{p_j}
 \quad\text{ in }\, \Om.
\]

Fix $j\in J$ and set
\[
 \si_j=\mu(\{p_j\}),\quad \ga_j=\si_j-\eta_j.
\]
Choose a small ball centered at $p_j$ containing
no other singular or concentration points.
Since $U\ge0$ and $-\Delta_NU=F\ge0$ in the punctured
ball, $U$ is $N$-superharmonic there.
Since a point has zero $N$-capacity,
\cite[Theorem~6.3(d)]{B06} implies that $U$ admits
an $N$-superharmonic extension across $p_j$.
Its bounded truncations $\min\{U,k\}$ belong to $W^{1,N}_{\rm loc}$
and are weak supersolutions. Testing them against a nonnegative
smooth function and passing to the limit using
$|\nabla U|^{N-1}\in L^1_{\rm loc}$ proves that $-\Delta_N U$
is a nonnegative distribution on the ball. Hence $\ga_j\ge0$.
Since $\eta_j\ge0$, we obtain
\be\label{eq:puncture-net-atom}
 0\le\ga_j\le\si_j.
\ee

If $\ga_j=0$, then $U$ satisfies
\[
 -\Delta_NU=F=bHe^U
\]
on the whole ball centered at $p_j$.
Since $F\in L^1(\Om)$, we have
\[
 \int_{B_R(p_j)}F\,\d x\to0\quad\text{as }R\to0.
\]
Applying \eqref{eq:BM} to solutions with truncated
right-hand sides $\min\{F,T\}$ and passing to the
limit by uniqueness for measure data, we obtain
$e^U\in L^a(B_R(p_j))$ for any prescribed finite
$a\ge1$, provided $R$ is sufficiently small.
Since $H$ is bounded near $p_j$, it follows that
$F\in L^t(B_R(p_j))$ for some $t>1$.

As in Step 2 of the proof of Proposition~\ref{prop:blowup},
\eqref{eq:monotonicity} implies that
$U\in W^{1,N}(B_R(p_j))$.
The local regularity estimates then yield
$U\in C^{1,\eta}$ near $p_j$ for some $\eta>0$.
If $b=0$, then $U$ is $N$-harmonic across $p_j$.

\textbf{Step 3:} The local Pohozaev inequality.

Assume $\si_j>0$. Choose $R>0$ sufficiently small and write
\[
 H(x)=|x-p_j|^{N\al_j}V_j(x),
 \quad V_j>0,\quad V_j\in C^1(B_R(p_j)).
\]
We apply the Pohozaev identity on
$B_R(p_j)\setminus\overline B_{\var_n}(p_j)$
with the vector field $x-p_j$.
Since $u_n=0$ and $\nabla u_n=(\pa_\nu u_n)\nu$
on $\pa B_{\var_n}(p_j)$, the integral over the inner
boundary is
\[
 -\var_n\int_{\pa B_{\var_n}(p_j)}
 \Big(b_nH+\frac{N-1}{N}|\nabla u_n|^N\Big)\,\d S
 \le0.
\]
Moreover,
\[
 \Big|\int_{B_R(p_j)\setminus\overline B_{\var_n}(p_j)}
 (x-p_j)\cdot\nabla\ln V_j\,f_n\,\d x\Big|
 \le CR\rho.
\]
For fixed small $R$, the local $C^1$ convergence in Step 1
allows us to pass to the limit in the integrals over
$\pa B_R(p_j)$.
Dropping the nonpositive integral over the inner boundary
and letting $n\to\infty$, we obtain
\[
 N(1+\al_j)\mu(B_R(p_j))
 \le R\int_{\pa B_R(p_j)}
 \Big[F+|\nabla U|^{N-2}(\pa_\nu U)^2
 -\frac1N|\nabla U|^N\Big]\,\d S+CR\rho.
\]

Set $a_j=(\ga_j/S_N)^{1/(N-1)}$.
If $\ga_j=0$, then $U$ is regular near $p_j$ by Step 2,
and the boundary integral above tends to zero as $R\to0$.
If $\ga_j>0$ and $b>0$, Lemma~\ref{lem:singularity}
implies $a_j<N(1+\al_j)$ and
\[
 U(x)=a_j\ln\frac1{|x-p_j|}+O(1),\quad
 |x-p_j|\nabla U(x)
 =-a_j\frac{x-p_j}{|x-p_j|}+o(1).
\]
In particular,
\[
 R\int_{\pa B_R(p_j)}F\,\d S
 \le CR^{N(1+\al_j)-a_j}\to0.
\]
If $\ga_j>0$ and $b=0$, then $F=0$, and the same
asymptotic expansions follow from
\cite[Theorem~A.5]{EM15}.

Therefore, letting $R\to0$, we obtain
\be\label{eq:puncture-pohozaev-inequality}
 N(1+\al_j)\si_j
 \le\frac{N-1}{N}S_Na_j^N,
 \quad a_j=(\ga_j/S_N)^{1/(N-1)}.
\ee

\textbf{Step 4:} Exclusion of concentration and strong convergence.

We first prove that $b=0$. Suppose $b>0$.
If $p_j\in S$ for some $j\in J$, then $\si_j>0$.
By \eqref{eq:puncture-pohozaev-inequality}, we have
$\ga_j>0$. Lemma~\ref{lem:singularity} then implies
$a_j<N(1+\al_j)$. Together with
\eqref{eq:puncture-net-atom}, we obtain
\[
 \frac{N-1}{N}S_Na_j^N
 =\frac{N-1}{N}a_j\ga_j
 <(N-1)(1+\al_j)\si_j
 <N(1+\al_j)\si_j,
\]
which contradicts \eqref{eq:puncture-pohozaev-inequality}.
For any concentration point outside $\{p_j:j\in J\}$,
the argument in Proposition~\ref{prop:blowup} also
rules out the case $b>0$. Since $S\ne\varnothing$,
we conclude that $b=0$.

We next rule out concentration at the points $p_j$,
$j\in J$. If $\si_j>0$, then
\eqref{eq:puncture-pohozaev-inequality} and
\eqref{eq:puncture-net-atom} imply
\[
 \begin{aligned}
 N(1+\al_j)\si_j
 &\le\frac{N-1}{N}S_N^{-1/(N-1)}\ga_j^{N/(N-1)}\\
 &\le\frac{N-1}{N}S_N^{-1/(N-1)}\si_j^{N/(N-1)}.
 \end{aligned}
\]
Consequently,
\[
 \rho\ge\si_j\ge c_N\beta_j,
\]
which contradicts \eqref{eq:puncture-mass-condition}.
Hence no concentration occurs at $p_j$ for any $j\in J$.

Since $b=0$, we have $F=bHe^U=0$.
Fix a remaining concentration point $p$ and a small ball
about $p$ separated from the holes. On this ball the original
equations hold for all sufficiently large $n$. The limiting
atom has no contribution from $\nu_n$, so the Pohozaev
calculation in the proof of Proposition~\ref{prop:blowup}
gives, for $\si=\mu(\{p\})$,
\[
 N(1+\al)\si=\frac{N-1}{N}S_N
                  \Big(\frac{\si}{S_N}\Big)^{N/(N-1)}.
\]
Here $\al=0$ at an ordinary point and $\al=\al_i$ at a
retained singularity $p_i$. Hence the corresponding mass is
$c_N$ or $c_N\beta_i$.
Since
\[
 \rho=\mu(\overline\Om)=\sum_{p\in S}\mu(\{p\}),
\]
we obtain $\rho\in\Si_H$, a contradiction.
Therefore, all solutions on $\Om_\var$ are uniformly
bounded in $L^\infty(\Om_\var)$ for sufficiently small
$\var>0$.

Now take any sequence of solutions with $\var_n\to0$.
Testing the equations with $u_n$ gives
\[
 \norm{u_n}_E^N
 =\int_{\Om_{\var_n}}f_nu_n\,\d x
 \le\rho\norm{u_n}_{L^\infty(\Om_{\var_n})}\le C.
\]
Up to a subsequence, $u_n\rightharpoonup U$ in $E$.
Since $u_n=K_{\var_n}(u_n)$,
Lemma~\ref{lem:puncture-resolvent} implies
\[
 u_n\to K_0(U)=U\quad\text{strongly in }E.
\]
Thus $U$ is a solution of \eqref{maineq} on $\Om$.
\end{proof}

By Proposition~\ref{prop:puncture-compactness}, all fixed
points of $K_\var$ lie in a fixed ball in $E$ for
sufficiently small $\var>0$. Together with
\eqref{eq:puncture-map-convergence} and the homotopy
invariance of the degree, we obtain the following result.

\begin{prop}\label{prop:puncture-degree}
Under the assumptions of Proposition~\ref{prop:puncture-compactness},
for all sufficiently small $\var>0$, we have
\[
 d_\rho(\Om,H)=d_\rho(\Om_\var,H|_{\Om_\var}).
\]
\end{prop}

\begin{proof}
By the compactness theorem on $\Om$ and
Proposition~\ref{prop:puncture-compactness}, we can choose
$R>0$ such that all fixed points of $K_0$ and $K_\var$
lie in $B_R^E$ for sufficiently small $\var>0$.
Since $K_0$ is compact and has no fixed point on
$\pa B_R^E$, we have
\[
 \inf_{u\in\pa B_R^E}\norm{u-K_0(u)}_E>0.
\]
By \eqref{eq:puncture-map-convergence}, for sufficiently
small $\var>0$, the homotopy
\[
 \mathcal K(t,u)=(1-t)K_0(u)+tK_\var(u),
 \quad t\in[0,1],
\]
has no fixed point on $\pa B_R^E$.
By the homotopy invariance of the Leray--Schauder degree,
\[
 \deg(\Id-K_0,B_R^E,0)
 =\deg(\Id-K_\var,B_R^E,0).
\]
Since $K_\var(E)\subset E_\var$, the reduction property
of the Leray--Schauder degree yields
\[
 \deg(\Id-K_\var,B_R^E,0)
 =\deg(\Id-K_\var|_{E_\var},B_R^E\cap E_\var,0).
\]
Here the finite-dimensional approximations of $K_\var$
can be chosen in $E_\var$ directly, without a projection
from $E$ onto $E_\var$.

Since $K_0$ and $K_\var|_{E_\var}$ are the solution maps
on $\Om$ and $\Om_\var$, respectively, we obtain
\[
 d_\rho(\Om,H)=d_\rho(\Om_\var,H|_{\Om_\var}).
\]
\end{proof}

We next compute the degree on a disconnected domain
when all connected components share the same normalizing
integral.

Let $\Om_i$, $1\le i\le\ell$, be pairwise disjoint
smooth bounded connected domains in $\R^N$ and set
$D=\bigsqcup_{i=1}^{\ell}\Om_i$. On each $\Om_i$, let
\[
 H_i=h_i\prod_{j=1}^{m_i}|x-p_{ij}|^{N\al_{ij}},
 \quad \al_{ij}>-1,
\]
where $h_i\in C^1(\overline\Om_i)$ is positive and
$p_{ij}$ are distinct points in $\Om_i$.
Define $H$ on $D$ by $H|_{\Om_i}=H_i$ and set
\[
 E_i=W_0^{1,N}(\Om_i),\quad
 E=\prod_{i=1}^{\ell}E_i=W_0^{1,N}(D).
\]
For $u=(u_1,\ldots,u_\ell)\in E$, define
\[
 M_i(u_i)=\int_{\Om_i}H_ie^{u_i}\,\d x,\quad
 M(u)=\sum_{i=1}^{\ell}M_i(u_i),\quad
 \beta_{ij}=(1+\al_{ij})^{N-1}.
\]
Unless otherwise stated, all sums and products over
the components are taken over $1\le i\le\ell$.

The solution map on $D$ is
\[
 u\mapsto(-\Delta_N)^{-1}
\Big(\frac{c_Nr}{M(u)}He^u\Big).
\]
As in \eqref{eq:solutionmap}, we denote the corresponding
Leray--Schauder degree by $d_{c_Nr}(D,H)$.

For each $\Om_i$, define
\[
 \Si_i=\Big\{k+\sum_{j\in I}\beta_{ij}:
 k\in\N_0,\ I\subset\{1,\ldots,m_i\}\Big\}.
\]
The set $\Si_i$ contains zero and all possible
quantized masses on $\Om_i$ divided by $c_N$.
For the whole domain $D$, set
\[
 \Si_D=\{s_1+\cdots+s_\ell:
 s_i\in\Si_i,\ 1\le i\le\ell\}.
\]
Each distinct value is counted only once.
The resonance set on $D$ is
$c_N(\Si_D\setminus\{0\})$.

The compactness theorem and Proposition~\ref{prop:blowup}
also apply to disconnected domains. Indeed, on each
component $\Om_i$, the equation takes the form
\[
 -\Delta_Nu_i=bH_ie^{u_i}
 \quad\text{ in }\, \Om_i,\quad
 u_i=0\quad\text{ on }\, \pa\Om_i,
\]
where
\[
 b=\frac{c_Nr}{M(u)},\quad
 \rho_i=bM_i(u_i),\quad
 \sum_{i=1}^{\ell}\rho_i=c_Nr.
\]
Consider a sequence of solutions with bounded total masses.
If blow-up occurs on one component,
Proposition~\ref{prop:blowup} implies $b\to0$.
For each component on which the solutions remain bounded,
we have $\rho_i\to0$. For each component on which
blow-up occurs, the limiting mass is positive and
belongs to $c_N\Si_i$.
After passing to a subsequence, the limiting total mass
therefore belongs to $c_N(\Si_D\setminus\{0\})$.
In particular, compactness holds when $r\notin\Si_D$.
Moreover, the convergence $b\to0$ along blow-up sequences
does not require the nonresonance condition.

We first compare the degree of the normalized equation
with the degree of the Gelfand equation under an upper
bound on the total mass.

\begin{lem}\label{lem:small-gelfand-cutoff}
Let $D$ be a smooth bounded domain or a finite disjoint
union of smooth bounded domains, and let $r>0$ satisfy
$c_Nr\notin\Si_H$. For $\lam>0$, define
\[
 T_\lam(u)=(-\Delta_N)^{-1}(\lam He^u),\quad
 m_\lam(u)=\frac{\lam M(u)}{c_N}.
\]
For sufficiently small $\lam>0$ and sufficiently
large $R>0$, depending on $\lam$ and $r$, set
\[
 \mathcal U_{\lam,r,R}
 =\{u\in E:\norm{u}_E<R,\ m_\lam(u)<r\}.
\]
Then $T_\lam$ has no fixed point on
$\pa\mathcal U_{\lam,r,R}$, and all solutions of
\[
 -\Delta_Nu=\lam He^u\quad\text{ in }\, D,\quad
 u=0\quad\text{ on }\, \pa D
\]
with $m_\lam(u)<r$ lie in $\mathcal U_{\lam,r,R}$.
Moreover,
\be\label{eq:small-gelfand-cutoff}
 d_{c_Nr}(D,H)
 =\deg(\Id-T_\lam,\mathcal U_{\lam,r,R},0).
\ee
\end{lem}

\begin{proof}
\textbf{Step 1:} Construction of the homotopy.

Since $r$ is nonresonant, all solutions of the normalized
equation at $r$ satisfy $0\le u\le C_r$. Hence
\[
 \frac{c_Nr}{M(u)}
 \ge\frac{c_Nr}{e^{C_r}\int_D H\,\d x}
 := b_*(r)>0.
\]
Fix $0<\lam<b_*(r)$ and define
\[
 K_t(u)=(-\Delta_N)^{-1}(b_t(u)He^u),
 \quad t\in[0,1],
\]
where
\[
 b_t(u)=(1-t)\frac{c_Nr}{M(u)}
 +t\min\Big\{\lam,\frac{c_Nr}{M(u)}\Big\}.
\]
By exponential compactness and the continuity of
$(-\Delta_N)^{-1}$, the map $(t,u)\mapsto K_t(u)$
is continuous, and the family $\{K_t\}_{t\in[0,1]}$
is collectively compact on bounded subsets of $E$.

\textbf{Step 2:} Uniform estimates along the homotopy.

Let $u=K_t(u)$. By the maximum principle, $u\ge0$.
If $M(u)\ge c_Nr/\lam$, then
\[
 b_t(u)=\frac{c_Nr}{M(u)}\le\lam<b_*(r).
\]
Thus $u$ is a solution of the normalized equation
at $r$, which contradicts the definition of $b_*(r)$.

Therefore, every fixed point of $K_t$ satisfies
$M(u)<c_Nr/\lam$. In this case,
\[
 b_t(u)\ge\lam,\quad
 b_t(u)M(u)=(1-t)c_Nr+t\lam M(u)\le c_Nr.
\]
Since the total mass is bounded and $b_t(u)\ge\lam>0$,
Proposition~\ref{prop:blowup} rules out blow-up.
Hence, for fixed $\lam>0$, all fixed points of $K_t$
are uniformly bounded in $L^\infty(D)$ for $t\in[0,1]$.

Testing the equation with $u$, we obtain
\[
 \norm u_E^N
 =b_t(u)\int_D He^uu\,\d x
 \le c_Nr\norm{u}_{L^\infty(D)}\le C.
\]
Thus, for sufficiently large $R>0$, all fixed points
of $K_t$ lie in $B_R^E$. By the homotopy invariance
of the Leray--Schauder degree,
\[
 d_{c_Nr}(D,H)
 =\deg(\Id-K_0,B_R^E,0)
 =\deg(\Id-K_1,B_R^E,0).
\]

\textbf{Step 3:} The degree of the Gelfand equation.

By Step 2, every fixed point of $K_1$ satisfies
$M(u)<c_Nr/\lam$, or equivalently, $m_\lam(u)<r$.
On this open set, we have
\[
 K_1(u)=T_\lam(u).
\]
Conversely, every fixed point of $T_\lam$ with
$m_\lam(u)<r$ is a fixed point of $K_1$.
Hence all such fixed points lie in
$\mathcal U_{\lam,r,R}$ for sufficiently large $R$.

Moreover, if $u=T_\lam(u)$ and $m_\lam(u)=r$, then
\[
 \lam=\frac{c_Nr}{M(u)}<b_*(r),
\]
which contradicts the definition of $b_*(r)$.
Thus $T_\lam$ has no fixed point on
$\pa\mathcal U_{\lam,r,R}$.

By the excision property of the Leray--Schauder degree,
\[
 d_{c_Nr}(D,H)
 =\deg(\Id-K_1,B_R^E,0)
 =\deg(\Id-K_1,\mathcal U_{\lam,r,R},0)
 =\deg(\Id-T_\lam,\mathcal U_{\lam,r,R},0).
\]
This proves \eqref{eq:small-gelfand-cutoff}.
\end{proof}


For small $\lam>0$, the masses of Gelfand solutions
with bounded total mass are close to the possible
limiting masses in $\Si_i$.

\begin{lem}\label{lem:gelfand-mass-windows}
Fix a component $\Om_i$ and let $B>0$.
For any $\delta>0$, there exists $\lam_0>0$ such that,
if $0<\lam<\lam_0$ and $u$ is a Gelfand solution on
$\Om_i$ with $m_\lam(u)\le B$, then
\[
 \dist(m_\lam(u),\Si_i\cap[0,B])<\delta.
\]
Moreover, if $\lam_n\to0$ and the corresponding
solutions $u_n$ are uniformly bounded in $L^\infty(\Om_i)$,
then $m_{\lam_n}(u_n)\to0$.
\end{lem}

\begin{proof}
Suppose that the conclusion does not hold. Then there
exist $\delta>0$, $\lam_n\to0$, and solutions $u_n$
on $\Om_i$ such that
\[
 m_{\lam_n}(u_n)\le B,\quad
 \dist(m_{\lam_n}(u_n),\Si_i\cap[0,B])\ge\delta.
\]
Up to a subsequence, assume $m_{\lam_n}(u_n)\to s\in[0,B]$.

If $u_n$ is uniformly bounded in $L^\infty(\Om_i)$, then
\[
 m_{\lam_n}(u_n)=\frac{\lam_n M_i(u_n)}{c_N}\to0.
\]
Hence $s=0\in\Si_i$.

Otherwise, up to a subsequence,
$\norm{u_n}_{L^\infty(\Om_i)}\to\infty$.
Since $u_n$ satisfies the normalized equation with
parameter $\rho_n=\lam_nM_i(u_n)$,
Proposition~\ref{prop:blowup} implies
$s\in\Si_i\setminus\{0\}$.

In both cases, $s\in\Si_i\cap[0,B]$, which contradicts
the distance estimate above. The proof is complete.
\end{proof}

For any $s\in\Si_i\setminus\{0\}$, choose $0<\var<s$
such that $[s-\var,s+\var]\cap\Si_i=\{s\}$. Define
\be\label{eq:component-jumps}
 a_i(0)=1,\quad
 a_i(s)=d_{c_N(s+\var)}(\Om_i,H_i)
       -d_{c_N(s-\var)}(\Om_i,H_i).
\ee
Thus $a_i(s)$ measures the change in the degree
as the parameter passes through $c_Ns$.
By homotopy invariance of the degree, $a_i(s)$ is
independent of the choice of sufficiently small $\var$.
Since the degree at zero is one and remains constant
below the first positive resonance, we set $a_i(0)=1$.
For any $r>0$ with $r\notin\Si_i$, adding the changes
in the degree below $r$ yields
\[
 d_{c_Nr}(\Om_i,H_i)
 =\sum_{\substack{s\in\Si_i\\s<r}}a_i(s).
\]
Since $\Si_i$ is locally finite, the sum is finite.

By the product property of the Leray--Schauder degree
for the Gelfand equation, we can express the degree on
$D$ in terms of the degree changes on each component.

\begin{prop}\label{prop:disconnected-convolution}
Let $r>0$ with $r\notin\Si_D$. Then
\[
 d_{c_Nr}(D,H)
 =\sum_{\substack{s_i\in\Si_i,\ 1\le i\le\ell\\
                 s_1+\cdots+s_\ell<r}}
 \prod_{i=1}^{\ell}a_i(s_i),
\]
where $a_i(s_i)$ is defined in \eqref{eq:component-jumps}.
\end{prop}

\begin{proof}
\textbf{Step 1:} Choice of mass intervals and $\lam$.

Choose $B>r+1$ with $B\notin\bigcup_{i=1}^{\ell}\Si_i$.
Since $S_i=\Si_i\cap[0,B]$ is finite, we can choose
pairwise disjoint intervals
\[
 W_{i,0}=[0,\var_{i,0}),\quad
 W_{i,s}=(s-\var_{i,s},s+\var_{i,s}),\quad
 s\in S_i\setminus\{0\},
\]
such that $\overline W_{i,s}\subset[0,B)$ and all
positive endpoints are nonresonant. We also choose
$\var_{i,0}$ smaller than the first positive value in $\Si_i$.

Since $r\notin\Si_D$, we have
$\sum_{i=1}^{\ell}s_i\ne r$ for any  choice of $s_i\in S_i$.
Since there are only finitely many such choices,
we can choose the intervals sufficiently small so that
\be\label{eq:window-stability}
 \sum_{i=1}^{\ell}t_i<r,\quad t_i\in W_{i,s_i},
\ee
holds for all $t_i\in W_{i,s_i}$ if
$\sum_{i=1}^{\ell}s_i<r$, and for none if
$\sum_{i=1}^{\ell}s_i>r$.

We now fix $\lam>0$ sufficiently small so that
Lemma~\ref{lem:small-gelfand-cutoff} applies at $r$ on $D$
and at all positive endpoints of $W_{i,s}$ on $\Om_i$.
By Lemma~\ref{lem:gelfand-mass-windows}, we may also
assume that the mass of every Gelfand solution on $\Om_i$
with mass at most $B$ lies in one of the intervals $W_{i,s}$.

\textbf{Step 2:} The degree on each component.

For each $i$, define
\[
 T_{i,\lam}(u_i)=(-\Delta_N)^{-1}(\lam H_ie^{u_i}),
 \quad m_{i,\lam}(u_i)=\frac{\lam M_i(u_i)}{c_N}.
\]
Since $\lam>0$ is fixed, Proposition~\ref{prop:blowup}
implies that all Gelfand solutions on $\Om_i$ with
$m_{i,\lam}(u_i)\le B$ are uniformly bounded in $E_i$.

Choose $R_i>0$ sufficiently large and set
\[
 U_{i,s}=\{u_i\in E_i:\|u_i\|_{E_i}<R_i,\ 
                     m_{i,\lam}(u_i)\in W_{i,s}\}.
\]
Since $m_{i,\lam}(u_i)>0$, the set $U_{i,0}$ is also open.
By the choice of $\lam$ and $R_i$, the map $T_{i,\lam}$
has no fixed point on $\pa U_{i,s}$.

For $s>0$, Lemma~\ref{lem:small-gelfand-cutoff}
and the additivity of the degree yield
\[
 \deg(\Id-T_{i,\lam},U_{i,s},0)
 =d_{c_N(s+\var_{i,s})}(\Om_i,H_i)-d_{c_N(s-\var_{i,s})}(\Om_i,H_i)\\
 =a_i(s).
\]
For $s=0$, Lemma~\ref{lem:small-gelfand-cutoff} gives
\[
 \deg(\Id-T_{i,\lam},U_{i,0},0)
 =d_{c_N\var_{i,0}}(\Om_i,H_i)=1=a_i(0).
\]

\textbf{Step 3:} The degree on $D$.

Since the Gelfand equation on $D$ is independent on
each component, we have
\[
 T_\lam(u)=
 (T_{1,\lam}(u_1),\ldots,T_{\ell,\lam}(u_\ell)),
 \quad
 m_\lam(u)=\sum_{i=1}^{\ell}m_{i,\lam}(u_i).
\]
If $u=T_\lam(u)$ and $m_\lam(u)<r$, then
$m_{i,\lam}(u_i)<r<B$ for any  $i$.
By Step 1, each $m_{i,\lam}(u_i)$ lies in one of the
intervals $W_{i,s_i}$. By \eqref{eq:window-stability},
we have $\sum_{i=1}^{\ell}s_i<r$.

Choose $R>0$ sufficiently large such that
Lemma~\ref{lem:small-gelfand-cutoff} applies on $D$ and
\[
 \prod_{i=1}^{\ell}B_{R_i}^{E_i}\subset B_R^E.
\]
By \eqref{eq:window-stability}, the sets
$\prod_{i=1}^{\ell}U_{i,s_i}$ with
$\sum_{i=1}^{\ell}s_i<r$ are pairwise disjoint,
are contained in $\mathcal U_{\lam,r,R}$, and contain
all fixed points of $T_\lam$ in $\mathcal U_{\lam,r,R}$.

By the excision and product properties of the
Leray--Schauder degree,
\[
\deg(\Id-T_\lam,\mathcal U_{\lam,r,R},0)=
 \sum_{\substack{s_i\in S_i,\ 1\le i\le\ell\\
                 \sum_{i=1}^{\ell}s_i<r}}
 \prod_{i=1}^{\ell}
 \deg(\Id-T_{i,\lam},U_{i,s_i},0).
\]
Combining \eqref{eq:small-gelfand-cutoff} with Step 2,
we obtain
\[
 d_{c_Nr}(D,H)
 =\sum_{\substack{s_i\in\Si_i,\ 1\le i\le\ell\\
                 \sum_{i=1}^{\ell}s_i<r}}
 \prod_{i=1}^{\ell}a_i(s_i).
\]
This completes the proof.
\end{proof}

For example, if $D$ consists of two balls without
singularities, then $a_i(0)=1$, $a_i(1)=-1$, and
$a_i(s)=0$ for $s>1$. By
Proposition~\ref{prop:disconnected-convolution},
\[
 d_{c_Nr}(D,H)=1-1-1=-1,\quad 1<r<2.
\]

\begin{rem}
The degree formulas also hold for a finite disjoint
union of smooth bounded domains.
Applying Theorem~\ref{thm:degree} to each component
and using Proposition~\ref{prop:disconnected-convolution},
we obtain the degree on the whole domain.
Since
\[
 \chi(D)=\sum_{i=1}^{\ell}\chi(\Om_i),
\]
and the numbers of singular points also add over
the components, Vandermonde's identity reduces the
resulting sum to the formulas
\eqref{eq:all-nonresonant-degree} and
\eqref{eq:odd-nonresonant-degree}.
Thus the connectedness assumption is not needed
for either degree formula.
\end{rem}

Using Proposition~\ref{prop:puncture-degree} and the
Pohozaev identity, we compute the degree on a ball
with one positive singularity. The following formula
will be used in Proposition~\ref{prop:disconnected-convolution}.

\begin{lem}\label{lem:positive-local-factor}
Let $B=\{x\in\R^N:|x|<1\}$ and
$H(x)=|x|^{N\al}$ with $\al>0$.
Set $\beta=(1+\al)^{N-1}$.
For any $\rho=c_Nr>0$ with $\rho\notin\Si_H$, we have
\be\label{eq:positive-local-factor}
 d_\rho(B,H)=
 \left\{\begin{aligned}
  &1&&\text{if $N$ is even and }0<r<\beta,\\
  &0&&\text{if $N$ is even and }r>\beta,\\
  &1&&\text{if $N$ is odd and }0<r<1,\\
  &-1&&\text{if $N$ is odd and }1<r<\min\{\beta,2\},\\
  &0&&\text{if $N$ is odd and }r>\min\{\beta,2\}.
 \end{aligned}\right.
\ee
The degree remains unchanged if we change the radius
of $B$ or multiply $H$ by a positive constant.
\end{lem}
\begin{proof}
We first consider $0<r<\beta$ with $c_Nr\notin\Si_H$.
By Proposition~\ref{prop:puncture-degree}, the degree
on $B$ is equal to the degree on an annulus without
singularities. Since the Euler characteristic of
the annulus is $1-(-1)^N$, \cite[Theorem~1.1]{ATW26}
yields
\[
 d_{c_Nr}(B,H)
 =\binom{\lfloor r\rfloor-1+(-1)^N}{\lfloor r\rfloor}.
\]
For even $N$, the degree is one for all nonresonant
$r<\beta$. For odd $N$, the degree is one for $0<r<1$,
minus one for $1<r<2$, and zero for $r>2$,
provided $r<\beta$ and $c_Nr\notin\Si_H$.

We next prove that no solution exists for $r\ge\beta$.
Let $u$ be a solution and set $b=\rho/M(u)>0$,
where $\rho=c_Nr$. Applying the Pohozaev identity
with $b|x|^{N\al}(e^u-1)$ as the primitive of
$b|x|^{N\al}e^u$, we obtain
\[
 N(1+\al)\Big(\rho-b\int_B|x|^{N\al}\,\d x\Big)
 =\frac{N-1}{N}\int_{\pa B}|\pa_\nu u|^N\,\d S.
\]
Integrating the equation over $B$ yields
\[
 \rho=\int_{\pa B}|\pa_\nu u|^{N-1}\,\d S.
\]
By H\"older's inequality,
\[
 \int_{\pa B}|\pa_\nu u|^N\,\d S
 \ge S_N^{-1/(N-1)}\rho^{N/(N-1)}.
\]
Since $b\int_B|x|^{N\al}\,\d x>0$, we obtain
\[
 N(1+\al)\rho
 >\frac{N-1}{N}S_N^{-1/(N-1)}\rho^{N/(N-1)},
\]
which implies
\[
 \rho<c_N(1+\al)^{N-1}=c_N\beta.
\]
Hence no solution exists for $r\ge\beta$,
and the degree is zero for all nonresonant
$r>\beta$. No radial symmetry assumption is needed.

Combining the degree formula for $r<\beta$ with
the nonexistence result for $r\ge\beta$, we obtain
\eqref{eq:positive-local-factor}.

In particular, for even $N$, the degree changes by
$-1$ at $r=\beta$. For odd $N$, the degree changes
by $-2$ at $r=1$ and by $1$ at $r=\min\{\beta,2\}$.
There are no other changes in the degree.
The distinction between these changes and the concentration
mass is illustrated in Remark~\ref{rem:radial-sharpness}.

Finally, a change of variables shows that the degree
is independent of the radius of the ball. Multiplying
$H$ by a positive constant does not change the
normalized equation, and hence does not change the degree.
\end{proof}

\begin{rem}
The restriction \eqref{eq:puncture-mass-condition} is needed for a
uniform comparison over nonresonant parameters. In the homogeneous
ball of Lemma~\ref{lem:positive-local-factor}, choose $1<\beta<r<2$.
The strict Pohozaev inequality excludes solutions on the filled
ball. On every concentric punctured ball the weight is positive
and $C^1$. By \cite[Theorem~1.1]{ATW26}, the degree there is
$1-\chi(\mathbb S^{N-1})=(-1)^N$, so a solution exists.
If solutions on annuli with radii tending to zero were uniformly
bounded in $L^\infty$, testing the equations with the solutions
would bound them in $E$. Lemma~\ref{lem:puncture-resolvent} would
then yield a limit solution on the filled ball. Thus both uniform
compactness and degree equality can fail for $r>\beta$.
\end{rem}

\begin{rem}\label{rem:radial-sharpness}
The masses in \eqref{eq:sigma} and the nonzero degree changes
in \eqref{eq:positive-local-factor} need not occur at the same
parameters. To see this directly, keep the homogeneous ball
and the exponent $\al>0$ from
Lemma~\ref{lem:positive-local-factor}. Set
$s=q(1+\al)$ and, for $A>0$, define
\be\label{eq:radial-check}
 u_A(x)=N\ln\frac{1+A}{1+A|x|^s},\quad 
 M(u_A)=\frac{S_N(1+A)}{N(1+\al)},\quad 
 \rho_A=c_N\beta\Big(\frac{A}{1+A}\Big)^{N-1}.
\ee
The integral in \eqref{eq:radial-check} follows by the
substitution $t=A|x|^s$. Direct radial differentiation yields
\[
 -\Delta_Nu_A=b_A|x|^{N\al}e^{u_A},\quad
 b_A=\frac{N(Ns)^{N-1}(1+\al)A^{N-1}}{(1+A)^N}
     =\frac{\rho_A}{M(u_A)}.
\]
Moreover, $u_A=0$ on $\pa B$ and
$u_A\in W^{1,N}_0(B)\cap C^1(\overline B)$.
The integral of
$|\nabla u_A|^{N-2}\pa_\nu u_A$ over $\pa B_\varepsilon$
tends to zero as $\varepsilon\downarrow0$, so the equation
also holds weakly across the origin.
Consequently, \eqref{eq:radial-check} provides a solution
for any  $0<\rho<c_N\beta$, and
\[
 u_A(0)=N\ln(1+A)\to\infty,\quad
 \rho_A\uparrow c_N\beta.
\]
For odd $N$ and $\beta>2$, the degree is zero on both
sides of $c_N\beta$ at nearby nonresonant parameters,
although this family blows up as $\rho_A\uparrow c_N\beta$.
For any nonresonant $2<r<\beta$, choosing
$A/(1+A)=(r/\beta)^{1/(N-1)}$ also gives an explicit
solution with degree zero. Thus the number $2$ in
\eqref{eq:positive-local-factor} describes a change in
the annular degree; it does not replace the concentration
mass $c_N\beta$.

The strict positivity in \eqref{eq:assumptions} also has a role
in identifying the exceptional set. Choose $N\al>1$ with
$\beta\notin\N$. If zeros of the background were allowed without
counting them among the prescribed points, one could set $m=0$
and $h(x)=|x|^{N\al}$. This function is nonnegative and $C^1$.
The exceptional set computed with $m=0$ would be $c_N\N$,
whereas \eqref{eq:radial-check} blows up at $c_N\beta\notin c_N\N$.
\end{rem}

\subsection{Annular comparison of low sublevels}

Throughout this subsection, $\lam>0$ denotes the mean
field parameter. Fix $p\in(1,s_0)$ and equip all density
spaces with the strong $L^p$ topology.
We extend densities on subdomains by zero to $\Om$.

Let $D\subset\Om$ be a smooth open set containing all
the singular points, possibly with several connected
components. For $v\in W_0^{1,N}(D)$ and
$\mu\in\mathscr P_p(D)$, define
\[
 Q_D(v)=\frac1N\int_D|\nabla v|^N\,\d x,\quad 
 \mathcal E_{\lam,D}(v)
 =\frac{Q_D(v)}{\lam}-\ln\int_DHe^v\,\d x,\quad 
 \Phi_{\lam,D}(\mu)
 =\operatorname{Ent}_H(\mu)
   -\frac1\lam Q_D^*(\lam\mu),
\]
where $\mathscr P_p(D)$ is the space of probability
densities in $L^p(D)$.

Let $\Theta$ and $\mathscr S_\lam$ be the maps defined
in Lemma~\ref{lem:nonlinear-entropy-duality}.
By the same lemma, we have
\be\label{eq:cech-duality-bounds}
 \Phi_{\lam,\Om}(\Theta u)\le\mathcal E_{\lam,\Om}(u),
 \quad
 \mathcal E_{\lam,\Om}(\mathscr S_\lam\mu)
 \le\Phi_{\lam,\Om}(\mu).
\ee
Moreover, for the homotopy
\[
 \mu_t=(1-t)\mu+t\Theta(\mathscr S_\lam\mu),
 \quad t\in[0,1],
\]
we have
\[
 \Phi_{\lam,\Om}(\mu_t)\le\Phi_{\lam,\Om}(\mu),
 \quad t\in[0,1].
\]
Thus the homotopy preserves any sublevel set of
$\Phi_{\lam,\Om}$.

Fix a closed interval
$I\Subset(0,\infty)\setminus\Si_H$ with nonempty interior.
In the following sublevel comparisons, all parameters
are taken in $\operatorname{int}I$.

For $\lam<\lam'$ in $\operatorname{int}I$, we call a
continuous map controlled if it maps the sublevel set
at level $-L-C$ for parameter $\lam$ into the sublevel
set at level $-L$ for parameter $\lam'$, where $C$
is independent of $L$.
A controlled homotopy remains in the sublevel set
at level $-L$ for parameter $\lam'$ for all $t\in[0,1]$.
All sublevel sets below are defined using the normalized
functionals $\mathcal E$ and $\Phi$.

For fixed $\lam<\lam'$, Lemma~\ref{lem:small-mass-trimming}
implies
\be\label{eq:cech-normalized-trimming}
 \Phi_{\lam',\Om}(T_w\mu)
 \le\Phi_{\lam,\Om}(\mu)+C_\eta,
\ee
provided
\[
 \Phi_{\lam,\Om}(\mu)\le0,\quad
 \int_\Om(1-w)\mu\,\d x\le\eta,
\]
where $\eta>0$ is sufficiently small depending on
$\lam$ and $\lam'$. For a family of cutoff functions
$w_t$, estimate \eqref{eq:cech-normalized-trimming}
holds uniformly for $t\in[0,1]$ if
\[
 \int_\Om(1-w_t)\mu\,\d x\le\eta,
 \quad t\in[0,1].
\]

By \eqref{eq:cech-normalized-trimming}, we have an
energy estimate after removing a small amount of mass.
We next study how the energy changes under a
diffeomorphism of the domain.

\begin{lem}
Let $D,D_0\subset\Om$ be bounded open sets and let
$T:D\to D_1\subset D_0$ be a $C^1$ diffeomorphism.
Assume that $DT$ and $(DT)^{-1}$ are bounded and
\[
 L_T=\sup_{x\in D}|\det DT(x)|\,\|(DT(x))^{-1}\|^N
 <\infty,
\]
\[
 \operatorname*{ess\,sup}_{x\in D}
 \Big|\ln\frac{H(Tx)|\det DT(x)|}{H(x)}\Big|
 \le C_T.
\]
Then, for any nonnegative $\nu\in L^p(D)$,
\[
 Q_{D_0}^*(T_\#\nu)
 \ge L_T^{-1/(N-1)}Q_D^*(\nu),
\]
where $T_\#\nu$ denotes the density obtained from
$\nu$ by $T$ and extended by zero to $D_0$.

Moreover, if $\lam'/\lam\ge L_T$, then for any 
$\mu\in\mathscr P_p(D)$,
\be\label{eq:cech-transport-bound}
 \Phi_{\lam',D_0}(T_\#\mu)
 \le\Phi_{\lam,D}(\mu)+C_T.
\ee
Both estimates are uniform along any $C^1$-continuous
family of diffeomorphisms for which $DT$, $(DT)^{-1}$,
$L_T$, and $C_T$ are uniformly bounded.
\end{lem}

\begin{proof}
For any $v\in W^{1,N}_0(D)$, we have
$v\circ T^{-1}\in W^{1,N}_0(D_1)$.
By a change of variables,
\[
 Q_{D_1}(v\circ T^{-1})
 \le\frac1N\int_D|\nabla v(x)|^N
 \|(DT(x))^{-1}\|^N|\det DT(x)|\,\d x
 \le L_TQ_D(v).
\]
Moreover,
\[
 \int_{D_1}(v\circ T^{-1})T_\#\nu\,\d x
 =\int_Dv\nu\,\d x.
\]
By the definition of the convex conjugate,
\[
 Q_{D_1}^*(T_\#\nu)
 \ge\sup_{v\in W^{1,N}_0(D)}
 \Big\{\int_Dv\nu\,\d x-L_TQ_D(v)\Big\}
 =L_TQ_D^*(\nu/L_T)
 =L_T^{-1/(N-1)}Q_D^*(\nu),
\]
where the last equality follows from the homogeneity
of $Q_D^*$. Since $D_1\subset D_0$, extending functions
by zero from $D_1$ to $D_0$ yields
\[
 Q_{D_0}^*(T_\#\nu)
 \ge Q_{D_1}^*(T_\#\nu)
 \ge L_T^{-1/(N-1)}Q_D^*(\nu).
\]

We next compare the entropy terms.
By a change of variables,
\[
 \operatorname{Ent}_H(T_\#\mu)
 =\operatorname{Ent}_H(\mu)
 -\int_D\mu(x)\ln
 \frac{H(Tx)|\det DT(x)|}{H(x)}\,\d x.
\]
Since $\int_D\mu\,\d x=1$, we obtain
\[
 \operatorname{Ent}_H(T_\#\mu)
 \le\operatorname{Ent}_H(\mu)+C_T.
\]
Moreover, by the preceding estimate and the homogeneity
of $Q_D^*$,
\[
 \frac1{\lam'}Q_{D_0}^*(\lam'T_\#\mu)
 \ge\frac1{\lam}
 \Big(\frac{\lam'}{\lam L_T}\Big)^{1/(N-1)}
 Q_D^*(\lam\mu)\ge\frac1{\lam}Q_D^*(\lam\mu),
\]
where we have used $\lam'/\lam\ge L_T$.
Combining the last two estimates yields
\eqref{eq:cech-transport-bound}.

Finally, along a $C^1$-continuous family of
diffeomorphisms, the uniform estimates on $DT$
and $(DT)^{-1}$ imply uniform boundedness of the
maps $\nu\mapsto T_\#\nu$ from $L^p(D)$ to $L^p(D_0)$.
Strong $L^p$ continuity follows first for smooth
compactly supported densities by a change of variables,
and then for general densities by approximation.
\end{proof}

Choose a singular point $p_i$ and $R>0$ such that
$\overline B_R(p_i)\Subset\Om$ contains no other
singular point. Let $Z_1,\ldots,Z_M$ be pairwise
disjoint closed annuli centered at $p_i$ and contained
in $B_R(p_i)$, with positive inner radii and positive
distances between them. We index the annuli from
the innermost to the outermost.

For each $Z_j$, choose nested annuli satisfying
\be\label{eq:cech-nested-annuli}
 Z'_j\Subset C_j\Subset\operatorname{int}Z_j
 \subset Z_j\Subset G_j\Subset O_j,
\ee
where $Z'_j$ is closed and $C_j$, $G_j$, and $O_j$
are open. We require that the annuli $O_j$ have
positive pairwise distances. All these annuli can
be chosen arbitrarily close to $Z_j$ by slightly
changing its inner and outer radii.

For any nonempty subset $J\subset\{1,\ldots,M\}$, set
\[
 D_J=\Om\setminus\bigcup_{j\in J}Z_j,\quad
 D'_J=\Om\setminus\bigcup_{j\in J}Z'_j.
\]
Both $D_J$ and $D'_J$ are smooth open sets containing
all the singular points, possibly with several
connected components. Moreover,
\[
 D_J\subset D'_J,\quad
 D_K\subset D_J\quad\text{if }J\subset K.
\]
For these inclusions, we identify densities and
Sobolev functions on the smaller domain with their
zero extensions to the larger domain.

Using \eqref{eq:cech-transport-bound}, we compare
the low sublevel sets on domains obtained by removing
annuli of slightly different widths.

\begin{lem}\label{lem:cech-cut-comparison}
Fix $a<b$ in $\operatorname{int}I$.
The annuli in \eqref{eq:cech-nested-annuli}
can be chosen such that the following conclusions hold.

\begin{enumerate}
\item[(i)] For any nonempty $J\subset\{1,\ldots,M\}$
and $\lam\in I$, if $\mu\in\mathscr P_p(\Om)$
satisfies $\mu=0$ on $G_j$ for any  $j\in J$, then
\[
 0\le\Phi_{\lam,D_J}(\mu)
 -\Phi_{\lam,\Om}(\mu)\le C.
\]
Similarly, if $\mu=0$ on $C_j$ for any  $j\in J$, then
\[
 0\le\Phi_{\lam,D'_J}(\mu)
 -\Phi_{\lam,\Om}(\mu)\le C.
\]
The constant $C>0$ is independent of $J$, $\lam$, and $\mu$.

\item[(ii)] For all sufficiently large $L$, the zero extension
\[
 \Phi_{b,D_J}^{-L}\to\Phi_{b,D'_J}^{-L}
\]
induces an isomorphism in homology for any  nonempty
$J\subset\{1,\ldots,M\}$. The choice of $L$ is independent
of $J$.
\end{enumerate}
\end{lem}

\begin{proof}
\textbf{Step 1:} Energy estimates on $D_J$ and $D'_J$.

If $\mu=0$ on $G_j$ for any  $j\in J$, then $\mu$
vanishes in a fixed neighborhood of each inner boundary
of $D_J$. Applying Lemma~\ref{lem:buffered-potential-comparison}
with $\nu=\lam\mu$, we obtain
\[
 0\le Q_\Om^*(\lam\mu)-Q_{D_J}^*(\lam\mu)
 \le C\lam^{N/(N-1)}.
\]
Since zero extension does not change the entropy,
\[
 0\le\Phi_{\lam,D_J}(\mu)-\Phi_{\lam,\Om}(\mu)
 =\frac1\lam
(Q_\Om^*(\lam\mu)-Q_{D_J}^*(\lam\mu)) 
 \le C\lam^{1/(N-1)}.
\]
Similarly, if $\mu=0$ on $C_j$ for any  $j\in J$,
then
\[
 0\le\Phi_{\lam,D'_J}(\mu)-\Phi_{\lam,\Om}(\mu)
 \le C\lam^{1/(N-1)}.
\]
Since $I$ is compact and there are only finitely many
choices of $J$, both estimates hold with a uniform constant.

\textbf{Step 2:} Construction of the homotopies.

In radial coordinates about $p_i$, write the radial intervals
of $Z_j$ and $Z'_j$ as
$[r_j^-,r_j^+]$ and $[r_j^-+\varepsilon_j^-,r_j^+-\varepsilon_j^+]$,
respectively. Choose smooth radial cutoffs $\chi_j^\pm$ equal
to one near the corresponding endpoints of $Z_j$ and $Z'_j$.
Take their supports in disjoint collars, away from the
singular points and $\pa\Om$. Define
\[
 R_{J,t}(s)=s+t\sum_{j\in J}
   \big(-\varepsilon_j^-\chi_j^-(s)
       +\varepsilon_j^+\chi_j^+(s)\big),
 \quad
 T_{J,t}(p_i+s\omega)=p_i+R_{J,t}(s)\omega .
\]
For small $\varepsilon_j^\pm$, the map is a diffeomorphism
of $\Om$ and is the identity near the singular points and
$\pa\Om$. Indeed $R_{J,t}'>0$ and the map is the identity
outside the collars. The image of each removed radial interval
contains that interval, for both $Z_j$ and $Z'_j$. At $t=1$
the image of $Z'_j$ is exactly $Z_j$. Taking complements gives
\[
 T_{J,0}=\Id,\quad T_{J,1}(D'_J)=D_J,\quad
 T_{J,t}(D_J)\subset D_J,\quad T_{J,t}(D'_J)\subset D'_J .
\]
The radial and tangential stretches $R_{J,t}'$ and $R_{J,t}/s$
converge uniformly to one as all $\varepsilon_j^\pm\to0$.
Thus the collars can be chosen so that
$\sup_{J,t} L_{T_{J,t}}<b/a$ in \eqref{eq:cech-transport-bound}.
Since $T_{J,t}$ is the identity near all singular
points, there exists $C_{\mathrm{tr}}\ge0$ such that
\[
 \operatorname*{ess\,sup}_{x\in D'_J}
 \Big|
 \ln\frac{H(T_{J,t}x)|\det DT_{J,t}(x)|}{H(x)}
 \Big|
 \le C_{\mathrm{tr}}
\]
for any  $J$ and $t\in[0,1]$.

By \eqref{eq:cech-transport-bound}, the map
\[
 R_J(\mu)=(T_{J,1})_\#\mu
\]
is continuous from $\Phi_{a,D'_J}^{-L-C_{\mathrm{tr}}}$
to $\Phi_{b,D_J}^{-L}$.
Let
\[
 i_a:\Phi_{a,D_J}^{-L-C_{\mathrm{tr}}}
 \to\Phi_{a,D'_J}^{-L-C_{\mathrm{tr}}},\quad 
 i_b:\Phi_{b,D_J}^{-L}
 \to\Phi_{b,D'_J}^{-L}
\]
be the zero extension maps. We also define the natural
inclusions
\[
 c_J:\Phi_{a,D_J}^{-L-C_{\mathrm{tr}}}
 \to\Phi_{b,D_J}^{-L},\quad 
 c'_J:\Phi_{a,D'_J}^{-L-C_{\mathrm{tr}}}
 \to\Phi_{b,D'_J}^{-L}.
\]
The family $(T_{J,t})_\#\mu$, together with
\eqref{eq:cech-transport-bound}, yields
\[
 R_J\circ i_a\simeq c_J,\quad
 i_b\circ R_J\simeq c'_J.
\]
The first homotopy takes values in
$\Phi_{b,D_J}^{-L}$ and the second in
$\Phi_{b,D'_J}^{-L}$.

\textbf{Step 3:} Homology isomorphisms.

We first compare the sublevel sets at different
parameters and levels on each fixed domain
$D\in\{\Om,D_J,D'_J\}$.

For a disconnected $D$, let $B$ run over its connected
components. The normalized component mass of a sequence
of solutions at parameters $\lam_n\in I$ is
\[
 \lam_{n,B}=b_n\int_B He^{u_n}\,\d x,\quad
 b_n=\frac{\lam_n}{\int_DHe^{u_n}\,\d x},\quad
 \sum_B\lam_{n,B}=\lam_n .
\]
On a component where blow-up occurs,
Proposition~\ref{prop:blowup} applied with parameter
$\lam_{n,B}$ gives $b_n\to0$ and quantizes the limiting
component mass. On each bounded component the same $b_n$
gives $\lam_{n,B}\to0$. After taking a subsequence on the
finitely many components, a blow-up limit of the total
parameter thus lies in the sum of the component mass sets,
each augmented by zero. All prescribed points are retained
exactly once and the integer masses can be assigned to any
one component. Hence this sum, after removing zero, is
exactly $\Si_H$ in \eqref{eq:sigma}. This contradicts
$\lam_n\in I$ and proves the compactness on $D$ needed below.

Let $\lambda<\lambda'$ lie in $\operatorname{int}I$, and choose
$\lambda_0\in\operatorname{int}I$ below $\lambda$. Set
$\vartheta=\lambda'/\lambda$. Then
\[
 \vartheta J_{\lambda,D}
   =J_{\lambda',D}+(\vartheta-1)Q_D,
 \quad
 \frac{-\lambda'L}{\vartheta}=-\lambda L.
\]
Thus the common descent argument of
\cite[Proposition~3.5]{ATW26}, with lower parameter $\lambda_0$,
shows that the inclusion
\[
 \mathcal E_{\lambda,D}^{-L}
 \hookrightarrow \mathcal E_{\lambda',D}^{-L}
\]
is a homotopy equivalence for sufficiently large $L$. The required inputs are compactness of the
critical points for parameters in $[\lambda_0,\lambda']$ and compactness
of bounded Palais--Smale sequences, as in
Proposition~\ref{prop:all-parameter-euler}. In the same way,
\cite[Proposition~3.3]{ATW26} applies to the closed strip between two
sufficiently low levels. In particular, for each fixed $C\ge0$,
\[
 \mathcal E_{\lambda,D}^{-L-C}
 \hookrightarrow \mathcal E_{\lambda,D}^{-L}
\]
is a homotopy equivalence for sufficiently large $L$.

By Lemma~\ref{lem:nonlinear-entropy-duality}, the map
\[
 \Theta_D(v)=\frac{He^v}{\int_DHe^v\,\d x}
\]
induces a homotopy equivalence between the corresponding
sublevel sets of $\mathcal E_{\lam,D}$ and $\Phi_{\lam,D}$.
Since $\Theta_D$ is independent of the parameter
and the sublevel value, the same homotopy equivalences
hold for the inclusions
\[
 \Phi_{\lam,D}^{-L}
 \hookrightarrow\Phi_{\lam',D}^{-L},
 \quad
 \Phi_{\lam,D}^{-L-C}
 \hookrightarrow\Phi_{\lam,D}^{-L}.
\]
Consequently, $(c_J)_*$ and $(c'_J)_*$ are isomorphisms
for sufficiently large $L$.

By Step 2, we have
\[
 (R_J)_*(i_a)_*=(c_J)_*,\quad
 (i_b)_*(R_J)_*=(c'_J)_*.
\]
Since $(c_J)_*$ is an isomorphism, $(i_a)_*$ is
injective. Since $(c'_J)_*$ is an isomorphism,
$(i_b)_*$ is surjective. Moreover,
\[
 (i_b)_*(c_J)_*=(c'_J)_*(i_a)_*.
\]
Thus $(i_b)_*$ is also injective and hence is an
isomorphism.

Since there are only finitely many choices of $J$,
the same sufficiently large $L$ works for any  $J$.
Finally, for $J\subset K$, the zero extension maps
commute with the inclusions induced by
$D_K\subset D_J$ and $D'_K\subset D'_J$.
\end{proof}

If a density vanishes on an annulus, the corresponding
potential is $N$-harmonic there. We use this fact to
estimate the Gibbs mass on smaller annuli.

\begin{lem}\label{lem:cech-harmonic-buffer}
Let $\lam\in I$ and $\mu\in\mathscr P_p(\Om)$ satisfy
\[
 \Phi_{\lam,\Om}(\mu)\le-L,\quad
 \mu=0\quad\text{ on }\, Z_j,\quad j\in J.
\]
Set $u=\mathscr S_\lam\mu$ and $\nu=\Theta(u)$.
For each $j\in J$, choose a smooth cutoff function
$\psi_j$ such that
\[
 0\le\psi_j\le1,\quad
 \psi_j=1\quad\text{ on }\, C_j,\quad
 \supp\psi_j\Subset\operatorname{int}Z_j.
\]
Then
\[
 \nu\Big(\bigcup_{j\in J}\supp\psi_j\Big)
 \le Ce^{-L},
\]
where $C>0$ is independent of $\mu$, $L$, $J$, and
$\lam\in I$.
\end{lem}

\begin{proof}
Since $\mu\in\mathscr P_p(\Om)$ and $\lam\in I$, we have
\[
 \int_\Om\lam\mu\,\d x=\lam\le\max I.
\]
By the global estimates for equations with measure data,
\[
 \norm{u}_{L^a(\Om)}\le C_a
\]
for any  finite $a\ge1$.
Since $\mu=0$ on $Z_j$ for $j\in J$, we have
\[
 -\Delta_Nu=0\quad\text{ in }\, \operatorname{int}Z_j.
\]
By the local estimates for $N$-harmonic functions and
$\supp\psi_j\Subset\operatorname{int}Z_j$, we obtain
\[
 \sup_{\supp\psi_j}u\le C,\quad j\in J.
\]
Since $\supp\psi_j$ contains no singular points,
$H$ is bounded on $\supp\psi_j$. Consequently,
\[
 \sum_{j\in J}\int_{\supp\psi_j}He^u\,\d x\le C.
\]

On the other hand, by \eqref{eq:cech-duality-bounds},
\[
 \frac{Q_\Om(u)}{\lam}-\ln M(u)
 =\mathcal E_{\lam,\Om}(u)
 \le\Phi_{\lam,\Om}(\mu)\le-L.
\]
Since $Q_\Om(u)\ge0$, we have $M(u)\ge e^L$.
Therefore,
\[
 \nu\Big(\bigcup_{j\in J}\supp\psi_j\Big)
 =\frac{
 \int_{\bigcup_{j\in J}\supp\psi_j}He^u\,\d x}{M(u)}
 \le Ce^{-L}.
\]
\end{proof}

Combining Lemmas~\ref{lem:cech-cut-comparison}
and \ref{lem:cech-harmonic-buffer}, we can express
the Euler characteristic of low sublevel sets on $\Om$
in terms of those on the domains $D_J$.

\begin{prop}\label{prop:annular-cech-excision}
Fix $\rho\in\operatorname{int}I$.
For all sufficiently large $M$, the annuli
$Z_1,\ldots,Z_M$ can be chosen such that
\be\label{eq:annular-cech-Euler}
 \chi(J_{\rho,\Om}^{-L})
 =\sum_{\varnothing\ne J\subset\{1,\ldots,M\}}
 (-1)^{|J|-1}\chi(J_{\rho,D_J}^{-L_J})
\ee
for sufficiently large $L$ and $L_J$.

On each $D_J$, we use the functional
\[
 J_{\rho,D_J}(v)
 =\sum_B Q_B(v_B)
 -\rho\ln\Big(\sum_B\int_B He^{v_B}\,\d x\Big),
\]
where both sums are taken over all connected components
$B$ of $D_J$ and $v_B=v|_B$.

In terms of the Leray--Schauder degree,
\eqref{eq:annular-cech-Euler} is equivalent to
\be\label{eq:annular-cech-degree}
 d_{\rho,\Om}
 =\sum_{\varnothing\ne J\subset\{1,\ldots,M\}}
 (-1)^{|J|-1}d_{\rho,D_J}.
\ee
Here $d_{\rho,D_J}$ is the degree on the whole domain
$D_J$, with a single normalizing integral over all
its connected components.
\end{prop}

\begin{proof}
\textbf{Step 1:} Choice of parameters, annuli, and cutoff functions.

Fix $a<\rho<b$ in $\operatorname{int}I$ and choose
$\eta>0$ as in Lemma~\ref{lem:small-mass-trimming}
for the parameters $(a,b)$.

Set $\bar\lam=\max I$. Since the number of points
in the support of each weighted barycenter at
$\bar\lam$ is bounded, Lemma~\ref{lem:concentration}
implies that there exists an integer $K\ge1$ with
the following property: for any $d>0$ and $\var>0$,
if $\mathcal E_{\bar\lam,\Om}(u)$ is sufficiently
negative, then there exist $x_1,\ldots,x_k\in\overline\Om$,
with $k\le K$, such that
\[
 \Theta(u)\Big(\bigcup_{i=1}^kB_d(x_i)\Big)
 \ge1-\var.
\]
The barycenter space assigns weight one to the positive
singularities. If this space is empty, the deep sublevels are empty
and any $K\ge1$ suffices.

Since
\[
 \mathcal E_{\bar\lam,\Om}(u)
 \le\mathcal E_{\lam,\Om}(u),\quad \lam\in I,
\]
the same concentration estimate holds uniformly for
all $\lam\in I$.

Choose $M>16K/\eta$ and set $\delta=2/M$.
By Lemma~\ref{lem:cech-cut-comparison}, choose the
annuli in \eqref{eq:cech-nested-annuli} for the
parameters $(a,b)$, with positive distances between
the outer annuli $O_j$. All comparison homotopies
below take values in sublevel sets at parameter $b$.

Choose
\[
 0<d<\frac14\min_{i\ne j}\dist(O_i,O_j),
 \quad 0<\var<\frac{\eta}{8}.
\]
For $u$ in a sufficiently low sublevel set, suppose
$A\subset\{1,\ldots,M\}$ satisfies
\[
 \Theta(u)(O_j)<\delta,\quad j\in A.
\]
Each ball of radius $d$ meets at most one of the
annuli $O_j$. Thus at most $K$ annuli meet the
concentration balls, and the total mass outside
the concentration balls is at most $\var$. Hence
\be\label{eq:cech-union-small-mass}
 \Theta(u)\Big(\bigcup_{j\in A}O_j\Big)
 \le K\delta+\var<\frac{\eta}{4}.
\ee

Let $C_{\mathrm{buf}}\ge0$ be the constant in the
energy estimates of Lemma~\ref{lem:cech-cut-comparison},
and let $C_{\mathrm{trim}}\ge0$ be the constant in
\eqref{eq:cech-normalized-trimming} for $(a,b)$ and $\eta$.
Set
\[
 C=C_{\mathrm{buf}}+C_{\mathrm{trim}}.
\]
Fix the annuli and $C$. Then choose $L$ large enough for
the concentration estimate and for Lemmas~\ref{lem:cech-cut-comparison}
and~\ref{lem:cech-harmonic-buffer} to apply uniformly.

We use a partition of unity subordinate to the open sets
\[
 U_j=\{u\in E:\Theta(u)(O_j)<\delta\},
 \quad 1\le j\le M.
\]
Since the annuli $O_j$ are pairwise disjoint and
$M\delta=2$, we have
\[
 \sum_{j=1}^M\Theta(u)(O_j)\le1<\frac32<M\delta.
\]
Hence every $u\in E$ lies in at least one $U_j$.
Define
\[
 p_j(u)=
 \frac{(3/(2M)-\Theta(u)(O_j))_+}
 {\sum_{i=1}^M(3/(2M)-\Theta(u)(O_i))_+}.
\]
The preceding inequality makes the denominator positive.
Each $p_j$ is continuous and $\sum_{j=1}^Mp_j(u)=1$.
Moreover, $\supp p_j\subset\{u:\Theta(u)(O_j)\le3/(2M)\}\subset U_j$.

Let
\[
 \Delta=\Big\{\xi\in[0,1]^M:
 \sum_{j=1}^M\xi_j=1\Big\}
\]
be the standard $(M-1)$-simplex, and let $\Delta_J$
be the face spanned by the vertices indexed by $J$.
Set $\tau=1/(2M)$ and define
\[
 r_\tau(\xi)_j=
 \frac{(\xi_j-\tau)_+}
 { \sum_{i=1}^M(\xi_i-\tau)_+},
 \quad
 b_j(\xi)=\min\{1,\xi_j/\tau\}.
\]
Since
\[
 \sum_{j=1}^M(\xi_j-\tau)_+
 \ge1-M\tau=\frac12,
\]
the map $r_\tau:\Delta\to\Delta$ is well defined
and continuous. Moreover,
\[
 r_\tau(\xi)_j>0
 \quad\Longrightarrow\quad b_j(\xi)=1.
\]

\textbf{Step 2:} Construction and homology of the
\v{C}ech realization.

For each nonempty $J\subset\{1,\ldots,M\}$, consider
the sublevel set $\Phi_{\lam,D_J}^{-L}$ and define
\[
 \mathcal W_{\lam,L}
 =\Big(\coprod_{J\ne\varnothing}
 \Phi_{\lam,D_J}^{-L}\times\Delta_J\Big)\big/\sim.
\]
For each nonempty $A\subset J$, we identify points on
$\Delta_A\subset\Delta_J$ by extending the density
by zero from $D_J$ to $D_A$. By domain monotonicity,
the zero extension does not increase $\Phi$.

Each point of $\mathcal W_{\lam,L}$ can be written
as $(\mu,\xi)$, where $\xi\in\Delta$ and
\[
 \mu\in\Phi_{\lam,D_J}^{-L},\quad
 J=\{j:\xi_j>0\}.
\]
We regard all densities as functions in $L^p(\Om)$
by zero extension. Since the conditions
\[
 \mu\ge0,\quad \int_\Om\mu\,\d x=1,\quad
 \mu=0\quad\text{ on }\,\Om\setminus D_J
\]
are closed in $L^p(\Om)$, continuity of
$\Phi_{\lam,D_J}$ implies that
$\Phi_{\lam,D_J}^{-L}$ is closed in $L^p(\Om)$.
Hence
\[
 \mathcal W_{\lam,L}
 =\bigcup_{J\ne\varnothing}
(\Phi_{\lam,D_J}^{-L}\times\Delta_J) 
\]
is a finite union of closed subsets of
$L^p(\Om)\times\Delta$.
Therefore, the quotient topology agrees with the
subspace topology. Similarly, we define
$\mathcal W'_{\lam,L}$ by replacing $D_J$ with $D'_J$.

We next compute the homology of $\mathcal W_{\lam,L}$.
No contractibility assumption on the density sublevels is
needed: the skeletal calculation retains their homology.
Let $\mathcal W^{(k)}$ be its $k$-skeleton and set
$\mathcal W^{(-1)}=\varnothing$. To obtain
$\mathcal W^{(k)}$ from $\mathcal W^{(k-1)}$, we attach
the products
\[
 \Phi_{\lam,D_J}^{-L}\times\Delta_J,\quad |J|=k+1,
\]
along their boundaries in the simplex variable,
\[
 \Phi_{\lam,D_J}^{-L}\times\pa\Delta_J.
\]
For $k\ge1$, the product inclusions are closed cofibrations,
since $(\Delta_J,\pa\Delta_J)$ is a neighborhood deformation
retract pair. Pushing their disjoint union out along the attaching
maps makes $\mathcal W^{(k-1)}\hookrightarrow\mathcal W^{(k)}$
a cofibration as well. The quotient by the previous skeleton
is the wedge of the corresponding product quotients (with a
disjoint base point for any empty lower member).
The relative K\"unneth formula over $\mathbb Q$ then gives
\be\label{eq:cech-skeletal-homology}
 H_i(\mathcal W^{(k)},\mathcal W^{(k-1)};\mathbb Q)
 \cong
\oplus_{|J|=k+1}
 H_{i-k}(\Phi_{\lam,D_J}^{-L};\mathbb Q).
\ee
For $k=0$, the formula follows directly from the
homology of a finite disjoint union.
The isomorphisms in \eqref{eq:cech-skeletal-homology}
are compatible with inclusions of domains, sublevel
sets, and parameters.

By Proposition~\ref{prop:all-parameter-euler} and
Lemma~\ref{lem:nonlinear-entropy-duality}, we have
\[
 \sum_{i\ge0}
 \dim H_i(\Phi_{\lam,D_J}^{-L};\mathbb Q)<\infty
\]
for any  nonempty $J$ and sufficiently large $L$.
Using \eqref{eq:cech-skeletal-homology} and the long
exact sequence of each pair
$(\mathcal W^{(k)},\mathcal W^{(k-1)})$, we obtain
\[
 \sum_{i\ge0}
 \dim H_i(\mathcal W_{\lam,L};\mathbb Q)<\infty.
\]

Moreover, if a compatible family of maps on the
density sublevel sets induces homology isomorphisms
for any  $J$, then induction on $k$ and the five
lemma show that the induced map on the realizations
is also a homology isomorphism.

For sufficiently large $L$, we apply this argument
to the inclusions between different parameters and
sublevel sets. By Lemma~\ref{lem:cech-cut-comparison},
the same argument shows that the zero extension
\[
 h_b:\mathcal W_{b,L}\to\mathcal W'_{b,L}
\]
induces an isomorphism in homology.

\textbf{Step 3:} Construction of the maps $F$ and $G_\lam$.

For each $j$, choose a smooth cutoff function $\theta_j$
such that
\[
 0\le\theta_j\le1,\quad
 \theta_j=1\quad\text{ on }\,G_j,\quad
 \supp\theta_j\Subset O_j.
\]
For $\xi\in\Delta$, define
\[
 T^\theta_\xi\mu
 =\frac{\prod_{j=1}^M(1-b_j(\xi)\theta_j)\mu}
 { \int_\Om
 \prod_{j=1}^M(1-b_j(\xi)\theta_j)\mu\,\d x},
\]
whenever the denominator is positive.

Set $X_{\lam,L}=\mathcal E_{\lam,\Om}^{-L}$.
For $u\in X_{a,L+C}$, define
\[
 F(u)=\big(T^\theta_{p(u)}\Theta(u),r_\tau(p(u))) .
\]
By \eqref{eq:cech-union-small-mass}, the mass removed
from $\Theta(u)$ is at most $\eta/4$.
In particular, the denominator in
$T^\theta_{p(u)}\Theta(u)$ is bounded below by
$1-\eta/4>0$.

Set
\[
 J(u)=\{j:r_\tau(p(u))_j>0\}.
\]
For every $j\in J(u)$, we have $b_j(p(u))=1$.
Since $\theta_j=1$ on $G_j$,
\[
 T^\theta_{p(u)}\Theta(u)=0
 \quad\text{ on }\,G_j,\quad j\in J(u).
\]
By Lemma~\ref{lem:cech-cut-comparison},
\eqref{eq:cech-normalized-trimming}, and
\eqref{eq:cech-duality-bounds}, we obtain
\[
 \begin{aligned}
 \Phi_{b,D_{J(u)}}(T^\theta_{p(u)}\Theta(u))
 &\le\Phi_{b,\Om}(T^\theta_{p(u)}\Theta(u))
      +C_{\mathrm{buf}}\\
 &\le\Phi_{a,\Om}(\Theta(u))
      +C_{\mathrm{trim}}+C_{\mathrm{buf}}\\
 &\le\mathcal E_{a,\Om}(u)+C\\
 &\le-L.
 \end{aligned}
\]
Thus
\[
 F:X_{a,L+C}\to\mathcal W_{b,L}
\]
is well defined. Since $p$, $r_\tau$, and $b_j$
are continuous, and the denominator is bounded away
from zero, $F$ is continuous with respect to the
subspace topology established in Step 2.

We next define
\[
 G_\lam:\mathcal W_{\lam,L}\to X_{\lam,L},
 \quad G_\lam(\mu,\xi)=\mathscr S_\lam\mu.
\]
For $(\mu,\xi)\in\mathcal W_{\lam,L}$, set
$J=\{j:\xi_j>0\}$. By domain monotonicity and
\eqref{eq:cech-duality-bounds},
\[
 \mathcal E_{\lam,\Om}(G_\lam(\mu,\xi))
 \le\Phi_{\lam,\Om}(\mu)
 \le\Phi_{\lam,D_J}(\mu)\le-L.
\]
Hence $G_\lam$ is well defined and continuous.

We now prove that $G_b\circ F$ is homotopic to the
inclusion $X_{a,L+C}\hookrightarrow X_{b,L}$.
For $u\in X_{a,L+C}$, set $\mu=\Theta(u)$ and consider
\[
 v_t=(1-t)u+t\mathscr S_b\mu,\quad t\in[0,1].
\]
Define
\[
 \mathcal L_{b,\Om}(v,\mu)
 =\frac{Q_\Om(v)}b-\int_\Om v\mu\,\d x
  +\operatorname{Ent}_H(\mu).
\]
Since $\mathscr S_b\mu$ minimizes
$\mathcal L_{b,\Om}(\cdot,\mu)$, convexity and
Lemma~\ref{lem:nonlinear-entropy-duality} imply
\[
 \mathcal E_{b,\Om}(v_t)
 \le\mathcal L_{b,\Om}(v_t,\mu)
 \le\mathcal L_{b,\Om}(u,\mu)
 =\mathcal E_{b,\Om}(u)
 \le\mathcal E_{a,\Om}(u)\le-L-C.
\]
Thus the line segment from $u$ to $\mathscr S_b\mu$
lies in $X_{b,L}$.

Next, for $t\in[0,1]$, define
\[
 w_t=\prod_{j=1}^M(1-tb_j(p(u))\theta_j),\quad
 \mu_t=\frac{w_t\mu}{\int_\Om w_t\mu\,\d x}.
\]
By \eqref{eq:cech-union-small-mass},
\[
 \int_\Om(1-w_t)\mu\,\d x\le\frac{\eta}{4},
 \quad t\in[0,1].
\]
Applying \eqref{eq:cech-normalized-trimming} and
\eqref{eq:cech-duality-bounds}, we obtain
\[
 \mathcal E_{b,\Om}(\mathscr S_b\mu_t)
 \le\Phi_{b,\Om}(\mu_t)
 \le\Phi_{a,\Om}(\mu)+C_{\mathrm{trim}}
 \le\mathcal E_{a,\Om}(u)+C_{\mathrm{trim}}
 \le-L.
\]
Since
\[
 \mu_0=\Theta(u),\quad
 \mu_1=T^\theta_{p(u)}\Theta(u),
\]
the second homotopy joins $\mathscr S_b\Theta(u)$
to $G_b(F(u))$ in $X_{b,L}$.

Combining the two homotopies, we obtain
\[
 G_b\circ F\simeq\iota:
 X_{a,L+C}\to X_{b,L},
\]
where $\iota$ is the inclusion map.

\textbf{Step 4:} The second composition on the realization.

For $(\mu,\xi)\in\mathcal W_{a,L+C}$, set
\[
 u=\mathscr S_a\mu,\quad \nu=\Theta(u),\quad
 \zeta=p(u),\quad
 \xi^*=r_\tau(\xi),\quad \zeta^*=r_\tau(\zeta).
\]
By \eqref{eq:cech-duality-bounds}, we have
$u\in X_{a,L+C}$.

For each $j$, choose a smooth cutoff function $\psi_j$
such that
\[
 0\le\psi_j\le1,\quad
 \psi_j=1\quad\text{ on }\,C_j,\quad
 \supp\psi_j\Subset\operatorname{int}Z_j.
\]
Since $\mu=0$ on $Z_j$ whenever $\xi_j>0$,
Lemma~\ref{lem:cech-harmonic-buffer} implies
\[
 \nu\Big(\bigcup_{\xi_j>0}\supp\psi_j\Big)
 \le C_{\mathrm h}e^{-L-C}<\frac{\eta}{4}
\]
for sufficiently large $L$.
Moreover, by \eqref{eq:cech-union-small-mass},
\[
 \nu\Big(\bigcup_{\zeta_j>0}\supp\theta_j\Big)
 <\frac{\eta}{4}.
\]
For simplicity, set
\[
 w_\xi=\prod_{j=1}^M(1-b_j(\xi)\psi_j),\quad
 w_\zeta=\prod_{j=1}^M(1-b_j(\zeta)\theta_j).
\]
We construct a homotopy in $\mathcal W'_{b,L}$
in the following five steps.

(i) Keep $\mu$ fixed and move $\xi$ to $\xi^*$ along
the line segment
\[
 (\mu,(1-t)\xi+t\xi^*),\quad t\in[0,1].
\]
Since $\supp\xi^*\subset\supp\xi$, domain monotonicity
and $a<b$ imply that the path lies in
$\mathcal W'_{b,L}$.

\smallskip
\noindent
(ii) Keep $\xi^*$ fixed and define
\[
 \mu_t=(1-t)\mu+t\nu,\quad t\in[0,1].
\]
Consider the path
\[
 (T_{w_\xi}\mu_t,\xi^*).
\]
Since $\mu=0$ on $\supp\psi_j$ whenever $\xi_j>0$,
we have $T_{w_\xi}\mu=\mu$. Moreover,
\[
 \int_\Om(1-w_\xi)\mu_t\,\d x
 \le C_{\mathrm h}e^{-L-C}<\frac{\eta}{4}.
\]
If $\xi_j^*>0$, then $b_j(\xi)=1$, so
$T_{w_\xi}\mu_t=0$ on $C_j$.

By Lemma~\ref{lem:nonlinear-entropy-duality},
\[
 \Phi_{a,\Om}(\mu_t)
 \le\Phi_{a,\Om}(\mu)
 \le\Phi_{a,D_{\supp\xi}}(\mu)\le-L-C.
\]
Applying \eqref{eq:cech-normalized-trimming} and
Lemma~\ref{lem:cech-cut-comparison}, we obtain
\[
 \Phi_{b,D'_{\supp\xi^*}}(T_{w_\xi}\mu_t)
 \le\Phi_{b,\Om}(T_{w_\xi}\mu_t)+C_{\mathrm{buf}}
 \le\Phi_{a,\Om}(\mu_t)
      +C_{\mathrm{trim}}+C_{\mathrm{buf}}
 \le-L.
\]
Thus the path joins $(\mu,\xi^*)$ to
$(T_{w_\xi}\nu,\xi^*)$ in $\mathcal W'_{b,L}$.

(iii) Keep $\xi^*$ fixed and introduce the cutoff
functions $\theta_j$ continuously. More precisely,
set
\[
 w_{\zeta,t}
 =\prod_{j=1}^M(1-tb_j(\zeta)\theta_j),
 \quad t\in[0,1],
\]
and consider the path
\[
 (T_{w_\xi w_{\zeta,t}}\nu,\xi^*).
\]
The density remains zero on $C_j$ whenever
$\xi_j^*>0$.

(iv) Keep both cutoff products and move $\xi^*$
to $\zeta^*$ along the line segment
\[
 (T_{w_\xi w_\zeta}\nu,(1-t)\xi^*+t\zeta^*),
 \quad t\in[0,1].
\]
If $\xi_j^*>0$ or $\zeta_j^*>0$, then
$T_{w_\xi w_\zeta}\nu=0$ on $C_j$.

(v) Keep $\zeta^*$ fixed and remove the cutoff
functions $\psi_j$ continuously. Set
\[
 w_{\xi,t}
 =\prod_{j=1}^M(1-(1-t)b_j(\xi)\psi_j),
 \quad t\in[0,1],
\]
and consider the path
\[
 (T_{w_{\xi,t}w_\zeta}\nu,\zeta^*).
\]
Since $C_j\subset G_j$, the density remains zero
on $C_j$ whenever $\zeta_j^*>0$.

In each of (iii)--(v), the mass removed from $\nu$
is at most
\[
 C_{\mathrm h}e^{-L-C}+\frac{\eta}{4}
 <\frac{\eta}{2}.
\]
Moreover,
\[
 \Phi_{a,\Om}(\nu)
 \le\Phi_{a,\Om}(\mu)\le-L-C.
\]
Hence, for any  cutoff product $w$ used in (iii)--(v),
and for the corresponding set $A$ of positive simplex
coordinates, we have
\[
 \begin{aligned}
 \Phi_{b,D'_A}(T_w\nu)
 &\le\Phi_{b,\Om}(T_w\nu)+C_{\mathrm{buf}}\\
 &\le\Phi_{a,\Om}(\nu)+C_{\mathrm{trim}}+C_{\mathrm{buf}}
 \le-L.
 \end{aligned}
\]
Thus all five paths lie in $\mathcal W'_{b,L}$.

The final point of the homotopy is
\[
 (T_{w_\zeta}\nu,\zeta^*)
 =h_bFG_a(\mu,\xi).
\]
Since the cutoff functions are fixed, their coefficients
vary continuously, and the normalizing integrals are
uniformly bounded away from zero, the homotopy is
continuous in $L^p(\Om)\times\Delta$.

Combining (i)--(v), we obtain
\be\label{eq:cech-second-composition}
 h_bFG_a\simeq h_bc_{\mathcal W}
 \quad\text{ in }\, \mathcal W'_{b,L},
\ee
where $c_{\mathcal W}:\mathcal W_{a,L+C}\to
\mathcal W_{b,L}$ is the inclusion obtained by
increasing the parameter from $a$ to $b$ and the
sublevel value from $-L-C$ to $-L$.
All estimates hold uniformly for sufficiently large $L$.

\textbf{Step 5:} Homology comparison and Euler characteristics.

Let $c_X:X_{a,L+C}\to X_{b,L}$ be the inclusion map.
By Step~3, we have $G_b\circ F\simeq c_X$.
For $c_X$, factor the inclusion through $X_{b,L+C}$.
Both inclusions in this factorization induce isomorphisms by the
parameter and level comparisons in the proof of
Lemma~\ref{lem:cech-cut-comparison}. For $c_{\mathcal W}$, use
the same factorization on each $D_J$ in the density variables,
then apply the skeletal comparison of Step~2 to the compatible
inclusions. The same comparison applied to the zero extensions
at $(b,L)$ in Lemma~\ref{lem:cech-cut-comparison}(ii) gives the
isomorphism for $(h_b)_*$.

Since
\[
 (G_b)_*F_*=(c_X)_*,
\]
the map $F_*$ is injective. Moreover, by
\eqref{eq:cech-second-composition},
\[
 (h_b)_*F_*(G_a)_*
 =(h_b)_*(c_{\mathcal W})_*.
\]
Since $(h_b)_*$ is an isomorphism, we obtain
\[
 F_*(G_a)_*=(c_{\mathcal W})_*.
\]
Thus $F_*$ is surjective and hence is an isomorphism.
Since the parameter and level inclusions also induce
homology isomorphisms at $\rho$, we obtain
\[
 H_i(X_{\rho,L};\mathbb Q)
 \cong H_i(\mathcal W_{\rho,L};\mathbb Q),
 \quad i\ge0.
\]

By the additivity of the Euler characteristic and
\eqref{eq:cech-skeletal-homology}, we have
\[
 \chi(\mathcal W_{\rho,L})
 =\sum_{k=0}^{M-1}
 \chi(\mathcal W^{(k)},\mathcal W^{(k-1)})
 =\sum_{J\ne\varnothing}(-1)^{|J|-1}
 \chi(\Phi_{\rho,D_J}^{-L}).
\]
By Lemma~\ref{lem:nonlinear-entropy-duality}
and the above homology isomorphism,
\[
 \chi(J_{\rho,\Om}^{-\rho L})
 =\sum_{J\ne\varnothing}(-1)^{|J|-1}
 \chi(J_{\rho,D_J}^{-\rho L}).
\]
Since the Euler characteristics are independent of
sufficiently negative sublevel values, we obtain
\eqref{eq:annular-cech-Euler}.

Finally, by Proposition~\ref{prop:all-parameter-euler}
and the identity
\[
 \sum_{J\ne\varnothing}(-1)^{|J|-1}
 =\sum_{k=1}^M(-1)^{k-1}\binom Mk=1,
\]
we obtain
\[
 d_{\rho,\Om}
 =1-\chi(J_{\rho,\Om}^{-\rho L})=\sum_{J\ne\varnothing}(-1)^{|J|-1}
(1-\chi(J_{\rho,D_J}^{-\rho L})) 
 =\sum_{J\ne\varnothing}(-1)^{|J|-1}d_{\rho,D_J}.
\]
This proves \eqref{eq:annular-cech-degree}.
\end{proof}

\section{Proof of main results}\label{sec:5}
We now prove the results stated in Section~\ref{sec:intro}.

\begin{proof}[Proof of Theorem~\ref{thm:degree}]
We prove the result by induction on the number of
positive singularities.

\textbf{Step 1:} The case without positive singularities.

If $I_+=\varnothing$, then
\eqref{eq:odd-nonresonant-degree} reduces to
\eqref{eq:all-nonresonant-degree}.
By \cite[Theorem~1.1]{ATW26} and the finite identity
\[
 \sum_{\nu=0}^{\ell}(-1)^\nu
 \binom{\chi(\Om)-m}{\nu}
 =\binom{\ell-\chi(\Om)+m}{\ell},
\]
we obtain \eqref{eq:all-nonresonant-degree}.

Taking the difference of the degree formulas on two
consecutive nonresonant intervals gives the jumps
$a_i(s)$ in \eqref{eq:component-jumps}.
If several masses have the same value, we add
their contributions.

\textbf{Step 2:} Calculation after removing annuli.

Assume that both degree formulas hold whenever there
are fewer positive singularities. Fix $r>0$ with
$c_Nr\notin\Si_H$, and choose $p\in P_+$.

Take $M>r$ sufficiently large so that
Proposition~\ref{prop:annular-cech-excision} applies.
Choose $M$ concentric annuli around $p$ as in that
proposition.

Suppose that $j$ annuli are removed. The resulting
domain has $j+1$ connected components:
an inner ball containing $p$, an outer component
containing all the other singular points, and
$j-1$ annuli between them. The outer component has
Euler characteristic
\[
 \chi(\Om)-(-1)^N
\]
and contains $m-1$ singular points.

By the induction hypothesis, we know the degree
on the outer component. Lemma~\ref{lem:positive-local-factor}
gives the degree on the inner ball. By
Corollary~\ref{cor:degreepositions}, we may replace
the positive background weight on the inner ball
by a constant.

Each of the $j-1$ remaining annuli has no singular
points. By \cite[Theorem~1.1]{ATW26}, its degree
jumps at integer masses $\nu\ge0$ are
\[
 (-1)^\nu\binom{1-(-1)^N}{\nu}.
\]

We now combine the degrees of all connected components.
By Proposition~\ref{prop:disconnected-convolution},
the degree on the whole disconnected domain is the
sum of products of the component degree jumps.

First consider the integer mass. If the total integer
mass is $\ell$, the contributions from the outer
component and the $j-1$ annuli satisfy
\[
 \sum_{\substack{\ell_0,\ldots,\ell_{j-1}\in\N_0\\
                  \ell_0+\cdots+\ell_{j-1}=\ell}}
 \binom{\chi(\Om)-m+1-(-1)^N}{\ell_0}
 \prod_{\nu=1}^{j-1}
 \binom{1-(-1)^N}{\ell_\nu}=
 \binom{\chi(\Om)-m+j(1-(-1)^N)}{\ell}.
\]
This is the finite Vandermonde identity.
The common factor from the integer masses is $(-1)^\ell$.
The generalized binomial coefficients retain their own signs.

We next combine the contributions from positive
singularities. For even $N$, each positive singularity
contributes $-1$ at mass $\beta_i$.

For odd $N$, each positive singularity has two
possible nonzero jumps: $-2$ at mass $1$ and $1$
at mass $\min\{\beta_i,2\}$. Choose a subset
$J\subset I_+$ whose elements contribute at mass
$\min\{\beta_i,2\}$. Among the remaining positive
singularities, choosing $k$ contributions at mass
one gives the factor
\[
 (-2)^k\binom{|I_+|-|J|}{k}.
\]
Pascal's identity combines the choices on the inner
ball and the outer component and gives exactly the
coefficient in \eqref{eq:odd-nonresonant-degree}.

Consequently, for a domain obtained by removing
$j$ annuli, the degree has the same form as
\eqref{eq:all-nonresonant-degree} in even dimensions
and \eqref{eq:odd-nonresonant-degree} in odd dimensions.
The only change is that the binomial factor
\[
 \binom{\chi(\Om)-m}{\ell}
\]
is replaced by
\[
 \binom{\chi(\Om)-m+j(1-(-1)^N)}{\ell}.
\]

\textbf{Step 3:} Completion of the proof.

For any integer $0\le\ell<r<M$, the function
\[
 P(j)=\binom{\chi(\Om)-m+j(1-(-1)^N)}{\ell}
\]
is a polynomial in $j$ of degree at most $\ell<M$.
Hence its $M$th finite difference is zero, namely,
\[
 \sum_{j=0}^M(-1)^j\binom MjP(j)=0.
\]
Separating the term with $j=0$, we obtain
\[
 \sum_{j=1}^M(-1)^{j-1}\binom Mj
 \binom{\chi(\Om)-m+j(1-(-1)^N)}{\ell}
 =\binom{\chi(\Om)-m}{\ell}.
\]
Since all mass contributions in Step~2 are nonnegative,
every integer mass $\ell$ satisfies $\ell<r<M$.
Applying the above identity to each term in
\eqref{eq:annular-cech-degree}, we obtain
\eqref{eq:all-nonresonant-degree} and
\eqref{eq:odd-nonresonant-degree}.

Finally, when $\rho=0$, we have $K_{0,H}=0$ and hence
$d_0(\Om,H)=1$.
\end{proof}

The degree computation yields a criterion for solvability of \eqref{maineq}.

\begin{proof}[Proof of Corollary~\ref{cor:degree-existence}]
Since $d_\rho(\Om,H)\ne0$, \eqref{eq:solutionmap}
implies that $K_{\rho,H}$ has a fixed point in $B_R^E$.
By the definition of $K_{\rho,H}$, this fixed point
is a weak solution of \eqref{maineq}.
\end{proof}

The degree formula can be simplified by summing over
the integer mass index $\ell$.

\begin{proof}[Proof of Corollary~\ref{cor:even-full-degree}]
If $I_+=\varnothing$, then
\eqref{eq:odd-nonresonant-degree} reduces to
\eqref{eq:all-nonresonant-degree}.
For each subset $I$ with $\beta_I<r$, the nonresonance
condition \eqref{eq:sigma} implies that
$r-\beta_I$ is not an integer. Hence
\[
 \ell+\beta_I<r
 \quad\text{if and only if}\quad
 0\le\ell\le\lfloor r-\beta_I\rfloor
\]
for $\ell\in\N_0$.
Setting $K=\lfloor r-\beta_I\rfloor$ and using
the identity
\[
 \sum_{\ell=0}^{K}(-1)^\ell
 \binom{\chi(\Om)-m}{\ell}
 =(-1)^K\binom{\chi(\Om)-m-1}{K}
 =\binom{K-\chi(\Om)+m}{K},
\]
we obtain \eqref{eq:full-degree-formula}.
\end{proof}

In the first supercritical intervals, the degree follows directly
from small-hole comparison and a homotopy of the exponents.

\begin{proof}[Proof of Corollary~\ref{prop:first-positive-band}]
Let $J=\{j\in I_+:\beta_j>r\}$. In either case, $1<r<2$.
The hypotheses exclude $r$ from the resonance set
\eqref{eq:sigma}: any sum of two ordinary or positive masses
is at least two. In case (ii), adding a nonempty negative
subset to an ordinary or positive mass gives at least
$1+\beta_*>r$, and each sum of negative weights is less
than $r$.

Remove mutually disjoint small closed balls centered at the
points indexed by $J$, and denote the resulting domain by $D$.
If $J=\varnothing$, set $D=\Om$. Each removed point satisfies
the mass condition in Proposition~\ref{prop:puncture-compactness},
so Proposition~\ref{prop:puncture-degree} gives
$d_\rho(\Om,H)=d_\rho(D,H|_D)$. Since $N\ge2$, the balls can
be chosen so that $D$ is connected. Euler characteristic
additivity for the closures gives
\[
 \chi(D)=\chi(\Om)-(-1)^N|J|.
\]

For each $i\in I_+\setminus J$, decrease $\al_i$ continuously
to zero. Its normalized mass stays in $[1,\beta_i]\subset[1,r)$.
The same mass comparisons exclude resonance throughout
the deformation. Proposition~\ref{prop:parametercompact}
then preserves the degree on $D$. At the final weight only
the negative singularities remain. By
\cite[Theorem~1.1]{ATW26}, the degree is
\[
 \sum_{\substack{I\subset I_-,\,\ell\in\N_0\\
                   \ell+\beta_I<r}}
 (-1)^{|I|+\ell}\binom{\chi(D)-|I_-|}{\ell}.
\]
In case (i), only $\ell=0,1$ occur, and the sum is
$1-\chi(D)$. In case (ii), all negative subsets occur
for $\ell=0$, so their signs sum to zero. For $\ell=1$,
only the empty negative subset occurs. The sum is therefore
$|I_-|-\chi(D)$. Substitution of $\chi(D)$ yields \eqref{eq:first-positive-band} and \eqref{eq:clean-supercritical-degree}.
\end{proof}

The proof of existence uses direct minimization below
the first concentration mass and a topological argument
when the weighted barycenter space is noncontractible.

\begin{proof}[Proof of Theorem~\ref{thm:existence}]
(i) If $\rho<c_N\beta_*$, the weighted exponential
estimate following \eqref{eq:MT} implies that $J_\rho$
is coercive on $E$.
Since $Q$ is weakly lower semicontinuous and $M$
is weakly continuous on bounded subsets of $E$,
the direct method yields a global minimizer $u$ of
$J_\rho$. By the results in Section~\ref{sec:2},
$u$ has the required regularity and is positive
for $\rho>0$. When $\rho=0$, the minimizer is $u=0$.

(ii) Suppose that $J_\rho$ has no critical point.
Choose $L$ and $\lam$ as in Proposition~\ref{prop:homotopy}
and set $X=\Xr(\Ga)$. Consider the cone
\[
 CX=(X\times[0,1])/(X\times\{0\})
\]
and define
\[
 F:CX\to E,\quad
 F([\si,t])=t\Phi_\lam(\si).
\]
Since $F([\si,0])=0$ and $F([\si,1])=\Phi_\lam(\si)$,
the map $F$ is a continuous extension of $\Phi_\lam$
to $CX$.

By compactness of $CX$ and continuity of $J_\rho$,
there exists $b>-L$ such that
\[
 F(CX)\subset J_\rho^b.
\]
By Lemma~\ref{lem:descent}, there exists a retraction
\[
 R:J_\rho^b\to J_\rho^{-L}
\]
which fixes every point of $J_\rho^{-L}$.
Since $\Phi_\lam(X)\subset J_\rho^{-L}$, the composition
$R\circ F$ extends $\Phi_\lam$ to $CX$ within
$J_\rho^{-L}$. Hence $\Phi_\lam$ is homotopic to a
constant map in $J_\rho^{-L}$.

On the other hand, Proposition~\ref{prop:homotopy} yields
\[
 \Psi\circ\Phi_\lam\simeq\Id_X.
\]
Consequently, $\Id_X$ is homotopic to a constant map,
which contradicts the noncontractibility of $X$.
Therefore, $J_\rho$ has a critical point.
Regularity and positivity follow from the results
in Section~\ref{sec:2}.
\end{proof}

For the family in Corollary~\ref{cor:eulerzero},
we compute the degree and the homotopy type of the
weighted barycenter space by considering all subsets
of singular points with total weight less than $r$.

\begin{proof}[Proof of Corollary~\ref{cor:eulerzero}]
Since $kb<r<kb+a<1$, no ordinary mass or positive
singular mass can occur below $r$.

A subset containing a point of weight $kb$ has
total weight at least $kb+a$ if it contains another
point. Hence only the singleton of weight $kb$
has weight below $r$ in this case.

For subsets without points of weight $kb$, the
possible weights are $jb$ and $jb+a$, where $j$
is the number of points of weight $b$.
Since $a\le b$, we have
\[
 \begin{aligned}
 j\le k-1&:\quad jb+a\le(k-1)b+a\le kb,\\
 j=k&:\quad jb=kb,\quad jb+a=kb+a,\\
 j\ge k+1&:\quad jb\ge(k+1)b\ge kb+a.
 \end{aligned}
\]
Thus no subset has total weight in $(kb,kb+a)$.
By \eqref{eq:sigma}, every parameter in
\eqref{eq:homologyinterval} is nonresonant.

We next determine the homotopy type of
$\Xr(\overline\Om)$. Applying
Proposition~\ref{thm:threshold-wedge} with the vertex
of weight $a$ as the lightest vertex, we find that
the faces outside $\operatorname{st}(1)$ are exactly
\[
 \begin{aligned}
 &\binom nk\text{ singleton faces of weight }kb,\\
 &\binom nk\text{ faces consisting of $k$ vertices
 of weight }b.
 \end{aligned}
\]
The singleton faces correspond to spheres of dimension
zero, while the faces with $k$ vertices correspond
to spheres of dimension $k-1$. This proves
\eqref{eq:eulerzerofamily}.

Since $k$ is even, we obtain
\[
 \chi(\Xr(\overline\Om))-1
 =\binom nk+\binom nk(-1)^{k-1}=0.
\]
Hence $\chi(\Xr(\overline\Om))=1$.

Decrease each positive exponent continuously to zero.
Throughout this deformation every such normalized mass is
at least one, whereas $r<1$, so \eqref{eq:homologyinterval}
remains nonresonant. Proposition~\ref{prop:parametercompact}
therefore preserves the degree. The final weight has only
negative singularities, and \cite[Theorem~1.1]{ATW26} yields
\[
 d_\rho(\Om,H)
 =\sum_{\substack{I\subset I_-\\\beta_I<r}}
   (-1)^{|I|}
 =1-\chi(\Xr(\overline\Om))=0 .
\]
The second equality follows from \eqref{eq:Euleridentity}
and \eqref{eq:spinehomotopy}.

Although the degree is zero, the space
$\Xr(\overline\Om)$ is noncontractible by
\eqref{eq:eulerzerofamily}. Therefore,
Theorem~\ref{thm:existence}(ii) proves the existence
of a solution. Finally, \eqref{eq:splithomology}
yields the sublevel homology estimates.
\end{proof}

\medskip

		\noindent\textbf{Acknowledgements.}~
	M. Ahmedou is partially supported by the
	Deutsche Forschungsgemeinschaft (DFG) under grant AH~156/2-1. 
 Z. Hu is partially supported by the National Natural Science Foundation of China (Grant 12601383) and
 the China Postdoctoral Science Foundation (Grant 2026M793392). 
H. Wang is partially supported by the China Postdoctoral Science Foundation (2026M793418).
L. Zhang is partially supported by the Simons Foundation (SFI-MPS-TSM-00013752).

		\medskip

\noindent\textbf{Data availability statement.}~	 Data sharing is not applicable to this article as no datasets were generated or analysed during the current study.

\medskip

\noindent\textbf{Declarations.}~ \textbf{Conflict of interest.} ~ On behalf of all authors, the corresponding author states that there is no conflict of interest.

\medskip

\noindent\textbf{AI assistance statement.} ~The authors   used OpenAI models  for language editing, expository refinement, and preliminary consistency checks. All mathematical statements and proofs were independently verified by the authors, who take full responsibility for the content of the manuscript.


\begin{thebibliography}{99}


\bibitem{AB17}
M. Ahmedou and M. Ben Ayed, Theory of ``critical points at infinity'' and a resonant
singular Liouville-type equation,
\emph{Adv. Nonlinear Stud.} \textbf{17} (2017), no.~1, 139--166.

\bibitem{ABL17}
M. Ahmedou, M. Ben Ayed, and M. Lucia, On a resonant mean field type equation: A
``critical point at infinity'' approach,
\emph{Discrete Contin. Dyn. Syst.} \textbf{37} (2017), no.~4, 1789--1818.

\bibitem{AHW26}
M. Ahmedou, Z. Hu, and H. Wang, A degree-counting formula for a Keller--Segel equation
on a surface with boundary,
\emph{Commun. Contemp. Math.} (2026), Paper No.~2650039.

\bibitem{ATW26}
M. Ahmedou, Z. Tang, and H. Wang, Topological degree for singular $N$-mean field
equations, preprint (2026),
\href{https://arxiv.org/abs/2610.20281v1}{arXiv:2610.20281v1}.

\bibitem{B89}
A. Bahri,
\emph{Critical Points at Infinity in Some Variational Problems}, Pitman Research Notes
in Mathematics Series, vol.~182, Longman Scientific \& Technical, Harlow, 1989.

\bibitem{B15}
L. Battaglia, Existence and multiplicity result for the singular Toda system,
\emph{J. Math. Anal. Appl.} \textbf{424} (2015), 49--85.


\bibitem{B06}
A. Bj\"orn, Removable singularities for bounded $p$-harmonic and quasi(super)harmonic
functions on metric spaces,
\emph{Ann. Acad. Sci. Fenn. Math.} \textbf{31} (2006), 71--95.

\bibitem{CLMP92}
E. Caglioti, P.-L. Lions, C. Marchioro, and M. Pulvirenti, A special class of stationary
flows for two-dimensional Euler equations: A statistical mechanics description,
\emph{Comm. Math. Phys.} \textbf{143} (1992), no.~3, 501--525.

\bibitem{CLMP95}
E. Caglioti, P.-L. Lions, C. Marchioro, and M. Pulvirenti, A special class of stationary
flows for two-dimensional Euler equations: A statistical mechanics description. Part II,
\emph{Comm. Math. Phys.} \textbf{174} (1995), no.~2, 229--260.

\bibitem{CM12}
A. Carlotto and A. Malchiodi, Weighted barycentric sets and singular Liouville equations
on compact surfaces,
\emph{J. Funct. Anal.} \textbf{262} (2012), 409--450.

\bibitem{CGY93}
S.-Y. A. Chang, M. J. Gursky, and P. C. Yang, The scalar curvature equation on 2- and
3-spheres,
\emph{Calc. Var. Partial Differential Equations} \textbf{1} (1993), no.~2, 205--229.

\bibitem{CY87}
S.-Y. A. Chang and P. C. Yang, Prescribing Gaussian curvature on $S^2$,
\emph{Acta Math.} \textbf{159} (1987), nos.~3--4, 215--259.

\bibitem{CK94}
S. Chanillo and M. Kiessling, Rotational symmetry of solutions of some nonlinear
problems in statistical mechanics and in geometry,
\emph{Comm. Math. Phys.} \textbf{160} (1994), no.~2, 217--238.

\bibitem{CL02}
C.-C. Chen and C.-S. Lin, Sharp estimates for solutions of multi-bubbles in compact
Riemann surfaces,
\emph{Comm. Pure Appl. Math.} \textbf{55} (2002), no.~6, 728--771.

\bibitem{CL03}
C.-C. Chen and C.-S. Lin, Topological degree for a mean field equation on Riemann
surfaces,
\emph{Comm. Pure Appl. Math.} \textbf{56} (2003), no.~12, 1667--1727.

\bibitem{CL15}
C.-C. Chen and C.-S. Lin, Mean field equation of Liouville type with singular data:
Topological degree,
\emph{Comm. Pure Appl. Math.} \textbf{68} (2015), 887--947.

\bibitem{CLW04}
C.-C. Chen, C.-S. Lin, and G. Wang, Concentration phenomena of two-vortex solutions in a
Chern--Simons model,
\emph{Ann. Sc. Norm. Super. Pisa Cl. Sci.} (5) \textbf{3} (2004), no.~2, 367--397.

\bibitem{DPEFM15}
M. del Pino, P. Esposito, P. Figueroa, and M. Musso, Nontopological condensates for the
self-dual Chern--Simons--Higgs model,
\emph{Comm. Pure Appl. Math.} \textbf{68} (2015), no.~7, 1191--1283.

\bibitem{DPEM10}
M. del Pino, P. Esposito, and M. Musso, Two-dimensional Euler flows with concentrated
vorticities,
\emph{Trans. Amer. Math. Soc.} \textbf{362} (2010), no.~12, 6381--6395.

\bibitem{DJLW97}
W. Ding, J. Jost, J. Li, and G. Wang, The differential equation $\Delta u=8\pi-8\pi
he^u$ on a compact Riemann surface,
\emph{Asian J. Math.} \textbf{1} (1997), 230--248.

\bibitem{DJLW98}
W. Ding, J. Jost, J. Li, and G. Wang, An analysis of the two-vortex case in the
Chern--Simons Higgs model,
\emph{Calc. Var. Partial Differential Equations} \textbf{7} (1998), no.~1, 87--97.

\bibitem{DJLW99}
W. Ding, J. Jost, J. Li, and G. Wang, Existence results for mean field equations,
\emph{Ann. Inst. H. Poincar\'e Anal. Non Lin\'eaire} \textbf{16} (1999), 653--666.

\bibitem{D08}
Z. Djadli, Existence result for the mean field problem on Riemann surfaces of all
genuses,
\emph{Commun. Contemp. Math.} \textbf{10} (2008), 205--220.

\bibitem{DM08}
Z. Djadli and A. Malchiodi, Existence of conformal metrics with constant $Q$-curvature,
\emph{Ann. of Math.} (2) \textbf{168} (2008), 813--858.

\bibitem{DHM00}
G. Dolzmann, N. Hungerb\"uhler, and S. M\"uller, Uniqueness and maximal regularity for
nonlinear elliptic systems of $n$-Laplace type with measure valued right hand side,
\emph{J. Reine Angew. Math.} \textbf{520} (2000), 1--35.

\bibitem{E21}
P. Esposito, Isolated singularities for the $n$-Liouville equation,
\emph{Calc. Var. Partial Differential Equations} \textbf{60} (2021), Paper No.~137.

\bibitem{EF14}
P. Esposito and P. Figueroa, Singular mean field equations on compact Riemann surfaces,
\emph{Nonlinear Anal.} \textbf{111} (2014), 33--65.

\bibitem{EM15}
P. Esposito and F. Morlando, On a quasilinear mean field equation with an exponential
nonlinearity,
\emph{J. Math. Pures Appl.} \textbf{104} (2015), 354--382. Numbering cited in the text
follows \href{https://arxiv.org/abs/1410.6779}{arXiv:1410.6779}.

\bibitem{JKM15}
A. Jevnikar, S. Kallel, and A. Malchiodi, A topological join construction and the Toda
system on compact surfaces of arbitrary genus,
\emph{Anal. PDE} \textbf{8} (2015), 1963--2027.

\bibitem{JWY18}
A. Jevnikar, J. Wei, and W. Yang, On the topological degree of the mean field equation
with two parameters,
\emph{Indiana Univ. Math. J.} \textbf{67} (2018), no.~1, 29--88.

\bibitem{K19}
S. Kallel, Formal barycenter spaces with weights: the Euler characteristic,
\emph{Topol. Methods Nonlinear Anal.} \textbf{53} (2019), 801--823.

\bibitem{KW74}
J. L. Kazdan and F. W. Warner, Curvature functions for compact 2-manifolds,
\emph{Ann. of Math.} (2) \textbf{99} (1974), no.~1, 14--47.

\bibitem{KV86}
S. Kichenassamy and L. V\'eron, Singular solutions of the $p$-Laplace equation,
\emph{Math. Ann.} \textbf{275} (1986), 599--615.

\bibitem{K93}
M. K.-H. Kiessling, Statistical mechanics of classical particles with logarithmic
interactions,
\emph{Comm. Pure Appl. Math.} \textbf{46} (1993), no.~1, 27--56.


\bibitem{LSY23}
J. Li, L. Sun, and Y. Yang, The boundary value problem for the mean field equation on a
compact Riemann surface,
\emph{Sci. China Math.} \textbf{66} (2023), 115--142.

\bibitem{LZ19}
J. Li and C. Zhu, The convergence of the mean field type flow at a critical case,
\emph{Calc. Var. Partial Differential Equations} \textbf{58} (2019), Paper No.~60.

\bibitem{L99}
Y. Y. Li, Harnack type inequality: The method of moving planes,
\emph{Comm. Math. Phys.} \textbf{200} (1999), no.~2, 421--444.

\bibitem{LW10}
C.-S. Lin and C.-L. Wang, Elliptic functions, Green functions and the mean field
equations on tori,
\emph{Ann. of Math.} (2) \textbf{172} (2010), no.~2, 911--954.

\bibitem{LY13}
C.-S. Lin and S. Yan, Existence of bubbling solutions for Chern--Simons model on a
torus,
\emph{Arch. Ration. Mech. Anal.} \textbf{207} (2013), no.~2, 353--392.

\bibitem{L07}
M. Lucia, A deformation lemma with an application to a mean field equation,
\emph{Topol. Methods Nonlinear Anal.} \textbf{30} (2007), 113--138.

\bibitem{M08}
A. Malchiodi, Morse theory and a scalar field equation on compact surfaces,
\emph{Adv. Differential Equations} \textbf{13} (2008), 1109--1129.

\bibitem{M26}
A. Malchiodi, On a degree-counting formula for elliptic equations of exponential type,
\emph{Discrete Contin. Dyn. Syst.} \textbf{55} (2026), 288--305.

\bibitem{M71}
J. Moser, A sharp form of an inequality by N. Trudinger,
\emph{Indiana Univ. Math. J.} \textbf{20} (1970/71), 1077--1092.

\bibitem{NT99}
M. Nolasco and G. Tarantello, Double vortex condensates in the Chern--Simons--Higgs
theory,
\emph{Calc. Var. Partial Differential Equations} \textbf{9} (1999), no.~1, 31--94.

\bibitem{PW13}
J. Pakianathan and T. Winfree, Threshold complexes and connections to number theory,
\emph{Turkish J. Math.} \textbf{37} (2013), 511--539.

\bibitem{S67}
J. T. Stuart, On finite amplitude oscillations in laminar mixing layers,
\emph{J. Fluid Mech.} \textbf{29} (1967), no.~3, 417--440.

\bibitem{SWY21}
L. Sun, Y. Wang, and Y. Yang, Existence results for a generalized mean field equation on
a closed Riemann surface, preprint (2021),
\href{https://arxiv.org/abs/2101.03859}{arXiv:2101.03859}.

\bibitem{SZ21}
L. Sun and J. Zhu, Global existence and convergence of a flow to Kazdan--Warner equation
with non-negative prescribed function,
\emph{Calc. Var. Partial Differential Equations} \textbf{60} (2021), Paper No.~42.

\bibitem{SZ24}
L. Sun and J. Zhu, An existence result for the Kazdan--Warner equation with a
sign-changing prescribed function,
\emph{Calc. Var. Partial Differential Equations} \textbf{63} (2024), Paper No.~52.

\bibitem{YY24}
P. Yu and Y. Yang, Existence of the heat flow with sign-changing prescribed function,
\emph{J. Math. Anal. Appl.} \textbf{535} (2024), Paper No.~128118.

\bibitem{Z23}
X. Zhu, A singular Kazdan--Warner problem on a compact Riemann surface,
\emph{Calc. Var. Partial Differential Equations} \textbf{62} (2023), Paper No.~164.

\end{thebibliography}
\end{document}